\documentclass[11pt]{article}

\usepackage[margin=1in]{geometry}
\usepackage{graphicx}
\usepackage{amsfonts}
\usepackage{amsmath}
\usepackage{amssymb}
\usepackage{amsthm}
\usepackage{enumitem}
\usepackage{mathrsfs}
\usepackage{bm}
\usepackage{hyperref}
\usepackage{tikz}
\usetikzlibrary{graphs, positioning}
\usepackage{array}
\usepackage{booktabs}
\usepackage{adjustbox}
\usepackage{subfig}
\usepackage{float}
\usepackage{cleveref}

\DeclareMathOperator{\diag}{diag}

\DeclareMathOperator{\tr}{tr}

\def\Acal{{\mathcal{A}}}
\def\Bcal{{\mathcal{B}}}
\def\Ccal{{\mathcal{C}}}
\def\Dcal{{\mathcal{D}}}
\def\Ecal{{\mathcal{E}}}

\def\Ical{{\mathcal{I}}}

\def\Kcal{{\mathcal{K}}}

\def\Ocal{{\mathcal{O}}}

\def\Tcal{{\mathcal{T}}}

\def\Rbb{{\mathbb{R}}}

\def\onebm{{\bm{1}}}

\def\b0{{\mathbf{0}}}

\theoremstyle{plain}
\newtheorem{theorem}{Theorem}
\newtheorem{proposition}[theorem]{Proposition}

\newtheorem{definition}[theorem]{Definition}

\title{On the number of small edge-weighted subgraphs}

\author{
    Feng Yu\thanks{Department of Mathematical Sciences, University of Texas at El Paso, El Paso, TX, USA (\texttt{fyu@utep.edu})}
    \and
    Mingao Yuan\thanks{Department of Mathematical Sciences, University of Texas at El Paso, El Paso, TX, USA (\texttt{myuan2@utep.edu})}
}

\date{}

\begin{document}

\maketitle

\begin{abstract}
Subgraph counting is a fundamental task that underpins several network analysis methodologies, including community detection and graph two-sample tests. Counting subgraphs is a computationally intensive problem. Substantial research has focused on developing efficient algorithms and strategies to make it feasible for larger unweighted graphs.  Implementing those algorithms can be a significant hurdle for data professionals or researchers with limited expertise in algorithmic principles and programming. Furthermore, many real-world networks are weighted. Computing the number of weighted subgraphs in weighted networks presents a computational challenge, as no efficient algorithm exists for the worst-case scenario. In this paper, we derive explicit formulas for counting small edge-weighted subgraphs using the weighted adjacency matrix. These formulas are applicable to unweighted networks, offering a simple and highly practical analytical tool for researchers across various scientific domains. In addition, we introduce a generalized methodology for calculating arbitrary weighted subgraphs.
\end{abstract}



\section{Introduction}

Networks or graphs offer a highly versatile and effective framework for representing a wide range of real-world phenomena. They represent entities in a system as nodes and their interactions as edges. Network analysis has thus found broad application in diverse scientific domains, including social and biological sciences. Consequently, analyzing and extracting information from networks has become a significant and independent field of study.

Subgraph counting is a fundamental task that lies at the core of several network analysis methodologies. For example, the degree variance test statistic is a function of the number of 2-paths, used for testing a planted community in dense random networks \cite{AV14}. Cycle-based and path-based statistics have been successfully employed to develop powerful tests for community structures in networks \cite{GL17,JKL21,YYS22,YLFS22,YS22,fan2022alma}. The frequency of short paths provides a measure for evaluating the uniformity of network nodes and the dimension of the network in the embedded space \cite{YR25,YY25,Y25}. Cliques provide a basis for developing higher-order clustering coefficients that quantify the extent of clustering beyond simple triangles \cite{yin2018higher}. Moreover, subgraph counts provide a flexible method for statistically comparing network samples, allowing us to test whether they originate from a predefined distribution, or share the same model as a separate network sample \cite{GGC17,maugis2020testing}.

Subgraph counting is a computationally expensive task. Subgraph occurrences can be surprisingly large even in smaller networks. Due to its practical importance and inherent difficulty, subgraph counting has received considerable attention from the research community \cite{ribeiro2021survey,gonen2011counting,ABG18,FFF15,pinar2017escape}. A range of polynomial-time algorithms have been developed for counting the frequency of network motifs. For instance,  sublinear-time algorithms have been developed to estimate the count of stars \cite{gonen2011counting,ABG18}. Deterministic and randomized algorithms for k-clique enumeration have been introduced and evaluated in \cite{FFF15}.

The practical application of the existing algorithms presents a formidable hurdle for those lacking a strong background in algorithmic theory and programming. Moreover, these algorithms do not apply to weighted networks. In practice, many real-world networks are weighted, with associated weights that record the strength or intensity of the connections between two nodes \cite{N04}. Analyzing the frequency of weighted subgraphs plays a vital role in understanding their structure. Computing the number of weighted subgraphs in weighted networks is computationally prohibitive due to the lack of an efficient algorithm for general instances \cite{vassilevska2009finding}.  A direct computation through nested for-loops
is computationally intense. It is therefore essential to derive formulas for the number of subgraphs based on the adjacency matrix. In the context of unweighted networks, the literature provides several established formulas, derived from the adjacency matrix, for particular subgraphs such as cycles and paths \cite{HM71,BR23}. However, to the best of our knowledge,  analogous formulas do not exist for other subgraphs or for weighted subgraphs in weighted networks.

In this article, we introduce a universal method for exactly computing arbitrary weighted subgraphs. Moreover, we present explicit formulas for the number of short edge-weighted subgraphs in weighted networks using the adjacency matrix. Specifically, we list the number of weighted subgraphs with up to five nodes.   These formulas offer an accessible and highly effective analytical tool for researchers across diverse scientific disciplines. The formulas are applicable to unweighted graphs and are significantly simpler in that context. Numerical simulations show that the proposed formulas achieve speedups from $1{,}800\times$ to $600{,}000\times$ for counting 5-node weighted subgraphs in a network of size $50$.

\subsection*{Notations}
Let $A \in \mathbb{R}^{n \times n}$ be a symmetric matrix with $\diag(A) = \bm{0}$, and let $\bm{1}_n = [1, \ldots, 1]^\top \in \mathbb{R}^n$. When the dimension is clear from context, we simply write $\bm{1}$. The symbol $\odot$ denotes the Hadamard product (element-wise product) of matrices (or tensors), and $A^{\circ k}$ represents the Hadamard power, with $[A^{\circ k}]_{ij} = A_{ij}^k$. Given positive integers $1 \le i_1, \ldots, i_m \le n$, the notation $\sum_{i_1 \neq i_2 \neq \cdots \neq i_m}$ indicates summation over distinct indices $i_1, i_2, \ldots, i_m$, whereas $\sum_{i_1, i_2, \ldots, i_m}$ denotes summation over all indices, which may coincide.

If $\Acal\in\Rbb^{p\times m\times q_1\times \ldots\times q_s},\Bcal\in\Rbb^{m\times q\times q_1\times \ldots\times q_s}$ are two tensors, their mode-(1,2) product denoted by $\Acal\times_{1,2}\Bcal$ is a tensor in $\Rbb^{p\times q\times q_1\times \ldots\times q_s}$ with elements
\begin{align*}
    [\Acal\times_{1,2}\Bcal]_{iji_1\ldots i_s} = \sum_{k=1}^m\Acal_{{ik}i_1\ldots i_s}\Bcal_{kji_1\ldots i_s}.
\end{align*}

For a graph $G=(V,E)$, the node set is $V=\{i_1,\ldots,i_m\}$ and $E$ is the edge set. If nodes $i_a$ and $i_b$ are connected, we write $(i_a,i_b)\in E$. We also use the compact notation $i_a-i_b$ to denote an edge. For example, if $E=\{i_1-i_2-i_3,i_4-i_5\}$, then this represents the edge set $E = \{(i_1,i_2),\,(i_2,i_3),\,(i_1,i_3),\,(i_4,i_5)\}$. If any edge occurs more than once, i.e., if $G$ is a multigraph, we indicate the multiplicity of the edge. For instance, $m(i_a,i_b)=2$ denotes that $i_a$ and $i_b$ are connected twice.

\section{Main Results}
In this section, we provide a general way to compute the labeled number of weighted subgraphs in weighted networks. We consider the summation of an index function over distinct indices and present general results in \Cref{sec:general_results}, along with explicit expressions for small subgraphs in \Cref{sec:small_graphs}. When the index function is associated with a graph, we propose a fast method to evaluate this sum in \Cref{sec:faster_contraction}, and the algorithm is summarized in \Cref{sec:summary}.

\subsection{General Results}\label{sec:general_results}
Given a symmetric matrix $A \in \mathbb{R}^{n \times n}$ with $\diag(A) = \bm{0}$, we are interested in evaluating the following quantity
\begin{align}\label{eq:target_sum}
     L = \sum_{i_1\neq i_2\neq \ldots\neq i_m}A_{i_{k_1}i_{k_2}}A_{i_{k_3}i_{k_4}}\ldots A_{i_{k_{2p-1}}i_{k_{2p}}}, \quad i_{k_\ell}\in\{i_1,\ldots,i_m\}, \forall 1\leq\ell\leq 2p.
\end{align}
The quantity $L$ is the number of labeled weighted subgraphs.
The number of weighted subgraphs can be expressed as a function of $L$. For example, $\sum_{i\neq j}A_{ij}$ is 2 times the total number of weighted edges, and  $\sum_{i\neq j\neq k}A_{ij}A_{jk}A_{ki}$ is 6 times the number of weighted triangles.

Note that $i_{k_\ell}\in\{i_1,\ldots,i_m\}, \forall 1\leq\ell\leq 2p$, the product $A_{i_{k_1}i_{k_2}}A_{i_{k_3}i_{k_4}}\ldots A_{i_{k_{2p-1}}i_{k_{2p}}}$ is an index function of $i_1,\ldots,i_m$, so we denote it by $F(i_1,\ldots,i_m)=A_{i_{k_1}i_{k_2}}A_{i_{k_3}i_{k_4}}\ldots A_{i_{k_{2p-1}}i_{k_{2p}}}$. Let $m\geq 2$ be fixed, and denote $S=\{i_1,\ldots,i_m\}$ by a set of $m$ indices. Let $G=(S,E)$ be a connected graph where $E=\{(i_{k_1},i_{k_2}),\ldots,(i_{k_{2p-1}},i_{k_{2p}})\}$ be the set of edges. Our results are constructed based on the concept of partition, which is provided in Definition~\ref{def:partition}.

\begin{definition}[Partition]\label{def:partition}
    Fix $k \leq m$. We call $\pi = \{B_1, B_2,\ldots, B_k\}$ a $k$-partition of $S$ if $B_1, B_2,\ldots, B_k$ are pairwise disjoint non-empty subsets of $S$ satisfying $\bigcup_{j=1}^kB_j=S$. We call each $B_j$ a block.
    For two partitions $\pi$ and $\sigma$, we say $\pi\preceq\sigma$ if each block of $\pi$ is a subset of a block of $\sigma$, and we say $\pi\prec\sigma$ if $\pi\preceq\sigma$ but $\pi \neq \sigma$. 
    Let $\Pi(S)$ denote the set of all partitions of $S$. Equipped with the relation $\preceq$, the set $\Pi(S)$ forms a \textit{partially ordered set (poset)}, denoted $(\Pi(S), \preceq)$.

\end{definition}

Let $\pi=\{B_1,\ldots,B_k\}$ be a partition of $S=\{i_1,\ldots,i_m\}$.
For each block $B\in\pi$, we assign a label $i_B\in S$.
Then, for each index $i_s\in S$, there exist a unique block $B\in\pi$ such that $i_s\in B$. This defines a surjective mapping $f_\pi:S\rightarrow S$ given by $f_\pi(i_s)=i_B$, which maps each element to the label of the block containing it.
We define two sums associated with the partition $\pi$:
\begin{align}
    S_{\pi}(F) & = \sum_{\substack{B\in\pi \\ i_B \text{(distict)}} } F(f_\pi(i_1),\ldots,f_\pi(i_m)), \label{eq:DBS}\\
    M_{\pi}(F) & = \sum_{B\in\pi } F(f_\pi(i_1),\ldots,f_\pi(i_m)),\label{eq:GBS}
\end{align}
where $S_\pi$ is referred to as the distinct-block sum (DBS) and $F_\pi$ as the general block sum (GBS). 
Both of $S_\pi$ and $M_\pi$ are the sum of $F$ over all tuples that are constant on each block of $\pi$. The difference is that, for $M_\pi$, the blocks can take any values, including repeated values across different blocks, whereas for $S_\pi$, the blocks must take distinct values. 

\noindent\textbf{Relation between DBS and GBS.} Let $\Pi(S)$ denote the set of all partitions of $S$. 
If $\pi\preceq\sigma$ for any two partitions $\pi$ and $\sigma$ in $\Pi(S)$, then any block in $\sigma$ is a union of the blocks of $\pi$. Consequently, every assignment counted in $S_\sigma$ is also included in $M_\pi$. Summing over all partitions $\sigma \succeq \pi$ yields the relation
\begin{align}
    M_\pi = \sum_{\sigma\succeq\pi}S_\sigma.
\end{align}
Conversely, $S_\pi$ can be expressed in terms of $M_\pi$ through the M\"obius inversion formula by \Cref{prop:mobius_inversion}.

\begin{proposition}[{\cite[Proposition 2]{rota1964foundations}}]\label{prop:mobius_inversion}
Suppose that $M_\pi$ and $S_\sigma$ are defined in \eqref{eq:DBS} and \eqref{eq:GBS}, respectively. Then it follows that
\begin{align*}
    S_\pi = \sum_{\sigma \succeq \pi} \mu(\pi, \sigma) M_\sigma,
\end{align*}
where $\mu(\pi, \sigma)$ is the M\"obius function on the lattice of partitions $\Pi(S)$.
For any $\pi, \sigma \in \Pi(S)$ with $\pi \prec \sigma$, it satisfies $\mu(\pi, \pi) = 1$, and is defined recursively by
$\mu(\pi, \sigma) = - \sum_{\pi \preceq \tau \prec \sigma} \mu(\pi, \tau)$.
\end{proposition}

Evaluating the summation in \eqref{eq:target_sum} can be interpreted as an inclusion–exclusion problem.
Note that a special partition of $S$ is $\hat{0} = \{\{i_1\}, \{i_2\}, \ldots, \{i_m\}\}$. The DBS corresponding to this partition, $S_{\hat{0}}$, equals the desired term $\sum_{i_1 \neq i_2 \neq \ldots \neq i_m} F(i_1, \ldots, i_m)$.
The motivation for expressing $S_\pi$ in terms of $M_\pi$ is that the GBS $M_\pi$ is sometimes more computationally efficient to evaluate. According to \Cref{prop:mobius_inversion}, since $\hat{0}$ is the finest partition, we have that $S_{\hat{0}}= \sum_{\sigma\in\Pi(S)} \mu(\hat{0}, \sigma) M_\sigma$ and we need the M\"obius function values $\mu(\hat{0}, \sigma)$ and $M_\sigma$ for all blocks in $\Pi(S)$ to express $\sum_{i_1 \neq i_2 \neq \ldots \neq i_m} F(i_1, \ldots, i_m)$, which is summarized in \Cref{thm:main}. The coefficients $\mu(\hat{0}, \sigma)$ are provided in \Cref{thm:mobius_values}. 

\begin{theorem}\label{thm:main}
    Let $S=\{i_1,\ldots,i_m\}$ and $\Pi(S)$ be the set of all possible partitions of $S$. We have the following 
    \begin{align*}
        \sum_{i_1 \neq i_2 \neq \ldots \neq i_m} F(i_1, \ldots, i_m) = \sum_{\sigma\in\Pi(S)}\mu(\hat{0},\sigma)M_\sigma^{(m)}.
    \end{align*}
\end{theorem}

\begin{theorem}\label{thm:mobius_values}
Let $\mu(\pi,\sigma)$ denote the M\"obius function on the partition lattice $\Pi(S)$.
For the finest partition $\hat{0}=\{\{i_1\},\ldots,\{i_m\}\}$ and any $\sigma\in\Pi(S)$ with blocks $B_1,\ldots,B_k$, we have
\begin{align}
    \mu(\hat{0},\sigma)
    = (-1)^{m-k}\prod_{B\in\sigma} (|B|-1)! .
\end{align}    
\end{theorem}

\begin{proof}

For any partition $\sigma$ with blocks $B_1,\ldots,B_k$, the interval $[\hat{0}, \sigma]$ consists of all partitions finer than $\sigma$. Refinements within distinct blocks are independent, so this interval decomposes as a Cartesian product of partition lattices on the blocks:
\begin{align}
    [\hat{0},\sigma] \cong \Pi(B_1)\times \Pi(B_2)\times\ldots\Pi(B_k),
\end{align}
Let $1_B$ denote the coarsest partition of $B$, i.e., all elements of $B$ are in a single block.
By the general product property of the M\"obius function for posets, it follows that
\begin{align}\label{eq:mobius_fun_prod}
    \mu(\hat{0},\tau)
    = \prod_{B\in\tau}\mu_{\Pi(B)}(\hat{0},1_B).
\end{align}
Because each factor $\Pi(B)$ depends only on the size $|B|$ of $B$, we write $1_{|B|}:=1_B$ (up to relabeling) and obtain the compact form of \eqref{eq:mobius_fun_prod} as
\[
\mu(\hat{0},\tau)=\prod_{B\in\tau}\mu(\hat{0},1_{|B|}).
\]
Thus it suffices to compute $\mu(\hat{0},1_r)$ for a single block of size $r$. For a fixed block $B$ with $|B|=r$, we have the base case $\mu(\hat{0},1_1)=1$ and the recurrence $\mu(\hat{0},1_s)=-(s-1)\mu(\hat{0},1_{s-1})$ for all $2\leq s<r$, which follows by grouping partitions according to the size of the block containing a fixed element. (details see \cite[Example 3.10.4]{stanley2011enumerative}). Thus, one obtains by induction that $\mu(\hat{0},1_s)=(-1)^{s-1}(s-1)!$. Therefore,
\begin{align}
    \mu(\hat{0},\tau)=\prod_{B\in\tau}\mu(\hat{0},1_{|B|}) = \prod_{B\in\tau}(-1)^{|B|-1}(|B|-1)! = (-1)^{m-k}\prod_{B\in\tau}(|B|-1)!,
\end{align}
where $m$ is the total number of elements of $S$ and $k$ the number of blocks of $\sigma$. This completes the proof.
\end{proof}

\subsection{Results for small subgraphs}\label{sec:small_graphs}
The next step is to evaluate $M_\sigma$ which needs to find all $\sigma\in\Pi(S)$. 
Let $\sigma=\{B_1,\ldots,B_k\}$ and let $\bm{n}_\sigma=(n_1,n_2,\ldots,n_k)$ denote the sorted sequence of block sizes, where $(n_1,n_2,\ldots,n_k)$ is a nonincreasing arrangement of $(|B_1|,|B_2|,\ldots,|B_k|)$.
We classify $M_\sigma$ by its $\bm{n}_\sigma$-type, represented by the subscript $n_1+n_2+\ldots+n_k$. This notation is convenient because, according to \Cref{thm:mobius_values}, all $M_\sigma$ with the same $\bm{n}_\sigma$ share the same coefficient. A general procedure for evaluating $M_\sigma$ of $\bm{n}_\sigma$-type is given in \Cref{prop:M_sig}.

\begin{proposition}[$M_\sigma$ of $\bm{n}_\sigma$-type]\label{prop:M_sig}
Let $m\geq 2$ is fixed and there exists $\bm{n}_\sigma=(n_1,n_2,\ldots,n_k)$ for some $1\leq k\leq m$ such that $n_1\geq n_2\ldots\geq n_k$ with $\sum_{i=1}^kn_k=m$. Then the GBS $M_\sigma^{(m)}$ can be evaluated as follows:
\begin{align}\label{eq:M_sig}
    M_\sigma^{(m)} & = \sum_{i_1,\ldots,i_k} \sum_{\substack{\bigcup^k_{q=1} S_q=[m]\\|S_q|=n_q,\forall 1\leq q\leq k}}F(x_1,\ldots,x_m):=\sum_{i_1,\ldots,i_k}F_{n_1+\ldots+n_k}^{(m)}(i_1,\ldots,i_k),
\end{align}
where $x_t=i_q$ if $t\in S_q$. We call $F_{n_1+\ldots+n_k}^{(m)}(i_1,\ldots,i_k)$ as the k-fold contraction of $F$ with multiplicities $\bm{n}_\sigma=(n_1,n_2,\ldots,n_k)$. 

\end{proposition}

\begin{proof}
    If $M_\sigma$ is of $\bm{n}_\sigma=(n_1,n_2,\ldots,n_k)$-type, then there exists a partition $\sigma=\{B_1,\ldots,B_k\}$ with $k$ blocks and its nonincreasing arrangement of $(|B_1|,|B_2|,\ldots,|B_k|)$ is $(n_1,n_2,\ldots,n_k)$. Since all indices within a block $B_q$ are identical, we assign them the common label $i_q$. Hence, summing $F$ over all $\sigma$ of the same type reduces to summing over all ways to choose which positions share each label, that is,
    \begin{align*}
        \sum_{\substack{(S_1,\ldots,S_k):\bigcup_{q=1}^kS_q=[m]\\ |S_q|=n_q,\forall 1\leq q\leq k}}F(x_1,\ldots,x_m), \qquad x_t = i_q \text{ if } t\in S_q.
    \end{align*}
\end{proof}

\noindent\textbf{Remark:} 
We observe that the function $F$ itself is not symmetric with respect to any pair of variables. However, when we sum over all indices, the resulting quantity becomes fully symmetric. For instance, if $F = F(i,j,k)$, in general $F(i,j,k) \neq F(i,k,j)$, yet $\sum_{i,j,k} F(i,j,k) = \sum_{i,j,k} F(i,k,j)$, since the summation is invariant under any permutation of the indices. Therefore, such functions can be viewed as belonging to the same equivalence class under index permutation, i.e., $F(i,j,k) \sim F(i,k,j)$ if $\sum_{i,j,k} F(i,j,k) = \sum_{i,j,k} F(i,k,j)$. This notion of equivalence also extends to the contraction function, since it is defined by selecting a subset of indices from $F$ to form a new function with fewer variables, while the overall summation remains invariant under the same permutations.

\Cref{prop:M_sig} provides a simple and systematic way to evaluate $M_\sigma$. For any fixed $m \ge 2$, we can enumerate all possible $\bm{n}_\sigma$-types of $M_\sigma$ and, using \eqref{eq:M_sig}, construct the corresponding functions $F_{n_1+\cdots+n_k}^{(m)}$ whose sums over $i_1,\ldots,i_k$ are typically more convenient to compute and often admit computationally efficient equivalent forms. Consequently, by applying \Cref{thm:main}, the term $\sum_{i_1 \neq i_2 \neq \cdots \neq i_m} F(i_1, \ldots, i_m)$ can be rewritten in a computationally efficient equivalent form.

We provide the complete formulas of $\sum_{i_1 \neq i_2 \neq \cdots \neq i_m} F(i_1, \ldots, i_m)$ up to $m=5$ in \Cref{prop:general_F}.

\begin{proposition}\label{prop:general_F}
The following identities hold
\begin{enumerate}[label=(\alph*)]
    \item $\displaystyle\sum_{i\not=j}F(i,j) = \sum_{i,j}F(i,j)-\sum_{i}F(i,i)$.
    \item $\displaystyle\sum_{i\not=j\not=k}F(i,j,k) = \sum_{i,j,k}F(i,j,k)-\sum_{i,j}F_{2+1}^{(3)}(i,j)+2\sum_{i}F_{3}^{(3)}(i)$,
    where $F_{2+1}^{(3)}(i,j) = F(i,i,j)+ F(i,j,i)+F(j,i,i)$ and $F_{3}^{(3)}(i) = F(i,i,i)$.
    \item $\displaystyle\sum_{i\not=j\not=k\not=l}F(i,j,k,l) = \sum_{i,j,k,l}F(i,j,k,l)-\sum_{i,j,k}F_{2+1+1}^{(4)}(i,j,k)+\sum_{i,j}F_{2+2}^{(4)}(i,j)+2\sum_{i,j}F_{3+1}^{(4)}(i,j)
    -6\sum_iF_{4}^{(4)}(i),$
    where $F_{2+1+1}^{(4)}(i,j,k) = F(i,i,j,k)+F(i,j,i,k) +F(i,j,k,i)+F(j,i,i,k)+F(j,i,k,i)+F(j,k,i,i)$, $F_{2+2}^{(4)}(i,j) = F(i,i,j,j) +F(i,j,i,j)+F(i,j,j,i)$, $F_{3+1}^{(4)}(i,j) = $ 
    
    $F(i,i,i,j)+F(i,i,j,i)+F(i,j,i,i)+F(j,i,i,i)$ and $F_{4}^{(4)}(i)=F(i,i,i,i)$.
    \item $\displaystyle\sum_{i\not=j\not=k\not=l\neq q}F(i,j,k,l,q)= \sum_{i,j,k,l,q}F(i,j,k,l,q)-\sum_{i,j,k,l}F_{2+1+1+1}^{(5)}(i,j,k,l)+\sum_{i,j,k}F_{2+2+1}^{(5)}(i,j,k)$ \\
    $+2\sum_{i,j,k}F_{3+1+1}^{(5)}(i,j,k)-2\sum_{i,j}F_{3+2}^{(5)}(i,j)-6\sum_{i,j}F_{4+1}^{{(5)}}(i,j)+24\sum_{i}F_{5}^{(5)}(i),$
    where
    \begin{align*}
        F_{2+1+1+1}^{(5)}(i,j,k,l) = & F(i,i,j,k,l)+F(i,j,i,k,l)+F(i,j,k,i,l)+F(i,j,k,l,i) + F(j,i,i,k,l)\\
        +& F(j,i,k,i,l)+F(j,i,k,l,i)+F(j,k,i,i,l)+F(j,k,i,l,i)+F(j,k,l,i,i) \\
        F_{2+2+1}^{(5)}(i,j,k) = & F(i,i,j,j,k)+F(i,i,j,k,j)+F(i,i,k,j,j)+F(i,k,i,j,j)+F(k,i,i,j,j) \\
        +& F(i,j,i,j,k)+F(i,j,i,k,j)+F(i,j,k,i,j)+F(i,k,j,i,j)+F(k,i,j,i,j) \\
        +& F(i,j,j,i,k)+F(i,j,j,k,i)+F(i,j,k,j,i)+F(i,k,j,j,i)+F(k,i,j,j,i) \\
        F_{3+1+1}^{(5)}(i,j,k) = & F(i,i,i,j,k)+F(i,i,j,i,k)+F(i,i,j,k,i) + F(i,j,i,i,k)+F(i,j,i,k,i) \\
        + &F(i,j,k,i,i)+F(j,i,i,i,k)+F(j,i,i,k,i)+F(j,i,k,i,i)+F(j,k,i,i,i) \\
        F_{3+2}^{(5)}(i,j) = & F(i,i,i,j,j) + F(i,i,j,i,j)+F(i,i,j,j,i)+F(i,j,i,i,j)+F(i,j,i,j,i)\\
        +&F(i,j,j,i,i) +F(j,i,i,i,j)+F(j,i,i,j,i)+F(j,i,j,i,i)+F(j,j,i,i,i)\\
        F_{4+1}^{(5)}(i,j) = & F(i,i,i,i,j)+F(i,i,i,j,i)+F(i,i,j,i,i)+F(i,j,i,i,i)+F(j,i,i,i,i) \\
        F_{5}^{(5)}(i) = & F(i,i,i,i,i)
    \end{align*}
\end{enumerate}
\end{proposition}

% {\bf \Large  \color{red} Is it possible to add results for $m=5$ here?} 
% { Yes, I have added it before, which is very long.}

\subsection{Contraction function for a graph}\label{sec:faster_contraction}

We may notice that computing the $k$-fold contraction $F_{n_1+\ldots+n_k}^{(m)}$ for a general function $F$ becomes increasingly complicated when $m \ge 4$. To simplify this process, we introduce a more straightforward approach. The contraction function $F_{n_1+\ldots+n_k}^{(m)}(i_1,\ldots,i_k)$ enumerates all possible combinations of the blocks $B_1,\ldots,B_k$, which can be evaluated by assigning distinct index values to each block. Specifically, for a contraction of type $\bm{n}_\sigma = (n_1,n_2,\ldots,n_k)$ corresponding to $M_\sigma$, we select $n_1$ indices from $S = \{i_1,\ldots,i_m\}$ and set them equal in $F$, then select $n_2$ indices from the remaining ones and set them equal, and so on. All possible combinations of such selections collectively define the contraction function $F_{n_1+\ldots+n_k}^{(m)}$. When $n_i=1$ for all $1\le i\le m$, no indices are identified, and the contraction function reduces to $F$ itself and this case is called a \textit{trivial contraction}. Any contraction with at least one $n_j>1$ is referred to as \textit{nontrivial}.

Recall that $F=A_{i_{k_1}i_{k_2}}A_{i_{k_3}i_{k_4}}\ldots A_{i_{k_{2p-1}}i_{k_{2p}}}, \forall i_{k_\ell}\in S,\forall 1\leq \ell\leq 2p$, where each term $A_{i_{\ell_1} i_{\ell_2}}$ in $F$ represents an edge connecting the nodes $i_{\ell_1}$ and $i_{\ell_2}$.
If both indices $i_{\ell_1}$ and $i_{\ell_2}$ belong to the same block $B$ during the selection process, they are forced to be equal, and the corresponding term becomes $A_{i_{\ell_1} i_{\ell_1}} = 0$ since $\operatorname{diag}(A) = 0$. Hence, such terms can be ignored. Consequently, any edge that exists in the graph $G$ should not be selected when evaluating the contraction function $F_{n_1+\ldots+n_k}^{(m)}$.

Now consider a specific block $B$ with size $|B|>1$. To obtain a nonzero contribution, we must find a pair of indices $(i_{\ell_1}, i_{\ell_2})$ that do not appear together in $B$; otherwise, the corresponding term vanishes. Equivalently, $(i_{\ell_1}, i_{\ell_2})$ must not be an edge in $K_B = (B, E_B)$, where $K_B$ denotes the complete graph (clique) generated by $B$ and $E_B$ is the set of all vertex pairs within $B$. In other words, we can instead construct the complement graph of $G$, denoted by $\mathcal{H}$, and search for all possible complete subgraphs contained in this complement. Each of these subgraphs corresponds to a nonzero term in the contraction function $F_{n_1+\ldots+n_k}^{(m)}$. 

If the partition $\sigma$ contains more than one block with sizes greater than 1, the corresponding cliques must be disjoint to each other in $\mathcal{H}$. When such graphs exist, the indices within each block are forced to be equal, and their nonzero contributions are then included in the computation of $F_{n_1+\ldots+n_k}^{(m)}$.

\subsection{Proposed Algorithm}\label{sec:summary}

In this section, we summarize the procedure for evaluating
\begin{align*}
    L = \sum_{i_1\neq i_2\neq \ldots\neq i_m}A_{i_{k_1}i_{k_2}}A_{i_{k_3}i_{k_4}}\ldots A_{i_{k_{2p-1}}i_{k_{2p}}}, \quad i_{k_\ell}\in\{i_1,\ldots,i_m\}, \forall 1\leq\ell\leq 2p.
\end{align*}
through the following steps:
\begin{enumerate}[label=(\alph*)]
    \item Construct the graph $G=(S,E)$ where the edge set is $E=\{(i_{k_{2\ell-1}},i_{k_{2\ell}})\}^{p}_{\ell=1}$ be the set of all edges. Then, form the complement graph $\mathcal{H}$ of $G$.
    \item To identify all nonzero contraction functions, proceed as follows:
    \begin{enumerate}[label=(\arabic*)]
        \item Search for any clique of size $n_1 \in \{m-1, \ldots, 2\}$ in $\mathcal{H}$.
        \item Each clique $\Kcal_1=(V_1,E_1)$ corresponds to a term of type $(n_1, 1, \ldots, 1)$ in $M_\sigma$. Its $(m - n_1 + 1)$-fold contraction $F_{n_1+1+\ldots+1}^{(m)}$ is obtained by setting all indices in $V_1$ equal in $F$.
        \item Next, search for any clique of size $n_2 \in \{n_1, \ldots, 2\}$ in the reduced graph $\mathcal{H}_1 = \mathcal{H} / K_1$. Each such clique $K_2 = (V_2, E_2)$ corresponds to the $(n_1, n_2, \ldots, 1)$-type of $M_\sigma$, yielding the contraction $F_{n_1 + n_2 + 1 + \ldots + 1}^{(m)}.$
        \item Continue this process for $n_3,...$ until all possible partitions $\bm{n}_\sigma = (n_1, \ldots, n_k)$ leading to nonzero contraction functions have been found. 
    \end{enumerate}
    \item Finally, compute $L=\sum_{\bm{n}_\sigma}\mu_{\bm{n}_\sigma}M_{\sigma}$,
    where $\mu_{\bm{n}_\sigma}=(-1)^{m-k}\prod_{i=1}^k(n_i-1)!$ and $M_\sigma$ is given by \eqref{eq:M_sig}. The sum is taken over all possible vectors $\bm{n}_\sigma$ satisfying the conditions identified in step (b).
\end{enumerate}

\section{Subgraph counts in weighted graphs}
In this section, we present computationally efficient equivalent forms for the counts of all weighted graphs in \eqref{eq:target_sum} for $m \le 5$.

\begin{table}
\centering
\begin{tabular}{>{\centering\arraybackslash}m{3cm} @{\hspace{-0.3cm}} >{\centering\arraybackslash}m{3cm} c c}
\toprule
\textbf{Graph} & \textbf{Complement} & \textbf{Sum}  & \textbf{Matrix Formula}\\
\midrule

% --- Row 1 ---
\adjustbox{valign=m}{
\begin{tikzpicture}[baseline=(current bounding box.center),every node/.style={circle,draw,inner sep=0pt,minimum size=12pt,text centered}, scale=0.8]
    \node (1) at (0,0) {k};
    \node (2) at (1,0) {j};
    \node (3) at (0.5,0.866) {i};
    \draw (1)--(2)--(3)--(1);
\end{tikzpicture}
}
&
\adjustbox{valign=m}{
\begin{tikzpicture}[baseline=(current bounding box.center),every node/.style={circle,draw,inner sep=0pt,minimum size=12pt,text centered}, scale=0.8]
    \node (1) at (0,0) {k};
    \node (2) at (1,0) {j};
    \node (3) at (0.5,0.866) {i};
    % \draw (1)--(2)--(3)--(1);
\end{tikzpicture}
}
& $\displaystyle\sum_{i\neq j\neq k}A_{ij}A_{jk}A_{ki}$ & $\tr(A^3)$\\

\midrule
% --- Row 2 ---
\adjustbox{valign=m}{
\begin{tikzpicture}[baseline=(current bounding box.center),every node/.style={circle,draw,inner sep=0pt,minimum size=12pt,text centered}, scale=0.8]
    \node (1) at (0,0) {k};
    \node (2) at (1,0) {j};
    \node (3) at (0.5,0.866) {i};
    \draw (1)--(3)--(2);
\end{tikzpicture}
}
& 
\adjustbox{valign=m}{
\begin{tikzpicture}[baseline=(current bounding box.center),every node/.style={circle,draw,inner sep=0pt,minimum size=12pt,text centered}, scale=0.8]
    \node (1) at (0,0) {k};
    \node (2) at (1,0) {j};
    \node (3) at (0.5,0.866) {i};
    \draw (1)--(2);
\end{tikzpicture}
}
& $\displaystyle\sum_{i\neq j\neq k}A_{ij}A_{ik}$ & $\onebm^\top (A^2- A^{\circ2})\onebm$ \\

\bottomrule
\end{tabular}
\caption{All connected 3-node subgraphs and their counting formulas given in an efficient form.}
\label{tab:result_3_nodes}
\end{table}

\subsection{Number of weighted subgraphs with 3 nodes}
There are two distinct connected graphs on three nodes, the counts of them are given by $\sum_{i\neq j\neq k}A_{ij}A_{jk}A_{ki}$ and $\sum_{i\neq j\neq k}A_{ij}A_{ik}$.
They can be computed by \Cref{prop:3_node} and are summarized in \Cref{tab:result_3_nodes}.
\begin{theorem}\label{prop:3_node}
The following identities hold:
\begin{enumerate}[label=(\alph*)]
    \item $\displaystyle\sum_{i\neq j\neq k}A_{ij}A_{jk}A_{ki}= \tr(A^3)$,
    \item $\displaystyle\sum_{i\neq j\neq k}A_{ij}A_{ik} = \onebm^\top (A^2- A^{\circ2})\onebm$.
\end{enumerate}
\end{theorem}
\begin{proof}
Let $V=\{i,j,k\}$ be the set of nodes.
\begin{enumerate}[label=(\alph*)]
\item All nodes $i,j,k$ are connected, so there is no edge existing in the complement graph and thus we conclude that $F_{2+1}^{(3)}=F_3^{(3)}=0$. Hence, 
\begin{align*}
    \sum_{i\neq j\neq k}A_{ij}A_{jk}A_{ki} & = \sum_{i,j,k}A_{ij}A_{jk}A_{ki} = \tr(A^3).
\end{align*}
\item The complement graph only contains one edge, $(j,k)$. So $F_3^{(3)}$ and there is only term is nonzero in $F_{2+1}^{(3)}$, to calculate it, we let $k=j$ and 
\begin{align*}
    {F}_{2+1}^{(3)}(i,j)&=F(i,j,j)= A_{ij}^2.
\end{align*}
Following the remaining steps in \Cref{sec:summary}, we have that 
\begin{align*}
    \sum_{i\not=j\not=k}A_{ij}A_{ik} & = \sum_{i,j,k}A_{ij}A_{ik} - \sum_{i,j}A_{ij}^2 = \onebm^\top A^2\onebm-\onebm^\top A^{\circ2}\onebm=\onebm^\top (A^2- A^{\circ2})\onebm.
\end{align*}
\end{enumerate}
    
\end{proof}

\subsection{Number of weighted subgraphs of 4 nodes}

\begin{table}[!ht]
\centering
\resizebox{\textwidth}{!}{
\begin{tabular}{>{\centering\arraybackslash}m{2.5cm} >{\centering\arraybackslash}m{2.5cm} c p{6cm}}
\toprule
\textbf{Graph} & \textbf{Complement} & \textbf{Sum}  & \textbf{Matrix Formula}\\
\midrule

% --- Row 1 ---
\adjustbox{valign=m}{
\begin{tikzpicture}[baseline=(current bounding box.center),every node/.style={circle,draw,inner sep=0pt,minimum size=12pt,text centered}, node distance=0.5cm, scale=0.8]
    \node (a1) at (0,0) {i};
    \node (a2) [above=of a1] {j};
    \node (a3) [right=of a2] {k};
    \node (a4) [below=of a3] {l};
    \draw (a1)--(a2)--(a3)--(a4);
\end{tikzpicture}
}
&
\adjustbox{valign=m}{
\begin{tikzpicture}[baseline=(current bounding box.center),every node/.style={circle,draw,inner sep=0pt,minimum size=12pt,text centered}, node distance=0.5cm, scale=0.8]
    \node (a1) at (0,0) {i};
    \node (a2) [above=of a1] {j};
    \node (a3) [right=of a2] {k};
    \node (a4) [below=of a3] {l};
    \draw (a1)--(a4)--(a2);
    \draw (a1)--(a3);
\end{tikzpicture}
}
& $\displaystyle\sum_{i\neq j\neq k\neq l}A_{ij}A_{jk}A_{kl}$ & $\onebm^\top (A^3-2AA^{\circ2}+A^{\circ3})\onebm-\tr(A^3)$\\

\midrule
% --- Row 2 ---
\adjustbox{valign=m}{
\begin{tikzpicture}[baseline=(current bounding box.center),every node/.style={circle,draw,inner sep=0pt,minimum size=12pt,text centered}, node distance=0.5cm, scale=0.8]
    \node (a1) at (0,0) {i};
    \node (a2) [above=of a1] {j};
    \node (a3) [right=of a2] {k};
    \node (a4) [below=of a3] {l};
    \draw (a1)--(a2)--(a3)--(a4)--(a1)--(a3);
\end{tikzpicture}
}
&
\adjustbox{valign=m}{
\begin{tikzpicture}[baseline=(current bounding box.center),every node/.style={circle,draw,inner sep=0pt,minimum size=12pt,text centered}, node distance=0.5cm, scale=0.8]
    \node (a1) at (0,0) {i};
    \node (a2) [above=of a1] {j};
    \node (a3) [right=of a2] {k};
    \node (a4) [below=of a3] {l};
    \draw (a4)--(a2);
\end{tikzpicture}
}
& $\displaystyle\sum_{i\neq j\neq k\neq l}A_{ij}A_{jk}A_{kl}A_{li}A_{ik}$ & $\onebm^\top [(A^2)^{\circ2}\odot A-(A^{\circ2})^2\odot A]\onebm$\\

\midrule
% --- Row 3 ---
\adjustbox{valign=m}{
\begin{tikzpicture}[baseline=(current bounding box.center),every node/.style={circle,draw,inner sep=0pt,minimum size=12pt,text centered}, node distance=0.5cm, scale=0.8]
    \node (a1) at (0,0) {i};
    \node (a2) [above=of a1] {j};
    \node (a3) [right=of a2] {k};
    \node (a4) [below=of a3] {l};
    \draw (a1)--(a2)--(a3)--(a4)--(a1);
\end{tikzpicture}
}
&
\adjustbox{valign=m}{
\begin{tikzpicture}[baseline=(current bounding box.center),every node/.style={circle,draw,inner sep=0pt,minimum size=12pt,text centered}, node distance=0.5cm, scale=0.8]
    \node (a1) at (0,0) {i};
    \node (a2) [above=of a1] {j};
    \node (a3) [right=of a2] {k};
    \node (a4) [below=of a3] {l};
    \draw (a4)--(a2);
    \draw (a1)--(a3);
\end{tikzpicture}
}
& $\displaystyle\sum_{i\neq j\neq k\neq l}A_{ij}A_{jk}A_{kl}A_{li}$ & $\tr(A^4)+\onebm^\top[A^{\circ4}-2(A^{\circ2})^2]\onebm$\\
\midrule

% --- Row 4 ---
\adjustbox{valign=m}{
\begin{tikzpicture}[baseline=(current bounding box.center),every node/.style={circle,draw,inner sep=0pt,minimum size=12pt,text centered}, node distance=0.5cm, scale=0.8]
    \node (a1) at (0,0) {i};
    \node (a2) [above=of a1] {j};
    \node (a3) [right=of a2] {k};
    \node (a4) [below=of a3] {l};
    \draw (a4)--(a1)--(a3);
    \draw (a2)--(a1);
\end{tikzpicture}
}
&
\adjustbox{valign=m}{
\begin{tikzpicture}[baseline=(current bounding box.center),every node/.style={circle,draw,inner sep=0pt,minimum size=12pt,text centered}, node distance=0.5cm, scale=0.8]
    \node (a1) at (0,0) {i};
    \node (a2) [above=of a1] {j};
    \node (a3) [right=of a2] {k};
    \node (a4) [below=of a3] {l};
    \draw (a2)--(a3)--(a4);
    \draw (a2)--(a4);
\end{tikzpicture}
}
& $\displaystyle\sum_{i\neq j\neq k\neq l}A_{ij}A_{ik}A_{il}$ & $\onebm^\top(A\onebm)^{\circ3}+\onebm^\top (2A^{\circ3}-3A^{\circ2}A)\onebm$\\

\midrule

% --- Row 5 ---
\adjustbox{valign=m}{
\begin{tikzpicture}[baseline=(current bounding box.center),every node/.style={circle,draw,inner sep=0pt,minimum size=12pt,text centered}, node distance=0.5cm, scale=0.8]
    \node (a1) at (0,0) {i};
    \node (a2) [above=of a1] {j};
    \node (a3) [right=of a2] {k};
    \node (a4) [below=of a3] {l};
    \draw (a1)--(a2)--(a3)--(a1);
    \draw (a1)--(a4);
\end{tikzpicture}
}
&
\adjustbox{valign=m}{
\begin{tikzpicture}[baseline=(current bounding box.center),every node/.style={circle,draw,inner sep=0pt,minimum size=12pt,text centered}, node distance=0.5cm, scale=0.8]
    \node (a1) at (0,0) {i};
    \node (a2) [above=of a1] {j};
    \node (a3) [right=of a2] {k};
    \node (a4) [below=of a3] {l};
    \draw (a3)--(a4)--(a2);
\end{tikzpicture}
}
& $\displaystyle\sum_{i\neq j\neq k\neq l}A_{ij}A_{jk}A_{ki}A_{il}$ & $\diag(A^3)^\top A\onebm-2\tr(A^{\circ2}A^2)$\\

\midrule
% --- Row 6 ---
\adjustbox{valign=m}{
\begin{tikzpicture}[baseline=(current bounding box.center),every node/.style={circle,draw,inner sep=0pt,minimum size=12pt,text centered}, node distance=0.5cm, scale=0.8]
    \node (a1) at (0,0) {i};
    \node (a2) [above=of a1] {j};
    \node (a3) [right=of a2] {k};
    \node (a4) [below=of a3] {l};
    \draw (a1)--(a2)--(a3)--(a4)--(a1)--(a3);
    \draw (a2)--(a4);
\end{tikzpicture}
}
&
\adjustbox{valign=m}{
\begin{tikzpicture}[baseline=(current bounding box.center),every node/.style={circle,draw,inner sep=0pt,minimum size=12pt,text centered}, node distance=0.5cm, scale=0.8]
    \node (a1) at (0,0) {i};
    \node (a2) [above=of a1] {j};
    \node (a3) [right=of a2] {k};
    \node (a4) [below=of a3] {l};
\end{tikzpicture}
}
& $\displaystyle\sum_{i\neq j\neq k\neq l}A_{ij}A_{jk}A_{kl}A_{li}A_{ik}A_{jl}$ & $\langle (\Acal_{(2)}\odot \Acal_{(3)})\times_{1,2}(\Acal_{(2)}\odot \Acal_{(3)})\times_{1,2}(\Acal_{(2)}\odot \Acal_{(3)}),\Ical_{n\times n\times n}\rangle$\\

\bottomrule
\end{tabular}
}
\caption{All connected 4-node subgraphs and their counting formulas given in an efficient form.}
\label{tab:result_4_nodes}
\end{table}

All distinct connected 4-node subgraphs and their counts are listed in \Cref{tab:result_4_nodes}, with their computation provided in \Cref{prop:4-node}.

\begin{theorem}\label{prop:4-node}
Let $\Acal_{(2)}, \Acal_{(3)}\in\Rbb^{n\times n\times n}$ be 3-way tensors satisfying $[\Acal_{(2)}]_{:k':}=[\Acal_{(3)}]_{::k'}=A$ for all $1\leq k'\leq n$. Let $\Ical_{n\times n\times n}\in\Rbb^{n\times n\times n}$ denote the tensor whose frontal slices are the identity matrix, i.e. $[\Ical_{n\times n\times n}]_{:,:,k'}=I_n$.

The following identities hold:
\begin{enumerate}[label=(\alph*)]
    \item $\displaystyle\sum_{i\neq j\neq k\neq l}A_{ij}A_{jk}A_{kl}= \onebm^\top (A^3-2AA^{\circ2}+A^{\circ3})\onebm-\tr(A^3)$,
    \item $\displaystyle\sum_{i\neq j\neq k\neq l}A_{ij}A_{jk}A_{kl}A_{li}A_{ik}= \onebm^\top [(A^2)^{\circ2}\odot A-(A^{\circ2})^2\odot A]\onebm$,
    \item $\displaystyle\sum_{i\neq j\neq k\neq l}A_{ij}A_{jk}A_{kl}A_{li}= \tr(A^4)+\onebm^\top[A^{\circ4}-2(A^{\circ2})^2]\onebm$, 
    \item $\displaystyle\sum_{i\neq j\neq k\neq l}A_{ij}A_{ik}A_{il} = \onebm^\top(A\onebm)^{\circ3}+\onebm^\top (2A^{\circ3}-3A^{\circ2}A)\onebm$,
    \item $\displaystyle\sum_{i\neq j\neq k\neq l}A_{ij}A_{jk}A_{ki}A_{il}=\diag(A^3)^\top A\onebm-2\tr(A^{\circ2}A^2)$,
    \item $\displaystyle\sum_{i\neq j\neq k\neq l}A_{ij}A_{jk}A_{kl}A_{li}A_{ik}A_{jl}=\langle (\Acal_{(2)}\odot \Acal_{(3)})\times_{1,2}(\Acal_{(2)}\odot \Acal_{(3)})\times_{1,2}(\Acal_{(2)}\odot \Acal_{(3)}),\Ical_{n\times n\times n}\rangle$.
\end{enumerate}
\end{theorem}

\begin{proof}
Let $V=\{i,j,k,l\}$ be the set of nodes.
\begin{enumerate}[label=(\alph*)]
\item Let $F=A_{ij}A_{jk}A_{kl}$. The corresponding graph of $A_{ij}A_{jk}A_{kl}$ contains three edge (see Row 1 of \Cref{tab:result_4_nodes}) and its complement is $H=\{V,E\}$ where $E=\{(j,l),(i,k),(i,l)\}$. Thus there are two nonzero nontrivial contraction function $F_{2+1+1}^{(4)}$ and $F_{2+2}^{(4)}$, which are given by 
\begin{align*}
    F_{2+1+1}^{(4)} & = F(i,j,k,j) + F(i,j,i,l) + F(i,j,k,i)\sim F(i,j,k,j) + F(i,j,i,k) + F(i,j,k,i) \\
    & = A_{ij}A_{jk}^2 + A_{ij}^2A_{ik} + A_{ij}A_{jk}A_{ki}, \\
    F_{2+2}^{(4)} & = F(i,j,i,j) = A_{ij}^3.
\end{align*}
Hence, following the formula in \Cref{sec:summary} one has that 
\begin{align*}
    & \sum_{i\neq j\neq k}A_{ij}A_{jk}A_{kl} = \sum_{i,j,k,l}A_{ij}A_{jk}A_{kl}-\sum_{i,j,k}(A_{ij}A_{jk}^2 + A_{ij}^2A_{ik} + A_{ij}A_{jk}A_{ki})+\sum_{i,j}A_{ij}^3\\
    & = \onebm^\top A^3\onebm-2\onebm^\top AA^{\circ2}\onebm -\tr(A^3)+\onebm^\top A^{\circ3}\onebm = \onebm^\top (A^3-2AA^{\circ2}+A^{\circ3})\onebm-\tr(A^3).
\end{align*}
\item Let $F=A_{ij}A_{jk}A_{kl}A_{li}A_{ik}$. The complement graph $H=(V,E)$ with $E=\{(j,l)\}$ and there is only one nonzero nontrivial contraction function $F_{2+1+1}^{(4)}$. By letting $l=k$, we have that
\begin{align*}
    F_{2+1+1}^{(4)} = F(i,j,k,j) = A_{ij}A_{jk}A_{kj}A_{ji}A_{ik}=A_{ij}^2A_{jk}^2A_{ik}.
\end{align*}
Hence, the sum can be evaluated as follows:
\begin{align*}
    &\sum_{i\neq j\neq k\neq l}A_{ij}A_{jk}A_{kl}A_{li}A_{ik} = \sum_{i,j,k,l}A_{ij}A_{jk}A_{kl}A_{li}A_{ik}-\sum_{i,j,k}A_{ij}^2A_{jk}^2A_{ik} \\
    & = \sum_{i,k}\left(\sum_{j}A_{ij}A_{jk} \right)A_{ik}\left(\sum_lA_{kl}A_{li} \right) - \sum_{i,k}\left( \sum_jA_{ij}^2A_{jk}^2\right)A_{ik} \\
    & = \sum_{i,k} [A^2\odot A\odot A^2]_{ik}-\sum_{i,k} [(A^{\circ2})^2]_{ik}A_{ik}=\onebm^\top [(A^2)^{\circ2}\odot A-(A^{\circ2})^2\odot A]\onebm.
\end{align*}
\item Let $F=A_{ij}A_{jk}A_{kl}A_{li}$. Similar to the part (a), we find the complement graph contains the edges, $(i,k)$ and $(j,l)$. Thus the contraction functions are given by
\begin{align*}
    F_{2+1+1}^{(4)} & = F(i,j,i,l) + F(i,j,k,j) \sim F(i,j,i,k) + F(i,j,k,j) \\
    & = A_{ij}^2A_{ik}^2+A_{ij}^2A_{jk}^2, \\
    F_{2+2}^{(4)} & = F(i,j,i,j) = A_{ij}^4,
\end{align*}
and the total sum can be found by
\begin{align*}
    \sum_{i\neq j\neq k\neq l}A_{ij}A_{jk}A_{kl}A_{li}& = \sum_{i,j,k,l}A_{ij}A_{jk}A_{kl}A_{li}-\sum_{i,j,k}(A_{ij}^2A_{ik}^2+A_{ij}^2A_{jk}^2)+ \sum_{i,j,k}A_{ij}^4\\
    & = \tr(A^4)-2\onebm^\top(A^{\circ2})^2\onebm+\onebm^\top(A^{\circ4})\onebm = \tr(A^4)+\onebm^\top[A^{\circ4}-2(A^{\circ2})^2]\onebm.
\end{align*}
\item Let $F=A_{ij}A_{ik}A_{il}$ The complement graph is $H=(V,E)$ where $E=\{j-l-k-j\}$ and it is also a clique with size 3. Thus, the two nontrivial contraction functions are given by
\begin{align*}
    F_{2+1+1}^{(4)} & = F(i,j,k,j)+F(i,j,k,k)+F(i,j,j,l)\sim F(i,j,k,j)+F(i,j,k,k)+F(i,j,j,k) \\
    & = 2A_{ij}^2A_{ik}+A_{ij}A_{ik}^2\sim 3A_{ij}^2A_{ik}, \\
    F_{3+1}^{(4)}& = F(i,j,j,j) = A_{ij}^3.
\end{align*}
Using the formula in \Cref{sec:summary}, we have that 
\begin{align*}
    &\sum_{i\neq j\neq k\neq l} A_{ij}A_{ik}A_{il} = \sum_{i,j,k,l}A_{ij}A_{ik}A_{il} - \sum_{i,j,k}3A_{ij}^2A_{ik}+2\sum_{i,j}A_{ij}^3 \\
    & = \sum_i\left(\sum_jA_{ij}\right)\left(\sum_kA_{ik}\right)\left(\sum_lA_{il}\right)-3\sum_i\left(\sum_jA_{ij}^2\right)\left(\sum_kA_{ik}\right)+2\onebm^\top A^{\circ3}\onebm \\
    & = \onebm^\top (A\onebm)^{\circ3}-3 \onebm^\top (A^{\circ2}A)\onebm+2\onebm^\top A^{\circ3}\onebm \\
    & = \onebm^\top(A\onebm)^{\circ3}+\onebm^\top (2A^{\circ3}-3A^{\circ2}A)\onebm.
\end{align*}
\item Let $F=A_{ij}A_{jk}A_{ki}A_{il}$ and the complement graph is $H=(V,E)$ where $E=\{(j,l),(k,l)\}$. So the only one nontrivial contraction function can be found by
\begin{align*}
    F_{2+1+1}^{(4)} & = F(i,j,k,j)+F(i,j,k,k) = A_{ij}^2A_{jk}A_{ki}+A_{ij}A_{jk}A_{ki}^2.
\end{align*}
Therefore,
\begin{align*}
    &\sum_{i\neq j\neq k\neq l} A_{ij}A_{jk}A_{ki}A_{il} = \sum_{i,j,k,l}A_{ij}A_{jk}A_{ki}A_{il} - \sum_{i,j,k}(A_{ij}^2A_{jk}A_{ki}+A_{ij}A_{jk}A_{ki}^2) \\
    & = \sum_i\left( \sum_{j,k}A_{ij}A_{jk}A_{ki} \right)\left( \sum_{l}A_{il} \right) - \sum_i\left( \sum_{j,k}A_{ij}^2A_{jk}A_{ki} \right)-\sum_i\left( \sum_{j,k}A_{ij}A_{jk}A_{ki}^2 \right) \\
    & = \sum_i[A^3]_{ii}[A\onebm]_i-\tr(A^{\circ2}A^2)-\tr(A^2 A^{\circ2}) \\
    & = \diag(A^3)^\top A\onebm-2\tr(A^{\circ2}A^2).
\end{align*}
\item Let $F=A_{ij}A_{jk}A_{kl}A_{li}A_{ik}A_{jl}$. Since the complement graph contains no cliques of size $\geq 2$, the sum over distinct indices is equal to the unrestricted sum.
To proceed, we fix $l$ and $F$ can be considered as the product of $A_{ij}A_{li},A_{jk}A_{lj}$ and $A_{ki}A_{lk}$. Define the matrix $[T^l]_{ij}=A_{ij}A_{li}$, then 
\begin{align*}
    \sum_{i,j,k} = A_{ij}A_{li}A_{jk}A_{lj}A_{ki}A_{lk} =\sum_{i,j,k}[T^l]_{ij}[T^l]_{jk}[T^l]_{ki}= \tr((T^l)^3),
\end{align*}
which implies that 
\begin{align*}
    \sum_{i,j,k,l}A_{ij}A_{li}A_{jk}A_{lj}A_{ki}A_{lk} =\sum_l\tr((T^l)^3).
\end{align*}
We further define two tensors $\Acal_{(2)}, \Acal_{(3)}\in\Rbb^{n\times n\times n}$ by letting all lateral slices of $\Acal_{(2)}$ be equal to $A$ and all frontal slices of $\Acal_{(3)}$ be equal to $A$; that is, $[\Acal_{(2)}]_{:k':}=[\Acal_{(3)}]_{::k'}=A$ for all $1\leq k'\leq n$. So $[\Acal_{(2)}]_{i,j,k}=A_{i,k}$ and $[\Acal_{(3)}]_{i,j,k}=A_{i,j}$, which implies that $[\Acal_{(2)}\odot \Acal_{(3)}]_{i,j,k}=[\Acal_{(3)}]_{i,j,k}[\Acal_{(2)}]_{i,j,k}=A_{ij}A_{ik}$, i.e., $[\Acal_{(2)}\odot \Acal_{(3)}]_{::l}=T^l$.
Let $\Ical_{n\times n\times n}\in\Rbb^{n\times n\times n}$ denote the tensor whose frontal slices are the identity matrix, i.e. $[\Ical_{n\times n\times n}]_{::k'}=I_n$. Then we have that 
\begin{align*}
    \sum_{i\neq j\neq k\neq l}A_{ij}A_{jk}A_{kl}A_{li}A_{ik}A_{jl} & = \sum_{i,j,k,l}A_{ij}A_{li}A_{jk}A_{lj}A_{ki}A_{lk} \\
    & = \langle (\Acal_{(2)}\odot \Acal_{(3)})\times_{1,2}(\Acal_{(2)}\odot \Acal_{(3)})\times_{1,2}(\Acal_{(2)}\odot \Acal_{(3)}),\Ical_{n\times n\times n}\rangle.
\end{align*}

\end{enumerate}

\end{proof}

\subsection{Number of weighted subgraphs of 5 nodes}
All distinct connected 5-node subgraphs and their counts are listed in \Cref{tab:result_5_nodes-1}-\Cref{tab:result_5_nodes-3}, with their computation provided in \Cref{prop:5-node}.

\begin{table}[!ht]
\centering
\resizebox{\textwidth}{!}{
\begin{tabular}{>{\centering\arraybackslash}m{3cm} >{\centering\arraybackslash}m{3cm} c p{7cm}}
\toprule
\textbf{Graph} & \textbf{Complement} & \textbf{Sum}  & \textbf{Matrix Formula}\\
\midrule

% --- Row 1 ---
\adjustbox{valign=m}{
\begin{tikzpicture}[every node/.style={circle, draw, inner sep=2pt, minimum size=8pt}]
\node (i) at (0,0) {i};
\node (j) at ([xshift=0.974cm, yshift=-0.709cm] i) {j};
\node (k) at ([xshift=-0.372cm, yshift=-0.927cm] j) {k};
\node (l) at ([xshift=-1.204cm, yshift=0cm] k) {l};
\node (q) at ([xshift=-0.372cm, yshift=0.927cm] l) {q};
\draw (i)--(j)--(k)--(l)--(q);

\end{tikzpicture}
}
&
\adjustbox{valign=m}{
\begin{tikzpicture}[every node/.style={circle, draw, inner sep=2pt, minimum size=8pt}]
\node (i) at (0,0) {i};
\node (j) at ([xshift=0.974cm, yshift=-0.709cm] i) {j};
\node (k) at ([xshift=-0.372cm, yshift=-0.927cm] j) {k};
\node (l) at ([xshift=-1.204cm, yshift=0cm] k) {l};
\node (q) at ([xshift=-0.372cm, yshift=0.927cm] l) {q};
\draw (i)--(k)--(q)--(j)--(l)--(i)--(q);

\end{tikzpicture}
}
& $\displaystyle\sum_{i\neq j\neq k\neq l\neq q}A_{ij}A_{jk}A_{kl}A_{lq}$ & $\onebm^\top \Big[A^4-2A^2A^{\circ2}-A\diag(A^{\circ2}\onebm)A+3(A^{\circ2})^2+2A^{\circ3}A-2A^{\circ4}\Big]\onebm-2\diag(A^3)^\top A\onebm-\tr(A^4)+3\tr(A^{\circ2}A^2)$\\

\midrule
% --- Row 2 ---
\adjustbox{valign=m}{
\begin{tikzpicture}[every node/.style={circle, draw, inner sep=2pt, minimum size=8pt}]
\node (i) at (0,0) {i};
\node (j) at ([xshift=0.974cm, yshift=-0.709cm] i) {j};
\node (k) at ([xshift=-0.372cm, yshift=-0.927cm] j) {k};
\node (l) at ([xshift=-1.204cm, yshift=0cm] k) {l};
\node (q) at ([xshift=-0.372cm, yshift=0.927cm] l) {q};
\draw (i)--(j)--(k)--(l)--(q)--(i);

\end{tikzpicture}
}
&
\adjustbox{valign=m}{
\begin{tikzpicture}[every node/.style={circle, draw, inner sep=2pt, minimum size=8pt}]
\node (i) at (0,0) {i};
\node (j) at ([xshift=0.974cm, yshift=-0.709cm] i) {j};
\node (k) at ([xshift=-0.372cm, yshift=-0.927cm] j) {k};
\node (l) at ([xshift=-1.204cm, yshift=0cm] k) {l};
\node (q) at ([xshift=-0.372cm, yshift=0.927cm] l) {q};
\draw (i)--(k)--(q)--(j)--(l)--(i);

\end{tikzpicture}
}
& $\displaystyle\sum_{i\neq j\neq k\neq l\neq q}A_{ij}A_{jk}A_{kl}A_{lq}A_{qi}$ & $\tr(A^5)-5\onebm^\top A^{\circ2}\diag(A^3)+5\tr(A^{\circ3}A^2)$\\

\midrule
% --- Row 3 ---
\adjustbox{valign=m}{
\begin{tikzpicture}[every node/.style={circle, draw, inner sep=2pt, minimum size=8pt}]
\node (i) at (0,0) {i};
\node (j) at ([xshift=0.974cm, yshift=-0.709cm] i) {j};
\node (k) at ([xshift=-0.372cm, yshift=-0.927cm] j) {k};
\node (l) at ([xshift=-1.204cm, yshift=0cm] k) {l};
\node (q) at ([xshift=-0.372cm, yshift=0.927cm] l) {q};
\draw (q)--(i)--(j)--(k);
\draw (l)--(i);

\end{tikzpicture}
}
&
\adjustbox{valign=m}{
\begin{tikzpicture}[every node/.style={circle, draw, inner sep=2pt, minimum size=8pt}]
\node (i) at (0,0) {i};
\node (j) at ([xshift=0.974cm, yshift=-0.709cm] i) {j};
\node (k) at ([xshift=-0.372cm, yshift=-0.927cm] j) {k};
\node (l) at ([xshift=-1.204cm, yshift=0cm] k) {l};
\node (q) at ([xshift=-0.372cm, yshift=0.927cm] l) {q};
\draw (i)--(k)--(l)--(q)--(j)--(l);
\draw (q)--(k);

\end{tikzpicture}
}
& $\displaystyle\sum_{i\neq j\neq k\neq l\neq q}A_{ij}A_{jk}A_{il}A_{qi}$ & $\onebm^\top \Big[ A\diag(A^{2}\onebm)A-A\diag(A^{\circ2}\onebm)A-2AA^{\circ2}A-A^{\circ2}A^2+4A^{\circ3}A+(A^{\circ2})^2-2A^{\circ4}\Big]\onebm-2\diag(A^3)^\top A\onebm+4\tr(A^{\circ2}A^2)$\\

\midrule
% --- Row 4 ---
\adjustbox{valign=m}{
\begin{tikzpicture}[every node/.style={circle, draw, inner sep=2pt, minimum size=8pt}]
\node (i) at (0,0) {i};
\node (j) at ([xshift=0.974cm, yshift=-0.709cm] i) {j};
\node (k) at ([xshift=-0.372cm, yshift=-0.927cm] j) {k};
\node (l) at ([xshift=-1.204cm, yshift=0cm] k) {l};
\node (q) at ([xshift=-0.372cm, yshift=0.927cm] l) {q};
\draw (i)--(j)--(k)--(l)--(q)--(j);

\end{tikzpicture}
}
&
\adjustbox{valign=m}{
\begin{tikzpicture}[every node/.style={circle, draw, inner sep=2pt, minimum size=8pt}]
\node (i) at (0,0) {i};
\node (j) at ([xshift=0.974cm, yshift=-0.709cm] i) {j};
\node (k) at ([xshift=-0.372cm, yshift=-0.927cm] j) {k};
\node (l) at ([xshift=-1.204cm, yshift=0cm] k) {l};
\node (q) at ([xshift=-0.372cm, yshift=0.927cm] l) {q};
\draw (q)--(i)--(k)--(q);
\draw (i)--(l)--(j);

\end{tikzpicture}
}
& $\displaystyle\sum_{i\neq j\neq k\neq l\neq q}A_{ij}A_{jk}A_{jq}A_{kl}A_{lq}$ & $\diag(A^4)^\top A \onebm+\onebm^\top \Big[ 4A^{\circ3}A^{\circ2}+A^{\circ4}A-(A^2)^{\circ2}\odot A-A\diag(A^{\circ2}\onebm)A^{\circ2}-A(A^{\circ2})^2-2A^{\circ5}\Big]\onebm+\tr(A(A^{\circ2})^2)-2\tr(A^{\circ2}A^3)$\\

\midrule
% --- Row 5 ---
\adjustbox{valign=m}{
\begin{tikzpicture}[every node/.style={circle, draw, inner sep=2pt, minimum size=8pt}]
\node (i) at (0,0) {i};
\node (j) at ([xshift=0.974cm, yshift=-0.709cm] i) {j};
\node (k) at ([xshift=-0.372cm, yshift=-0.927cm] j) {k};
\node (l) at ([xshift=-1.204cm, yshift=0cm] k) {l};
\node (q) at ([xshift=-0.372cm, yshift=0.927cm] l) {q};
\draw (i)--(j)--(k)--(q)--(l)--(j);
\draw (q)--(i);

\end{tikzpicture}
}
&
\adjustbox{valign=m}{
\begin{tikzpicture}[every node/.style={circle, draw, inner sep=2pt, minimum size=8pt}]
\node (i) at (0,0) {i};
\node (j) at ([xshift=0.974cm, yshift=-0.709cm] i) {j};
\node (k) at ([xshift=-0.372cm, yshift=-0.927cm] j) {k};
\node (l) at ([xshift=-1.204cm, yshift=0cm] k) {l};
\node (q) at ([xshift=-0.372cm, yshift=0.927cm] l) {q};
\draw (i)--(k)--(l)--(i);
\draw (q)--(j);

\end{tikzpicture}
}
& $\displaystyle\sum_{i\neq j\neq k\neq l\neq q}A_{ij}A_{jk}A_{jl}A_{kq}A_{lq}A_{qi}$ & $\onebm^\top \Big[ (A^2)^{\circ3}+3A^{\circ4}A^{\circ2})+2(A^{\circ3})^2-2A^{\circ6}\Big]\onebm-\onebm^\top(A^{\circ2}\onebm)^{\circ3}-3\tr((A^{\circ2})^2A^2)$\\

\midrule
% --- Row 6 ---
\adjustbox{valign=m}{
\begin{tikzpicture}[every node/.style={circle, draw, inner sep=2pt, minimum size=8pt}]
\node (i) at (0,0) {i};
\node (j) at ([xshift=0.974cm, yshift=-0.709cm] i) {j};
\node (k) at ([xshift=-0.372cm, yshift=-0.927cm] j) {k};
\node (l) at ([xshift=-1.204cm, yshift=0cm] k) {l};
\node (q) at ([xshift=-0.372cm, yshift=0.927cm] l) {q};
\draw (j)--(q)--(i)--(j)--(k)--(l);
\draw (q)--(i);

\end{tikzpicture}
}
&
\adjustbox{valign=m}{
\begin{tikzpicture}[every node/.style={circle, draw, inner sep=2pt, minimum size=8pt}]
\node (i) at (0,0) {i};
\node (j) at ([xshift=0.974cm, yshift=-0.709cm] i) {j};
\node (k) at ([xshift=-0.372cm, yshift=-0.927cm] j) {k};
\node (l) at ([xshift=-1.204cm, yshift=0cm] k) {l};
\node (q) at ([xshift=-0.372cm, yshift=0.927cm] l) {q};
\draw (j)--(l)--(q)--(k)--(i)--(l);

\end{tikzpicture}
}
& $\displaystyle\sum_{i\neq j\neq k\neq l\neq q}A_{ij}A_{jk}A_{jq}A_{kl}A_{qi}$ & $\diag(A^3)^\top A^2\onebm-2\diag(A^{\circ2}A^2)^\top A\onebm-\diag(A^3)^\top A^{\circ2}\onebm-2\onebm^\top ((A^2)^{\circ2}\odot A)\onebm+2\tr(A^{\circ3}A^2)+2\tr((A^{\circ2})^2A)$\\

\midrule
% --- Row 7 ---
\adjustbox{valign=m}{
\begin{tikzpicture}[every node/.style={circle, draw, inner sep=2pt, minimum size=8pt}]
\node (i) at (0,0) {i};
\node (j) at ([xshift=0.974cm, yshift=-0.709cm] i) {j};
\node (k) at ([xshift=-0.372cm, yshift=-0.927cm] j) {k};
\node (l) at ([xshift=-1.204cm, yshift=0cm] k) {l};
\node (q) at ([xshift=-0.372cm, yshift=0.927cm] l) {q};
\draw (j)--(q)--(i)--(j)--(k);
\draw (l)--(q);

\end{tikzpicture}
}
&
\adjustbox{valign=m}{
\begin{tikzpicture}[every node/.style={circle, draw, inner sep=2pt, minimum size=8pt}]
\node (i) at (0,0) {i};
\node (j) at ([xshift=0.974cm, yshift=-0.709cm] i) {j};
\node (k) at ([xshift=-0.372cm, yshift=-0.927cm] j) {k};
\node (l) at ([xshift=-1.204cm, yshift=0cm] k) {l};
\node (q) at ([xshift=-0.372cm, yshift=0.927cm] l) {q};
\draw (i)--(k)--(l)--(i);
\draw (q)--(k);
\draw (j)--(l);
\end{tikzpicture}
}
& $\displaystyle\sum_{i\neq j\neq k\neq l\neq q}A_{ij}A_{jk}A_{jq}A_{lq}A_{qi}$ & $\onebm^\top \Big[ A(A^2\odot A)A-(A^2)^{\circ2}\odot A\Big]\onebm{-2\diag(A^{\circ2}A^2)^\top A\onebm-2\diag(AA^{\circ2}A)^\top A\onebm}+4\tr((A^{\circ2})^2A)+\tr(A^{\circ3}A^2)$\\

\bottomrule
\end{tabular}
}
\caption{All connected 5-node subgraphs and their counting formulas given in an efficient form, Part I.}
\label{tab:result_5_nodes-1}
\end{table}

\begin{table}[!ht]
\centering
\resizebox{\textwidth}{!}{
\begin{tabular}{>{\centering\arraybackslash}m{3cm} >{\centering\arraybackslash}m{3cm} c p{7cm}}
\toprule
\textbf{Graph} & \textbf{Complement} & \textbf{Sum}  & \textbf{Matrix Formula}\\
\midrule

% --- Row 1 ---
\adjustbox{valign=m}{
\begin{tikzpicture}[every node/.style={circle, draw, inner sep=2pt, minimum size=8pt}]
\node (i) at (0,0) {i};
\node (j) at ([xshift=0.974cm, yshift=-0.709cm] i) {j};
\node (k) at ([xshift=-0.372cm, yshift=-0.927cm] j) {k};
\node (l) at ([xshift=-1.204cm, yshift=0cm] k) {l};
\node (q) at ([xshift=-0.372cm, yshift=0.927cm] l) {q};
\draw (j)--(k)--(l)--(q)--(i)--(j)--(q);

\end{tikzpicture}
}
&
\adjustbox{valign=m}{
\begin{tikzpicture}[every node/.style={circle, draw, inner sep=2pt, minimum size=8pt}]
\node (i) at (0,0) {i};
\node (j) at ([xshift=0.974cm, yshift=-0.709cm] i) {j};
\node (k) at ([xshift=-0.372cm, yshift=-0.927cm] j) {k};
\node (l) at ([xshift=-1.204cm, yshift=0cm] k) {l};
\node (q) at ([xshift=-0.372cm, yshift=0.927cm] l) {q};
\draw (j)--(l)--(i)--(k)--(q);

\end{tikzpicture}
}
& $\displaystyle\sum_{i\neq j\neq k\neq l\neq q}A_{ij}A_{jk}A_{jq}A_{kl}A_{lq}A_{qi}$ & $\onebm^\top \Big[ A\odot A^2\odot A^3-2A\odot A^2\odot (A^{\circ2}A)\Big]\onebm+2\tr(AA^{\circ2}A^{\circ3})+\tr(A^{\circ4}A^2)-2\diag(A^{\circ2}A^2)^\top A^{\circ2}\onebm$\\

\midrule
% --- Row 2 ---
\adjustbox{valign=m}{
\begin{tikzpicture}[every node/.style={circle, draw, inner sep=2pt, minimum size=8pt}]
\node (i) at (0,0) {i};
\node (j) at ([xshift=0.974cm, yshift=-0.709cm] i) {j};
\node (k) at ([xshift=-0.372cm, yshift=-0.927cm] j) {k};
\node (l) at ([xshift=-1.204cm, yshift=0cm] k) {l};
\node (q) at ([xshift=-0.372cm, yshift=0.927cm] l) {q};
\draw (i)--(q)--(l)--(k)--(j)--(q);
\draw (j)--(l);

\end{tikzpicture}
}
&
\adjustbox{valign=m}{
\begin{tikzpicture}[every node/.style={circle, draw, inner sep=2pt, minimum size=8pt}]
\node (i) at (0,0) {i};
\node (j) at ([xshift=0.974cm, yshift=-0.709cm] i) {j};
\node (k) at ([xshift=-0.372cm, yshift=-0.927cm] j) {k};
\node (l) at ([xshift=-1.204cm, yshift=0cm] k) {l};
\node (q) at ([xshift=-0.372cm, yshift=0.927cm] l) {q};
\draw (l)--(i)--(k)--(q);
\draw (i)--(j);

\end{tikzpicture}
}
& $\displaystyle\sum_{i\neq j\neq k\neq l\neq q}A_{jk}A_{jl}A_{jq}A_{kl}A_{lq}A_{qi}$ & $\diag(A(A\odot A^2)A)^\top A\onebm-\diag(A^{\circ2}AA^{\circ2})^\top A\onebm -2\tr(AA^{\circ2}(A\odot A^2))+2\tr(AA^{\circ2}A^{\circ3})-\langle (\Acal_{(2)}\odot \Acal_{(3)})\times_{1,2}(\Acal_{(2)}\odot \Acal_{(3)})\times_{1,2}(\Acal_{(2)}\odot \Acal_{(3)}),\Ical_{n\times n\times n}\rangle$\\

\midrule
% --- Row 3 ---
\adjustbox{valign=m}{
\begin{tikzpicture}[every node/.style={circle, draw, inner sep=2pt, minimum size=8pt}]
\node (i) at (0,0) {i};
\node (j) at ([xshift=0.974cm, yshift=-0.709cm] i) {j};
\node (k) at ([xshift=-0.372cm, yshift=-0.927cm] j) {k};
\node (l) at ([xshift=-1.204cm, yshift=0cm] k) {l};
\node (q) at ([xshift=-0.372cm, yshift=0.927cm] l) {q};
\draw (i)--(j)--(k)--(l)--(q)--(i);
\draw (j)--(l);
\draw (k)--(q);

\end{tikzpicture}
}
&
\adjustbox{valign=m}{
\begin{tikzpicture}[every node/.style={circle, draw, inner sep=2pt, minimum size=8pt}]
\node (i) at (0,0) {i};
\node (j) at ([xshift=0.974cm, yshift=-0.709cm] i) {j};
\node (k) at ([xshift=-0.372cm, yshift=-0.927cm] j) {k};
\node (l) at ([xshift=-1.204cm, yshift=0cm] k) {l};
\node (q) at ([xshift=-0.372cm, yshift=0.927cm] l) {q};
\draw (l)--(i)--(k);
\draw (q)--(j);

\end{tikzpicture}
}
& $\displaystyle\sum_{i\neq j\neq k\neq l\neq q}A_{ij}A_{jk}A_{jl}A_{kl}A_{kq}A_{lq}A_{qi}$ & $\langle (\Acal_{(2)}\odot\Acal_{(3)})\times_{1,2}(\Acal_{(2)}\odot\Acal_{(3)})\times_{1,2}(\Bcal_{(2)}\odot\Acal_{(3)}),\Ical_{n\times n\times n}\rangle-2\tr([(A^{\circ2}A)\odot A]AA^{\circ2})-\diag(A^{\circ2}AA^{\circ2})^\top A^{\circ2}\onebm+2\tr(AA^{\circ2}A^{\circ4})$\\

\midrule
% % --- Row 4 ---
\adjustbox{valign=m}{
\begin{tikzpicture}[every node/.style={circle, draw, inner sep=2pt, minimum size=8pt}]
\node (i) at (0,0) {i};
\node (j) at ([xshift=0.974cm, yshift=-0.709cm] i) {j};
\node (k) at ([xshift=-0.372cm, yshift=-0.927cm] j) {k};
\node (l) at ([xshift=-1.204cm, yshift=0cm] k) {l};
\node (q) at ([xshift=-0.372cm, yshift=0.927cm] l) {q};
\draw (j)--(i)--(k);
\draw (l)--(i)--(q);

\end{tikzpicture}
}
&
\adjustbox{valign=m}{
\begin{tikzpicture}[every node/.style={circle, draw, inner sep=2pt, minimum size=8pt}]
\node (i) at (0,0) {i};
\node (j) at ([xshift=0.974cm, yshift=-0.709cm] i) {j};
\node (k) at ([xshift=-0.372cm, yshift=-0.927cm] j) {k};
\node (l) at ([xshift=-1.204cm, yshift=0cm] k) {l};
\node (q) at ([xshift=-0.372cm, yshift=0.927cm] l) {q};
\draw (j)--(k)--(l)--(q)--(j)--(l);
\draw (q)--(k);

\end{tikzpicture}
}
& $\displaystyle\sum_{i\neq j\neq k\neq l\neq q}A_{ij}A_{ik}A_{il}A_{qi}$ & $\onebm^\top(A\onebm)^{\circ4}-6\onebm^\top A^{\circ2}(A\onebm)^{\circ2}+\onebm^\top \Big[ 3(A^{\circ2})^2+8A^{\circ3}A- 6A^{\circ4} \Big]\onebm$\\

\midrule
% --- Row 5 ---
\adjustbox{valign=m}{
\begin{tikzpicture}[every node/.style={circle, draw, inner sep=2pt, minimum size=8pt}]
\node (i) at (0,0) {i};
\node (j) at ([xshift=0.974cm, yshift=-0.709cm] i) {j};
\node (k) at ([xshift=-0.372cm, yshift=-0.927cm] j) {k};
\node (l) at ([xshift=-1.204cm, yshift=0cm] k) {l};
\node (q) at ([xshift=-0.372cm, yshift=0.927cm] l) {q};
\draw (k)--(j)--(i)--(k);
\draw (l)--(i)--(q);

\end{tikzpicture}
}
&
\adjustbox{valign=m}{
\begin{tikzpicture}[every node/.style={circle, draw, inner sep=2pt, minimum size=8pt}]
\node (i) at (0,0) {i};
\node (j) at ([xshift=0.974cm, yshift=-0.709cm] i) {j};
\node (k) at ([xshift=-0.372cm, yshift=-0.927cm] j) {k};
\node (l) at ([xshift=-1.204cm, yshift=0cm] k) {l};
\node (q) at ([xshift=-0.372cm, yshift=0.927cm] l) {q};
\draw (k)--(l)--(q)--(j)--(l);
\draw (q)--(k);

\end{tikzpicture}
}
& $\displaystyle\sum_{i\neq j\neq k\neq l\neq q}A_{ij}A_{ik}A_{jk}A_{il}A_{qi}$ & $\onebm^\top (A\odot A^2) (A\onebm)^{\circ2}-4\diag(A^2A^{\circ2})^\top A\onebm-\diag(A^3)^\top A^{\circ2}\onebm+2\tr((A^{\circ2})^2A)+4\tr(A^{\circ3}A^2)$\\

\midrule
% --- Row 6 ---
\adjustbox{valign=m}{
\begin{tikzpicture}[every node/.style={circle, draw, inner sep=2pt, minimum size=8pt}]
\node (i) at (0,0) {i};
\node (j) at ([xshift=0.974cm, yshift=-0.709cm] i) {j};
\node (k) at ([xshift=-0.372cm, yshift=-0.927cm] j) {k};
\node (l) at ([xshift=-1.204cm, yshift=0cm] k) {l};
\node (q) at ([xshift=-0.372cm, yshift=0.927cm] l) {q};
\draw (i)--(k)--(j)--(i)--(l)--(q)--(i);
\draw (q)--(i);

\end{tikzpicture}
}
&
\adjustbox{valign=m}{
\begin{tikzpicture}[every node/.style={circle, draw, inner sep=2pt, minimum size=8pt}]
\node (i) at (0,0) {i};
\node (j) at ([xshift=0.974cm, yshift=-0.709cm] i) {j};
\node (k) at ([xshift=-0.372cm, yshift=-0.927cm] j) {k};
\node (l) at ([xshift=-1.204cm, yshift=0cm] k) {l};
\node (q) at ([xshift=-0.372cm, yshift=0.927cm] l) {q};
\draw (k)--(l)--(j)--(q)--(k);

\end{tikzpicture}
}
& $\displaystyle\sum_{i\neq j\neq k\neq l\neq q}A_{ij}A_{ik}A_{jk}A_{il}A_{lq}A_{qi}$ & $\tr((A^3)^{\circ2})-4\onebm^\top ((A^2)^{\circ2}\odot A^{\circ2})\onebm+2\tr((A^{\circ2})^3)$\\

\midrule
% --- Row 7 ---
\adjustbox{valign=m}{
\begin{tikzpicture}[every node/.style={circle, draw, inner sep=2pt, minimum size=8pt}]
\node (i) at (0,0) {i};
\node (j) at ([xshift=0.974cm, yshift=-0.709cm] i) {j};
\node (k) at ([xshift=-0.372cm, yshift=-0.927cm] j) {k};
\node (l) at ([xshift=-1.204cm, yshift=0cm] k) {l};
\node (q) at ([xshift=-0.372cm, yshift=0.927cm] l) {q};
\draw (q)--(l)--(j)--(q)--(k)--(j);
\draw (i)--(j);

\end{tikzpicture}
}
&
\adjustbox{valign=m}{
\begin{tikzpicture}[every node/.style={circle, draw, inner sep=2pt, minimum size=8pt}]
\node (i) at (0,0) {i};
\node (j) at ([xshift=0.974cm, yshift=-0.709cm] i) {j};
\node (k) at ([xshift=-0.372cm, yshift=-0.927cm] j) {k};
\node (l) at ([xshift=-1.204cm, yshift=0cm] k) {l};
\node (q) at ([xshift=-0.372cm, yshift=0.927cm] l) {q};
\draw (i)--(k)--(l)--(i)--(q);

\end{tikzpicture}
}
& $\displaystyle\sum_{i\neq j\neq k\neq l\neq q}A_{ij}A_{jk}A_{jl}A_{jq}A_{kq}A_{lq}$ & $\onebm^\top \Big[A((A^2)^{\circ2}\odot A)-2(A^{\circ2}A)\odot A^2\odot A-A^{\circ2}\odot (A^2)^{\circ2}\Big]\onebm-\onebm^\top A\diag ((A^{\circ2})^2 A)+\tr\Big[(A^{\circ2})^3+2A^{\circ3}A^{\circ2}A\Big]$\\

\bottomrule
\end{tabular}
}
\caption{All connected 5-node subgraphs and their counting formulas given in an efficient form, Part II.}
\label{tab:result_5_nodes-2}
\end{table}

\begin{table}[!ht]
\centering
\resizebox{\textwidth}{!}{
\begin{tabular}{>{\centering\arraybackslash}m{3cm} >{\centering\arraybackslash}m{3cm} c p{7cm}}
\toprule
\textbf{Graph} & \textbf{Complement} & \textbf{Sum}  & \textbf{Matrix Formula}\\
\midrule

% --- Row 1 ---
\adjustbox{valign=m}{
\begin{tikzpicture}[every node/.style={circle, draw, inner sep=2pt, minimum size=8pt}]
\node (i) at (0,0) {i};
\node (j) at ([xshift=0.974cm, yshift=-0.709cm] i) {j};
\node (k) at ([xshift=-0.372cm, yshift=-0.927cm] j) {k};
\node (l) at ([xshift=-1.204cm, yshift=0cm] k) {l};
\node (q) at ([xshift=-0.372cm, yshift=0.927cm] l) {q};
\draw (j)--(i)--(q)--(j)--(l)--(k)--(q)--(l);

\end{tikzpicture}
}
&
\adjustbox{valign=m}{
\begin{tikzpicture}[every node/.style={circle, draw, inner sep=2pt, minimum size=8pt}]
\node (i) at (0,0) {i};
\node (j) at ([xshift=0.974cm, yshift=-0.709cm] i) {j};
\node (k) at ([xshift=-0.372cm, yshift=-0.927cm] j) {k};
\node (l) at ([xshift=-1.204cm, yshift=0cm] k) {l};
\node (q) at ([xshift=-0.372cm, yshift=0.927cm] l) {q};
\draw (l)--(i)--(k)--(j);

\end{tikzpicture}
}
& $\displaystyle\sum_{i\neq j\neq k\neq l\neq q}A_{ij}A_{jl}A_{jq}A_{kl}A_{kq}A_{lq}A_{qi}$ & $\tr\Big[A(A\odot A^2)^2-AA^{\circ2}(A^{\circ2}\odot A^2)-A^2((A^{\circ2}A)\odot A^{\circ2})+A^{\circ3}(A^{\circ2})^2\Big]-\langle (\Acal_{(2)}\odot\Acal_{(3)})\times_{1,2}(\Acal_{(2)}\odot\Acal_{(3)})\times_{1,2}(\Acal_{(2)}\odot\Ccal_{(3)}),\Ical_{n\times n\times n}\rangle $\\

\midrule
% --- Row 2 ---
\adjustbox{valign=m}{
\begin{tikzpicture}[every node/.style={circle, draw, inner sep=2pt, minimum size=8pt}]
\node (i) at (0,0) {i};
\node (j) at ([xshift=0.974cm, yshift=-0.709cm] i) {j};
\node (k) at ([xshift=-0.372cm, yshift=-0.927cm] j) {k};
\node (l) at ([xshift=-1.204cm, yshift=0cm] k) {l};
\node (q) at ([xshift=-0.372cm, yshift=0.927cm] l) {q};
\draw (i)--(j)--(k)--(i)--(q)--(l)--(i);
\draw (j)--(l);
\draw (k)--(q);

\end{tikzpicture}
}
&
\adjustbox{valign=m}{
\begin{tikzpicture}[every node/.style={circle, draw, inner sep=2pt, minimum size=8pt}]
\node (i) at (0,0) {i};
\node (j) at ([xshift=0.974cm, yshift=-0.709cm] i) {j};
\node (k) at ([xshift=-0.372cm, yshift=-0.927cm] j) {k};
\node (l) at ([xshift=-1.204cm, yshift=0cm] k) {l};
\node (q) at ([xshift=-0.372cm, yshift=0.927cm] l) {q};
\draw (q)--(j);
\draw (l)--(k);

\end{tikzpicture}
}
& $\displaystyle\sum_{i\neq j\neq k\neq l\neq q}A_{ij}A_{jk}A_{jl}A_{ki}A_{kq}A_{lq}A_{li}A_{qi}$ & $\langle (\mathscr{A}_{(3,4)}\odot \mathscr{A}_{(2,4)}\odot \mathscr{A}_{(2,3)}\odot \mathscr{A}_{(1,2)})\times_{1,2}(\mathscr{A}_{(3,4)}\odot \mathscr{A}_{(2,4)})\times_{1,2}(\mathscr{A}_{(3,4)}\odot \mathscr{A}_{(2,3)}),\Ical_{n\times n\times n\times n}\rangle+\tr(A^{\circ4}(A^{\circ2})^2)-2\onebm^\top((AA^{\circ2})^{\circ2}\odot A^{\circ2})\onebm$\\

\midrule
% --- Row 3 ---
\adjustbox{valign=m}{
\begin{tikzpicture}[every node/.style={circle, draw, inner sep=2pt, minimum size=8pt}]
\node (i) at (0,0) {i};
\node (j) at ([xshift=0.974cm, yshift=-0.709cm] i) {j};
\node (k) at ([xshift=-0.372cm, yshift=-0.927cm] j) {k};
\node (l) at ([xshift=-1.204cm, yshift=0cm] k) {l};
\node (q) at ([xshift=-0.372cm, yshift=0.927cm] l) {q};
\draw (i)--(j)--(k)--(l)--(q)--(j)--(l);
\draw (q)--(k);

\end{tikzpicture}
}
&
\adjustbox{valign=m}{
\begin{tikzpicture}[every node/.style={circle, draw, inner sep=2pt, minimum size=8pt}]
\node (i) at (0,0) {i};
\node (j) at ([xshift=0.974cm, yshift=-0.709cm] i) {j};
\node (k) at ([xshift=-0.372cm, yshift=-0.927cm] j) {k};
\node (l) at ([xshift=-1.204cm, yshift=0cm] k) {l};
\node (q) at ([xshift=-0.372cm, yshift=0.927cm] l) {q};
\draw (l)--(i)--(k);
\draw (q)--(i);

\end{tikzpicture}
}
& $\displaystyle\sum_{i\neq j\neq k\neq l\neq q}A_{ij}A_{jk}A_{jl}A_{jq}A_{kl}A_{kq}A_{lq}$ & $\langle (\Acal_{(2)}\odot\Acal_{(3)})\times_{1,2}(\Acal_{(2)}\odot\Acal_{(3)})\times_{1,2}(\Acal_{(2)}\odot(\Acal_{(3)}\odot\Dcal_{(2,3)}-3\Ccal_{(3)})),\Ical_{n\times n\times n}\rangle$\\

\midrule
% % --- Row 4 ---
\adjustbox{valign=m}{
\begin{tikzpicture}[every node/.style={circle, draw, inner sep=2pt, minimum size=8pt}]
\node (i) at (0,0) {i};
\node (j) at ([xshift=0.974cm, yshift=-0.709cm] i) {j};
\node (k) at ([xshift=-0.372cm, yshift=-0.927cm] j) {k};
\node (l) at ([xshift=-1.204cm, yshift=0cm] k) {l};
\node (q) at ([xshift=-0.372cm, yshift=0.927cm] l) {q};
\draw (j)--(i)--(q)--(j)--(l)--(q)--(k)--(j);

\end{tikzpicture}
}
&
\adjustbox{valign=m}{
\begin{tikzpicture}[every node/.style={circle, draw, inner sep=2pt, minimum size=8pt}]
\node (i) at (0,0) {i};
\node (j) at ([xshift=0.974cm, yshift=-0.709cm] i) {j};
\node (k) at ([xshift=-0.372cm, yshift=-0.927cm] j) {k};
\node (l) at ([xshift=-1.204cm, yshift=0cm] k) {l};
\node (q) at ([xshift=-0.372cm, yshift=0.927cm] l) {q};
\draw (i)--(k)--(l)--(i);

\end{tikzpicture}
}
& $\displaystyle\sum_{i\neq j\neq k\neq l\neq q}A_{ij}A_{jk}A_{jl}A_{jq}A_{kq}A_{lq}A_{qi}$ & $\onebm^{\top} \Big[(A^2)^{\circ3}\odot A-3A^2\odot A\odot (A^{\circ2})^2\Big]\onebm+2\tr((A^{\circ3})^2A)$\\

\midrule
% --- Row 5 ---
\adjustbox{valign=m}{
\begin{tikzpicture}[every node/.style={circle, draw, inner sep=2pt, minimum size=8pt}]
\node (i) at (0,0) {i};
\node (j) at ([xshift=0.974cm, yshift=-0.709cm] i) {j};
\node (k) at ([xshift=-0.372cm, yshift=-0.927cm] j) {k};
\node (l) at ([xshift=-1.204cm, yshift=0cm] k) {l};
\node (q) at ([xshift=-0.372cm, yshift=0.927cm] l) {q};
\draw (i)--(j)--(k)--(l)--(q)--(j)--(l);
\draw (i)--(q)--(k);

\end{tikzpicture}
}
&
\adjustbox{valign=m}{
\begin{tikzpicture}[every node/.style={circle, draw, inner sep=2pt, minimum size=8pt}]
\node (i) at (0,0) {i};
\node (j) at ([xshift=0.974cm, yshift=-0.709cm] i) {j};
\node (k) at ([xshift=-0.372cm, yshift=-0.927cm] j) {k};
\node (l) at ([xshift=-1.204cm, yshift=0cm] k) {l};
\node (q) at ([xshift=-0.372cm, yshift=0.927cm] l) {q};
\draw (k)--(i)--(l);

\end{tikzpicture}
}
& $\displaystyle\sum_{i\neq j\neq k\neq l\neq q}A_{ij}A_{jk}A_{jl}A_{jq}A_{kl}A_{kq}A_{lq}A_{qi}$ & $\langle (\Acal_{(2)}\odot\Acal_{(3)})\times_{1,2}\Big[(\Acal_{(2)}\odot\Acal_{(3)})\times_{1,2}(\Acal_{(2)}\odot\Ecal_{(3)})-2(\Acal_{(2)}\odot\Ccal_{(3)})\times_{1,2}(\Acal_{(2)}\odot\Ccal_{(3)}) \Big],\Ical_{n\times n\times n}\rangle$\\

\midrule
% --- Row 6 ---
\adjustbox{valign=m}{
\begin{tikzpicture}[every node/.style={circle, draw, inner sep=2pt, minimum size=8pt}]
\node (i) at (0,0) {i};
\node (j) at ([xshift=0.974cm, yshift=-0.709cm] i) {j};
\node (k) at ([xshift=-0.372cm, yshift=-0.927cm] j) {k};
\node (l) at ([xshift=-1.204cm, yshift=0cm] k) {l};
\node (q) at ([xshift=-0.372cm, yshift=0.927cm] l) {q};
\draw (i)--(k)--(q)--(j)--(l)--(i);
\draw (k)--(j)--(i)--(q)--(l);

\end{tikzpicture}
}
&
\adjustbox{valign=m}{
\begin{tikzpicture}[every node/.style={circle, draw, inner sep=2pt, minimum size=8pt}]
\node (i) at (0,0) {i};
\node (j) at ([xshift=0.974cm, yshift=-0.709cm] i) {j};
\node (k) at ([xshift=-0.372cm, yshift=-0.927cm] j) {k};
\node (l) at ([xshift=-1.204cm, yshift=0cm] k) {l};
\node (q) at ([xshift=-0.372cm, yshift=0.927cm] l) {q};
\draw (l)--(k);

\end{tikzpicture}
}
& $\displaystyle\sum_{i\neq j\neq k\neq l\neq q}A_{ij}A_{ik}A_{il}A_{iq}A_{jk}A_{jl}A_{jq}A_{kq}A_{lq}$ & $\langle (\mathscr{A}_{(3,4)}\odot \mathscr{A}_{(2,4)}\odot \mathscr{A}_{(2,3)}\odot \mathscr{A}_{(1,2)})\times_{1,2}(\mathscr{A}_{(3,4)}\odot \mathscr{A}_{(2,4)}\odot \mathscr{A}_{(2,3)})\times_{1,2}(\mathscr{A}_{(3,4)}\odot \mathscr{A}_{(2,3)}),\Ical_{n\times n\times n\times n}\rangle-\langle (\Ccal_{(2)}\odot\Acal_{(3)})\times_{1,2}(\Ccal_{(2)}\odot\Acal_{(3)})\times_{1,2}(\Ccal_{(2)}\odot\Acal_{(3)}),\Ical_{n\times n\times n}\rangle$\\

\midrule
% --- Row 7 ---
\adjustbox{valign=m}{
\begin{tikzpicture}[every node/.style={circle, draw, inner sep=2pt, minimum size=8pt}]
\node (i) at (0,0) {i};
\node (j) at ([xshift=0.974cm, yshift=-0.709cm] i) {j};
\node (k) at ([xshift=-0.372cm, yshift=-0.927cm] j) {k};
\node (l) at ([xshift=-1.204cm, yshift=0cm] k) {l};
\node (q) at ([xshift=-0.372cm, yshift=0.927cm] l) {q};
\draw (i)--(k)--(q)--(j)--(l)--(i);
\draw (k)--(j)--(i)--(q)--(l)--(k);

\end{tikzpicture}
}
&
\adjustbox{valign=m}{
\begin{tikzpicture}[every node/.style={circle, draw, inner sep=2pt, minimum size=8pt}]
\node (i) at (0,0) {i};
\node (j) at ([xshift=0.974cm, yshift=-0.709cm] i) {j};
\node (k) at ([xshift=-0.372cm, yshift=-0.927cm] j) {k};
\node (l) at ([xshift=-1.204cm, yshift=0cm] k) {l};
\node (q) at ([xshift=-0.372cm, yshift=0.927cm] l) {q};

\end{tikzpicture}
}
& $\displaystyle\sum_{i\neq j\neq k\neq l\neq q}A_{ij}A_{ik}A_{il}A_{iq}A_{jk}A_{jl}A_{jq}A_{kl}A_{kq}A_{lq}$ & $\langle (\mathscr{A}_{(3,4)}\odot \mathscr{A}_{(2,4)}\odot \mathscr{A}_{(2,3)}\odot \mathscr{A}_{(1,2)})\times_{1,2}(\mathscr{A}_{(3,4)}\odot \mathscr{A}_{(2,4)}\odot \mathscr{A}_{(2,3)})\times_{1,2}(\mathscr{A}_{(3,4)}\odot \mathscr{A}_{(2,4)}\odot \mathscr{A}_{(2,3)}),\Ical_{n\times n\times n\times n}\rangle$\\

\bottomrule
\end{tabular}
}
\caption{All connected 5-node subgraphs and their counting formulas given in an efficient form, Part III.}
\label{tab:result_5_nodes-3}
\end{table}

\begin{theorem}\label{prop:5-node}
Let $\Acal_{(2)}, \Acal_{(3)},\Bcal_{(2)},\Bcal_{(3)},\Ccal_{(2)},\Ccal_{(3)},\Dcal_{(2,3)},\Ecal_{(3)}\in\Rbb^{n\times n\times n}$ be 3-way tensors satisfying $[\Acal_{(2)}]_{:k':}=[\Acal_{(3)}]_{::k'}=A, [\Bcal_{(2)}]_{:k':}=[\Bcal_{(3)}]_{::k'}=A^2, [\Ccal_{(2)}]_{:k':}=[\Ccal_{(3)}]_{::k'}=A^{\circ2},[\Ecal_{(3)}]_{::k'}=A\odot A^2$ for all $1\leq k'\leq n$ and that $[\Dcal_{(2,3)}]_{:k'l'}=A\onebm$ for all $1\leq k',l'\leq n$. Define 4-way tensors $\mathscr{A}_{(1,2)},\mathscr{A}_{(2,3)},\mathscr{A}_{(2,4)},\mathscr{A}_{(3,4)}\in\Rbb^{n\times n\times n\times n}$ satisfying  $[\mathscr{A}_{(1,2)}]_{k'l'::}=[\mathscr{A}_{(2,3)}]_{:k'l':}=[\mathscr{A}_{(2,4)}]_{:k':l'}=[\mathscr{A}_{(3,4)}]_{::k'l'}=A$ for all $1\leq k',l'\leq n$. We further define $\Ical_{n\times n\times n}\in\Rbb^{n\times n\times n},\Ical_{n\times n\times n\times n}\in\Rbb^{n\times n\times n\times n}$ satisfying $[\Ical_{n\times n\times n}]_{::k'}=I_n,[\Ical_{n\times n\times n\times n}]_{::k'l'}=I_n$.

The following identities hold:
\begin{enumerate}[label=(\alph*)]
    \item $\displaystyle\sum_{i\neq j\neq k\neq l\neq q}A_{ij}A_{jk}A_{kl}A_{lq}=\onebm^\top \Big[A^4-2A^2A^{\circ2}-A\diag(A^{\circ2}\onebm)A+3(A^{\circ2})^2+2A^{\circ3}A-2A^{\circ4}\Big]\onebm-2\diag(A^3)^\top A\onebm-\tr(A^4)+3\tr(A^{\circ2}A^2)$,
    \item $\displaystyle\sum_{i\neq j\neq k\neq l\neq q}A_{ij}A_{jk}A_{kl}A_{lq}A_{qi} = \tr(A^5)-5\onebm^\top A^{\circ2}\diag(A^3)+5\tr(A^{\circ3}A^2)$,
    \item $\displaystyle\sum_{i\neq j\neq k\neq l\neq q}A_{ij}A_{jk}A_{il}A_{qi}=\onebm^\top \Big[ A\diag(A^{2}\onebm)A-A\diag(A^{\circ2}\onebm)A-2AA^{\circ2}A-A^{\circ2}A^2+4A^{\circ3}A+(A^{\circ2})^2-2A^{\circ4}\Big]\onebm-2\diag(A^3)^\top A\onebm+4\tr(A^{\circ2}A^2)$,
    \item $\displaystyle\sum_{i\neq j\neq k\neq l\neq q}A_{ij}A_{jk}A_{jq}A_{kl}A_{lq}=\diag(A^4)^\top A \onebm+\onebm^\top \Big[ 4A^{\circ3}A^{\circ2}+A^{\circ4}A-(A^2)^{\circ2}\odot A$
    
    $-A\diag(A^{\circ2}\onebm)A^{\circ2}-A(A^{\circ2})^2-2A^{\circ5}\Big]\onebm+\tr(A(A^{\circ2})^2)-2\tr(A^{\circ2}A^3)$,
    \item $\displaystyle\sum_{i\neq j\neq k\neq l\neq q}A_{ij}A_{jk}A_{jl}A_{kq}A_{lq}A_{qi}=\onebm^\top \Big[ (A^2)^{\circ3}+3A^{\circ4}A^{\circ2})+2(A^{\circ3})^2-2A^{\circ6}\Big]\onebm-\onebm^\top(A^{\circ2}\onebm)^{\circ3}-3\tr((A^{\circ2})^2A^2)$,
    \item $\displaystyle\sum_{i\neq j\neq k\neq l\neq q}A_{ij}A_{jk}A_{jq}A_{kl}A_{qi}=\diag(A^3)^\top A^2\onebm-2\diag(A^{\circ2}A^2)^\top A\onebm-\diag(A^3)^\top A^{\circ2}\onebm$
    
    $-2\onebm^\top ((A^2)^{\circ2}\odot A)\onebm+2\tr(A^{\circ3}A^2)+2\tr((A^{\circ2})^2A)$,
    \item $\displaystyle\sum_{i\neq j\neq k\neq l\neq q}A_{ij}A_{jk}A_{jq}A_{lq}A_{qi}=\onebm^\top \Big[ A(A^2\odot A)A-(A^2)^{\circ2}\odot A\Big]\onebm{-2\diag(A^{\circ2}A^2)^\top A\onebm-2\diag(AA^{\circ2}A)^\top A\onebm}$
    
    $+4\tr((A^{\circ2})^2A)+\tr(A^{\circ3}A^2)$,
    \item $\displaystyle\sum_{i\neq j\neq k\neq l\neq q}A_{ij}A_{jk}A_{jq}A_{kl}A_{lq}A_{qi}=\onebm^\top \Big[ A\odot A^2\odot A^3-2A\odot A^2\odot (A^{\circ2}A)\Big]\onebm+2\tr(AA^{\circ2}A^{\circ3})+\tr(A^{\circ4}A^2)-2\diag(A^{\circ2}A^2)^\top A^{\circ2}\onebm$,
    \item $\displaystyle\sum_{i\neq j\neq k\neq l\neq q}A_{jk}A_{jl}A_{jq}A_{kl}A_{lq}A_{qi}=\diag(A(A\odot A^2)A)^\top A\onebm-\diag(A^{\circ2}AA^{\circ2})^\top A\onebm -2\tr(AA^{\circ2}(A\odot A^2))+2\tr(AA^{\circ2}A^{\circ3})-\langle (\Acal_{(2)}\odot \Acal_{(3)})\times_{1,2}(\Acal_{(2)}\odot \Acal_{(3)})\times_{1,2}(\Acal_{(2)}\odot \Acal_{(3)}),\Ical_{n\times n\times n}\rangle$,
    \item $\displaystyle\sum_{i\neq j\neq k\neq l\neq q}A_{ij}A_{jk}A_{jl}A_{kl}A_{kq}A_{lq}A_{qi}=\langle (\Acal_{(2)}\odot\Acal_{(3)})\times_{1,2}(\Acal_{(2)}\odot\Acal_{(3)})\times_{1,2}(\Bcal_{(2)}\odot\Acal_{(3)}),\Ical_{n\times n\times n}\rangle-2\tr([(A^{\circ2}A)\odot A]AA^{\circ2})-\diag(A^{\circ2}AA^{\circ2})^\top A^{\circ2}\onebm+2\tr(AA^{\circ2}A^{\circ4})$,
    \item $\displaystyle\sum_{i\neq j\neq k\neq l\neq q}A_{ij}A_{ik}A_{il}A_{qi}=\onebm^\top(A\onebm)^{\circ4}-6\onebm^\top A^{\circ2}(A\onebm)^{\circ2}+\onebm^\top \Big[ 3(A^{\circ2})^2+8A^{\circ3}A- 6A^{\circ4} \Big]\onebm$,
    \item $\displaystyle\sum_{i\neq j\neq k\neq l\neq q}A_{ij}A_{ik}A_{jk}A_{il}A_{qi}=\onebm^\top (A\odot A^2) (A\onebm)^{\circ2}-4\diag(A^2A^{\circ2})^\top A\onebm-\diag(A^3)^\top A^{\circ2}\onebm+2\tr((A^{\circ2})^2A)+4\tr(A^{\circ3}A^2)$,
    \item $\displaystyle\sum_{i\neq j\neq k\neq l\neq q}A_{ij}A_{ik}A_{jk}A_{il}A_{lq}A_{qi}=\tr((A^3)^{\circ2})-4\onebm^\top ((A^2)^{\circ2}\odot A^{\circ2})\onebm+2\tr((A^{\circ2})^3)$,
    \item $\displaystyle\sum_{i\neq j\neq k\neq l\neq q}A_{ij}A_{jk}A_{jl}A_{jq}A_{kq}A_{lq}=\onebm^\top \Big[A((A^2)^{\circ2}\odot A)-2(A^{\circ2}A)\odot A^2\odot A-A^{\circ2}\odot (A^2)^{\circ2}\Big]\onebm-\onebm^\top A\diag ((A^{\circ2})^2A)+\tr\Big[(A^{\circ2})^3+2A^{\circ3}A^{\circ2}A\Big]$,
    \item $\displaystyle\sum_{i\neq j\neq k\neq l\neq q}A_{ij}A_{jl}A_{jq}A_{kl}A_{kq}A_{lq}A_{qi}=\tr\Big[A(A\odot A^2)^2-AA^{\circ2}(A^{\circ2}\odot A^2)-A^2((A^{\circ2}A)\odot A^{\circ2})+A^{\circ3}(A^{\circ2})^2\Big]-\langle (\Acal_{(2)}\odot\Acal_{(3)})\times_{1,2}(\Acal_{(2)}\odot\Acal_{(3)})\times_{1,2}(\Acal_{(2)}\odot\Ccal_{(3)}),\Ical_{n\times n\times n}\rangle $,
    \item $\displaystyle\sum_{i\neq j\neq k\neq l\neq q}A_{ij}A_{jk}A_{jl}A_{ki}A_{kq}A_{lq}A_{li}A_{qi}=\langle (\mathscr{A}_{(3,4)}\odot \mathscr{A}_{(2,4)}\odot \mathscr{A}_{(2,3)}\odot \mathscr{A}_{(1,2)})\times_{1,2}(\mathscr{A}_{(3,4)}\odot \mathscr{A}_{(2,4)})\times_{1,2}(\mathscr{A}_{(3,4)}\odot \mathscr{A}_{(2,3)}),\Ical_{n\times n\times n\times n}\rangle+\tr(A^{\circ4}(A^{\circ2})^2)-2\onebm^\top((AA^{\circ2})^{\circ2}\odot A^{\circ2})\onebm$,
    \item $\displaystyle\sum_{i\neq j\neq k\neq l\neq q}A_{ij}A_{jk}A_{jl}A_{jq}A_{kl}A_{kq}A_{lq}=\langle (\Acal_{(2)}\odot\Acal_{(3)})\times_{1,2}(\Acal_{(2)}\odot\Acal_{(3)})\times_{1,2}(\Acal_{(2)}\odot(\Acal_{(3)}\odot\Dcal_{(2,3)}-3\Ccal_{(3)})),\Ical_{n\times n\times n}\rangle$,
    \item $\displaystyle\sum_{i\neq j\neq k\neq l\neq q}A_{ij}A_{jk}A_{jl}A_{jq}A_{kq}A_{lq}A_{qi}=\onebm^{\top} \Big[(A^2)^{\circ3}\odot A-3A^2\odot A\odot (A^{\circ2})^2\Big]\onebm+2\tr((A^{\circ3})^2A)$,
    \item $\displaystyle\sum_{i\neq j\neq k\neq l\neq q}A_{ij}A_{jk}A_{jl}A_{jq}A_{kl}A_{kq}A_{lq}A_{qi}=\langle (\Acal_{(2)}\odot\Acal_{(3)})\times_{1,2}\Big[(\Acal_{(2)}\odot\Acal_{(3)})\times_{1,2}(\Acal_{(2)}\odot\Ecal_{(3)})-2(\Acal_{(2)}\odot\Ccal_{(3)})\times_{1,2}(\Acal_{(2)}\odot\Ccal_{(3)}) \Big],\Ical_{n\times n\times n}\rangle$,
    \item $\displaystyle\sum_{i\neq j\neq k\neq l\neq q}A_{ij}A_{ik}A_{il}A_{iq}A_{jk}A_{jl}A_{jq}A_{kq}A_{lq}=\langle (\mathscr{A}_{(3,4)}\odot \mathscr{A}_{(2,4)}\odot \mathscr{A}_{(2,3)}\odot \mathscr{A}_{(1,2)})\times_{1,2}(\mathscr{A}_{(3,4)}\odot \mathscr{A}_{(2,4)}\odot \mathscr{A}_{(2,3)})\times_{1,2}(\mathscr{A}_{(3,4)}\odot \mathscr{A}_{(2,3)}),\Ical_{n\times n\times n\times n}\rangle-\langle (\Ccal_{(2)}\odot\Acal_{(3)})\times_{1,2}(\Ccal_{(2)}\odot\Acal_{(3)})\times_{1,2}(\Ccal_{(2)}\odot\Acal_{(3)}),\Ical_{n\times n\times n}\rangle$,
    \item $\displaystyle\sum_{i\neq j\neq k\neq l\neq q}A_{ij}A_{ik}A_{il}A_{iq}A_{jk}A_{jl}A_{jq}A_{kl}A_{kq}A_{lq}=\langle (\mathscr{A}_{(3,4)}\odot \mathscr{A}_{(2,4)}\odot \mathscr{A}_{(2,3)}\odot \mathscr{A}_{(1,2)})\times_{1,2}(\mathscr{A}_{(3,4)}\odot \mathscr{A}_{(2,4)}\odot \mathscr{A}_{(2,3)})\times_{1,2}(\mathscr{A}_{(3,4)}\odot \mathscr{A}_{(2,4)}\odot \mathscr{A}_{(2,3)}),\Ical_{n\times n\times n\times n}\rangle$.
\end{enumerate}
\end{theorem}

\begin{proof}
Suppose $V=\{i,j,k,l,q\}$ be the set of nodes.
\begin{enumerate}[label=(\alph*)]
\item Let $F = A_{ij}A_{jk}A_{kl}A_{lq}$. Consider the graph $G = (V, E_G)$ with $E_G = \{i - j - k - l - q\}$ and its complement graph $H = (V, E_H)$ with $E_H = \{i - q - k - i - l - j - q\}$ (see Row 1 in \Cref{tab:result_5_nodes-1}).  

The maximum clique contained in $H$ is the 3-clique $i - q - k - i$. Hence, the nonzero contractions are $F_{2+1+1+1}^{(5)}, F_{2+2+1}^{(5)},F_{3+1+1}^{(5)}$, and $F_{3+2}^{(5)}$.  

First, we consider the contraction $F_{2+1+1+1}^{(5)}$, which has six components corresponding to the six edges in $E_H$. These components can be obtained directly by setting the indices of each edge equal, one by one, or alternatively, by reducing the graph $G$. For example, since $i - q \in E_H$, setting $q = i$ reduces $G$ to $G' = (V', E_G')$ with $V' = \{i,j,k,l\}$ and $E_G' = \{i - j - k - l - i\}$, corresponding to $F(i,j,k,l,i) = A_{ij}A_{jk}A_{kl}A_{li}$. All the reduced graphs are shown in \Cref{fig:reduced_graph_1}.  

Note that edges in a reduced graph can appear multiple times. For instance, if $q = k$, then $F(i,j,k,l,k) = A_{ij}A_{jk}A_{kl}^2$, where the exponent 2 indicates that the edge $k - l$ appears twice in the reduced graph.  

Furthermore, the reduced graphs for the cases $q = k$ and $k = i$, as well as $i = l$ and $j = q$, are of the same type; hence, the associated functions are equivalent. Summing the associated functions of these reduced graphs yields
\begin{align}\label{eq:5-a-F2111}
    F_{2+1+1+1}^{(5)} & = A_{ij}A_{jk}A_{kl}A_{li}+2A_{ij}A_{jk}A_{kl}^2+2A_{ij}A_{jk}A_{kl}A_{jl}+A_{ij}A_{jk}^2A_{jl}, \nonumber \\
    \sum_{i,j,k,l}F_{2+1+1+1}^{(5)} & = \sum_{i,j,k,l}\left(A_{ij}A_{jk}A_{kl}A_{li}+2A_{ij}A_{jk}A_{kl}^2+2A_{ij}A_{jk}A_{kl}A_{jl}+A_{ij}A_{jk}^2A_{jl}\right)\nonumber \\
    & = \tr(A^4)+2\onebm^\top A^2A^{\circ2}\onebm+2\diag(A^3)^\top A\onebm +\onebm^\top A \diag(A^{\circ2}\onebm) A\onebm.
\end{align}
\begin{figure}[!ht]
    \centering
    \subfloat[$i=q$]{%
        \begin{tikzpicture}[every node/.style={circle, draw, inner sep=0pt, minimum size=12pt,text centered},node distance=0.7cm]
            \node (l) at (0,2) {l};
            \node (i) [right=of l] {i};
            \node (j) [below=of i] {j};
            \node (k) [below=of l] {k};
            \draw (i)--(j)--(k)--(l)--(i);
        \end{tikzpicture}
    }
    \hspace{0.8cm} % space between subfigures
    \subfloat[$q=k$]{%
        \begin{tikzpicture}[every node/.style={circle, draw, inner sep=0pt, minimum size=12pt,text centered},node distance=0.7cm]
            \node (l) at (0,2) {l};
            \node (i) [right=of l] {i};
            \node (j) [below=of i] {j};
            \node (k) [below=of l] {k};
                % Single edges
            \draw (i)--(j)--(k);
        
            % Double edge between i and j
            \draw (k) to[bend left=15] (l);
            \draw (k) to[bend right=15] (l);
        \end{tikzpicture}
    }
    \hspace{0.8cm} % space between subfigures
    \subfloat[$k=i$]{%
        \begin{tikzpicture}[every node/.style={circle, draw, inner sep=0pt, minimum size=12pt,text centered},node distance=0.7cm]
            \node (q) at (0,2) {q};
            \node (i) [right=of l] {i};
            \node (j) [below=of i] {j};
            \node (l) [below=of l] {l};
            \draw (q)--(l)--(i);
            \draw (i) to[bend left=15] (j);
            \draw (i) to[bend right=15] (j);
        \end{tikzpicture}
    }
    \hspace{0.8cm} % space between subfigures
    \subfloat[$i=l$]{%
        \begin{tikzpicture}[every node/.style={circle, draw, inner sep=0pt, minimum size=12pt,text centered},node distance=0.7cm]
            \node (q) at (0,2) {q};
            \node (i) [right=of q] {i};
            \node (j) [below=of i] {j};
            \node (k) [below=of q] {k};
            \draw (q)--(i)--(j)--(k)--(i);
        \end{tikzpicture}
    }
    \hspace{0.8cm} % space between subfigures
    \subfloat[$l=j$]{%
        \begin{tikzpicture}[every node/.style={circle, draw, inner sep=0pt, minimum size=12pt,text centered},node distance=0.7cm]
            \node (q) at (0,2) {q};
            \node (i) [right=of q] {i};
            \node (j) [below=of i] {j};
            \node (k) [below=of q] {k};
            \draw (q)--(j)--(i);
            \draw (k) to[bend left=15] (j);
            \draw (k) to[bend right=15] (j);
        \end{tikzpicture}
    }
    \hspace{0.8cm} % space between subfigures
    \subfloat[$j=q$]{%
        \begin{tikzpicture}[every node/.style={circle, draw, inner sep=0pt, minimum size=12pt,text centered},node distance=0.7cm]
            \node (l) at (0,2) {l};
            \node (i) [right=of l] {i};
            \node (j) [below=of i] {j};
            \node (k) [below=of l] {k};
            \draw (i)--(j)--(k)--(l)--(j);
        \end{tikzpicture}
    }
    \caption{Reduced graphs of ($G$, Row 1 \Cref{tab:result_5_nodes-1}) by letting nodes of the edge in $E_H$ equal.}
    \label{fig:reduced_graph_1}
\end{figure}
Next, we consider the contraction $F_{2+2+1}^{(5)}$. There are six disjoint subgraphs contained in $H$, given by $E^1_H=\{i-q,j-l\}, E^2_H=\{i-l,q-k\},E^3_H=\{i-l,q-j\},E^4_H=\{i-k,q-j\},E^5_H=\{i-k,j-l\},E^6_H=\{j-l,q-k\}$. Using the same arguments as above, we then have that
\begin{align}\label{eq:5-a-F221}
    F_{2+2+1}^{(5)} & = A_{ij}^2A_{jk}^2+3A_{ij}^2A_{jk}A_{ki}+2A_{ij}^3A_{jl}, \nonumber\\
    \sum_{i,j,k}F_{2+2+1}^{(5)} & = \sum_{i,j,k}\left(A_{ij}^2A_{jk}^2+3A_{ij}^2A_{jk}A_{ki}+2A_{ij}^3A_{jl}\right) = \onebm^\top(A^{\circ2})^2\onebm+3\tr(A^{\circ2}A^2)+2\onebm^\top A^{\circ3}A\onebm.
\end{align}
The last two nonzero contractions can be obtained by
\begin{align*}
    F_{3+2}^{(5)} &= F(i,j,i,j,i)=A_{ij}^4, \quad F_{3+1+1}^{(5)} = F(i,j,i,k,i)=A_{ij}^2A_{ik}^2
\end{align*}
and their sums over indices are
\begin{align}
    \sum_{i,j}F_{3+2}^{(5)} &=\onebm^\top A^{\circ4}\onebm,\quad \sum_{i,j,k}F_{3+1+1}^{(5)} = \onebm^\top (A^{\circ2})^2\onebm \label{eq:5-a-F311}.
\end{align}

Using the formula in \Cref{sec:summary} and combining \eqref{eq:5-a-F2111}-\eqref{eq:5-a-F311}, we have that 
\begin{align*}
    &\sum_{i\neq j\neq k\neq l\neq q}A_{ij}A_{jk}A_{kl}A_{lq} \\
    =&  \sum_{i,j,k,l,q}F-\sum_{i,j,k,l}F_{2+1+1+1}^{(5)}+\sum_{i,j,k}F_{2+2+1}^{(5)}+2\sum_{i,j,k}F_{3+1+1}^{(5)}-2\sum_{i,j}F_{3+2}^{(5)}\\
    =& \onebm^\top \Big[A^4-2A^2A^{\circ2}-A\diag(A^{\circ2}\onebm)A+3(A^{\circ2})^2+2A^{\circ3}A-2A^{\circ4}\Big]\onebm\\
    & -2\diag(A^3)^\top A\onebm-\tr(A^4)+3\tr(A^{\circ2}A^2).
\end{align*}

\item Let $F = A_{ij}A_{jk}A_{kl}A_{lq}A_{qi}$, and let the complement graph $H$ be as shown in Row 2 of \Cref{tab:result_5_nodes-1}. Clearly, there are only two nonzero, nontrivial contractions: $F_{2+1+1+1}^{(5)}$ and $F_{2+2+1}^{(5)}$.

For $F_{2+1+1+1}^{(5)}$, we observe that $H$ contains five edges, corresponding to five components. Using the same arguments in Part (a), the reduced graphs of $G$ by forcing the nodes of the edges in $H$ equal are of the same type, so these five components are equivalent. Thus, we have that 
\begin{align*}
    F_{2+1+1+1}^{(5)} &\sim 5F(i,j,i,q,l) = 5A_{ij}^2A_{il}A_{lq}A_{qi}.
\end{align*}
Similarly, we can also verify that there are five equivalent components in $F_{2+2+1}^{(5)}$. Then by setting $l=i$ and $q=j$ we find that
\begin{align*}
    F_{2+2+1}^{(5)} &\sim 5F(i,j,k,i,j) = 5A_{ij}^3A_{jk}A_{ki}.
\end{align*}
Using the formula in \Cref{sec:summary}, we have that 
\begin{align*}
    \sum_{i\neq j\neq k\neq l\neq q}A_{ij}A_{jk}A_{kl}A_{lq}A_{qi} & =  \sum_{i,j,k,l,q}A_{ij}A_{jk}A_{kl}A_{lq}A_{qi}-\sum_{i,j,l,q}5A_{ij}^2A_{il}A_{lq}A_{qi}+\sum_{i,j,k}5A_{ij}^3A_{jk}A_{ki} \\
    & = \tr(A^5)-5\sum_i\sum_jA_{ij}^2\sum_{l,q}A_{il}A_{lq}A_{qi}+5\tr(A^{\circ3}A^2) \\
    & = \tr(A^5)-5\onebm^\top A^{\circ2}\diag(A^3)+5\tr(A^{\circ3}A^2).
\end{align*}

\item Let $F = A_{ij}A_{jk}A_{il}A_{iq}$.
We first compute the general sum of $F$ as 
\begin{align}\label{eq:5-c-F}
    \sum_{i,j,k,l,q}A_{ij}A_{jk}A_{il}A_{iq} = \sum_{i}\sum_k\sum_jA_{ij}A_{jk}\sum_lA_{il}\sum_qA_{iq} = \onebm^\top A\diag(A^{2}\onebm)A\onebm.
\end{align}
Consider the graph $G = (V, E_G)$ with $E_G = \{q-i - j - k, i-l\}$ and its complement graph $H = (V, E_H)$ with $E_H = \{i-k-l-q-q,q-j-l\}$ (see Row 3 in \Cref{tab:result_5_nodes-1}). Using the same arguments in Part (a), the contraction $F_{2+1+1+1}^{(5)}$ can be obtained from the reduced graphs of $G$, shown in \Cref{fig:reduced_graph_3}, as follows:
\begin{align*}
    F_{2+1+1+1}^{(5)} & = A_{ij}^2A_{ik}A_{iq}+2A_{ij}A_{jk}A_{ki}A_{il}+2A_{ij}^2A_{jk}A_{il}+A_{ij}A_{jk}A_{il}^2,
\end{align*}
which has the sum 
\begin{align}\label{eq:5-c-F2111}
    \sum_{i,j,k,l}F_{2+1+1+1}^{(5)} = \onebm^\top A\diag(A^{\circ2}\onebm)A\onebm+2\diag(A^3)^\top A\onebm+2\onebm^\top AA^{\circ2}A\onebm+\onebm^\top A^{\circ2}A^2\onebm.
\end{align}
\begin{figure}[!ht]
    \centering
    \subfloat[$i=k$]{%
        \begin{tikzpicture}[every node/.style={circle, draw, inner sep=0pt, minimum size=12pt,text centered},node distance=0.7cm]
            \node (q) at (0,2) {q};
            \node (i) [right=of q] {i};
            \node (j) [below=of i] {j};
            \node (l) [below=of q] {l};
            \draw (q)--(i)--(l);
            \draw (i) to[bend left=15] (j);
            \draw (i) to[bend right=15] (j);
        \end{tikzpicture}
    }
    \hspace{0.8cm} % space between subfigures
    \subfloat[$j=l$]{%
        \begin{tikzpicture}[every node/.style={circle, draw, inner sep=0pt, minimum size=12pt,text centered},node distance=0.7cm]
            \node (q) at (0,2) {q};
            \node (i) [right=of q] {i};
            \node (j) [below=of i] {j};
            \node (k) [below=of q] {k};
                % Single edges
            \draw (i)--(q);
            \draw (j)--(k);
            \draw (i) to[bend left=15] (j);
            \draw (i) to[bend right=15] (j);
        \end{tikzpicture}
    }
    \hspace{0.8cm} % space between subfigures
    \subfloat[$j=q$]{%
        \begin{tikzpicture}[every node/.style={circle, draw, inner sep=0pt, minimum size=12pt,text centered},node distance=0.7cm]
            \node (l) at (0,2) {l};
            \node (i) [right=of l] {i};
            \node (j) [below=of i] {j};
            \node (k) [below=of l] {k};
            \draw (i)--(l);
            \draw (k)--(j);
            \draw (i) to[bend left=15] (j);
            \draw (i) to[bend right=15] (j);
        \end{tikzpicture}
    }
    \hspace{0.8cm} % space between subfigures
    \subfloat[$k=l$]{%
        \begin{tikzpicture}[every node/.style={circle, draw, inner sep=0pt, minimum size=12pt,text centered},node distance=0.7cm]
            \node (q) at (0,2) {q};
            \node (i) [right=of q] {i};
            \node (j) [below=of i] {j};
            \node (l) [below=of q] {k};
                % Single edges
            \draw (q)--(i)--(l)--(j)--(i);
        \end{tikzpicture}
    }
    \hspace{0.8cm} % space between subfigures
    \subfloat[$k=q$]{%
        \begin{tikzpicture}[every node/.style={circle, draw, inner sep=0pt, minimum size=12pt,text centered},node distance=0.7cm]
            \node (q) at (0,2) {l};
            \node (i) [right=of q] {i};
            \node (j) [below=of i] {j};
            \node (l) [below=of q] {k};
                % Single edges
            \draw (q)--(i)--(l)--(j)--(i);
        \end{tikzpicture}
    }
    \hspace{0.8cm} % space between subfigures
    \subfloat[$j=q$]{%
        \begin{tikzpicture}[every node/.style={circle, draw, inner sep=0pt, minimum size=12pt,text centered},node distance=0.7cm]
            \node (l) at (0,2) {l};
            \node (i) [right=of l] {i};
            \node (j) [below=of i] {j};
            \node (k) [below=of l] {k};
            % \draw (l)--(i);
            \draw (k)--(j)--(i);
            \draw (i) to[bend left=15] (l);
            \draw (i) to[bend right=15] (l);
        \end{tikzpicture}
    }
    \caption{Reduced graphs of ($G$, Row 3 \Cref{tab:result_5_nodes-1}) by letting nodes of the edge in $E_H$ equal.}
    \label{fig:reduced_graph_3}
\end{figure}

For $F_{2+2+1}^{(5)}$, there are five disjoint subgraphs contained in $H$, given by $E_H^1=\{i-k,j-l\},E_H^2=\{i-k,j-q\},E_H^3=\{i-k,q-l\},E_H^4=\{j-l,k-q\},E_H^5=\{j-q,k-l\}$, which lead to 
\begin{align*}
    F_{2+2+1}^{(5)} & = 2A_{ij}^3A_{ik}+A_{ij}^2A_{ik}^2+2A_{ij}^2A_{jk}A_{ki}
\end{align*}
and
\begin{align}
    \sum_{i,j,k}F_{2+2+1}^{(5)} & = 2\onebm^\top AA^{\circ3}\onebm+\onebm^\top (A^{\circ2})^2\onebm+2\tr(A^{\circ2}A^2)
\end{align}
For $F_{3+1+1}^{(5)}$, there are two disjoint subgraphs, $E^1_H=\{q-k-l-q\}$ and $E^2_H=\{q-j-l-q\}$ so that 
\begin{align*}
    F_{3+1+1}^{(5)} & = A_{ij}^3A_{jk}+A_{ij}A_{jk}A_{ki}^2
\end{align*}
and that
\begin{align}
    \sum_{i,j,k}F_{3+1+1}^{(5)} & = \onebm^\top A^{\circ3}A\onebm+\tr(A^{\circ2}A^2).
\end{align}
The last nonzero contraction $F_{3+2}^{(5)}$ is obtained by setting $q=l=j$ and $i=k$ that $F_{3+2}^{(5)} = A_{ij}^4$ and
\begin{align}\label{eq:5-c-F32}
    \sum_{i,j}F_{3+2}^{(5)} = \onebm^\top A^{\circ4}\onebm.
\end{align}
Using the formula in \Cref{sec:summary} and combining \eqref{eq:5-c-F}-\eqref{eq:5-c-F32}, we have that 
\begin{align*}
    &\sum_{i\neq j\neq k\neq l\neq q}A_{ij}A_{jk}A_{il}A_{iq} \\
    =& \sum_{i,j,k,l,q}F-\sum_{i,j,k,l}F_{2+1+1+1}^{(5)}+\sum_{i,j,k}F_{2+2+1}^{(5)}+2\sum_{i,j,k}F_{3+1+1}^{(5)}-2\sum_{i,j}F_{3+2}^{(5)}\\
    = & \onebm^\top \Big[ A\diag(A^{2}\onebm)A-A\diag(A^{\circ2}\onebm)A-2AA^{\circ2}A-A^{\circ2}A^2+4A^{\circ3}A+(A^{\circ2})^2-2A^{\circ4}\Big]\onebm \\
    &-2\diag(A^3)^\top A\onebm+4\tr(A^{\circ2}A^2).
\end{align*}

\item Let $F = A_{ij}A_{jk}A_{jq}A_{kl}A_{lq}$.
We first compute the general sum of $F$ as 
\begin{align}\label{eq:5-d-F}
    \sum_{i,j,k,l,q}A_{ij}A_{ij}A_{jk}A_{jq}A_{kl}A_{lq} = \sum_{i,j}A_{ij}\sum_{k,l,q}A_{ij}A_{jk}A_{jk}A_{kl}A_{lj}=\sum_{i,j}A_{ij}[A^4]_{jj}=\diag(A^4)^\top A\onebm.
\end{align}
Consider the graph $G = (V, E_G)$ with $E_G = \{i-j-k-l-q-j\}$ and its complement graph $H = (V, E_H)$ with $E_H = \{j-l-i-k-q-i\}$ (see Row 4 in \Cref{tab:result_5_nodes-1}). Using the same arguments in Part (a), the contraction $F_{2+1+1+1}^{(5)}$ can be obtained from the reduced graphs of $G$, shown in \Cref{fig:reduced_graph_4}, as follows:
\begin{align*}
    F_{2+1+1+1}^{(5)} & = 2A_{ij}^2A_{jk}A_{kl}A_{li}+A_{ij}A_{jk}A_{ki}A_{kl}A_{li}+A_{ij}A_{jk}^2A_{jl}^2+A_{ij}A_{jk}^2A_{kl}^2,
\end{align*}
which has the sum 
\begin{align}
    \sum_{i,j,k,l}F_{2+1+1+1}^{(5)} = 2\tr(A^{\circ2}A^3)+\onebm^\top ((A^2)^{\circ2}\odot A)\onebm+\onebm^\top A \diag(A^{\circ2}\onebm)A^{\circ2}\onebm+\onebm^\top A(A^{\circ2})^2\onebm.
\end{align}

\begin{figure}[!ht]
    \centering
    \subfloat[$i=k$]{%
        \begin{tikzpicture}[every node/.style={circle, draw, inner sep=0pt, minimum size=12pt,text centered},node distance=0.7cm]
            \node (l) at (0,2) {l};
            \node (i) [right=of l] {i};
            \node (j) [below=of i] {j};
            \node (q) [below=of l] {q};
            \draw (i)--(l)--(q)--(j);
            \draw (i) to[bend left=15] (j);
            \draw (i) to[bend right=15] (j);
        \end{tikzpicture}
    }
    \hspace{0.8cm} % space between subfigures
    \subfloat[$i=l$]{%
        \begin{tikzpicture}[every node/.style={circle, draw, inner sep=0pt, minimum size=12pt,text centered},node distance=0.7cm]
            \node (q) at (0,2) {q};
            \node (j) [right=of q] {j};
            \node (k) [below=of j] {k};
            \node (l) [below=of q] {l};
                % Single edges
            \draw (i)--(k)--(l)--(q)--(j)--(l);
            
        \end{tikzpicture}
    }
    \hspace{0.8cm} % space between subfigures
    \subfloat[$i=q$]{%
        \begin{tikzpicture}[every node/.style={circle, draw, inner sep=0pt, minimum size=12pt,text centered},node distance=0.7cm]
            \node (j) at (0,2) {j};
            \node (k) [below=of j] {k};
            \node (l) [left=of k] {l};
            \node (q) [above=of l] {q};
            \draw (j)--(k)--(l)--(q);
            \draw (q) to[bend left=15] (j);
            \draw (q) to[bend right=15] (j);
        \end{tikzpicture}
    }
    \hspace{0.8cm} % space between subfigures
    \subfloat[$j=l$]{%
        \begin{tikzpicture}[every node/.style={circle, draw, inner sep=0pt, minimum size=12pt,text centered},node distance=0.7cm]
            \node (i) at (0,2) {i};
            \node (j) [below=of i] {j};
            \node (k) [left=of j] {k};
            \node (q) [above=of k] {q};
            \draw (i)--(j);
            \draw (k) to[bend left=15] (j);
            \draw (k) to[bend right=15] (j);
            \draw (q) to[bend left=15] (j);
            \draw (q) to[bend right=15] (j);
        \end{tikzpicture}
    }
    \hspace{0.8cm} % space between subfigures
    \subfloat[$k=q$]{%
        \begin{tikzpicture}[every node/.style={circle, draw, inner sep=0pt, minimum size=12pt,text centered},node distance=0.7cm]
            \node (i) at (0,2) {i};
            \node (j) [below=of i] {j};
            \node (k) [left=of j] {k};
            \node (l) [above=of k] {l};
            \draw (i)--(j);
            \draw (k) to[bend left=15] (j);
            \draw (k) to[bend right=15] (j);
            \draw (k) to[bend left=15] (l);
            \draw (k) to[bend right=15] (l);
            
        \end{tikzpicture}
    }
    \caption{Reduced graphs of ($G$, Row 4 \Cref{tab:result_5_nodes-1}) by letting nodes of the edge in $E_H$ equal.}
    \label{fig:reduced_graph_4}
\end{figure}
For $F_{2+2+1}^{(5)}$, there are four disjoint subgraphs contained in $H$, given by $E_H^1=\{i-k,j-l\},E_H^2=\{i-l,k-q\},E_H^3=\{i-q,j-l\},E_H^4=\{j-l,k-q\}$, which lead to 
\begin{align*}
    F_{2+2+1}^{(5)} & = 2A_{ij}^3A_{jk}^2+A_{ij}A_{jk}^2A_{ki}^2+A_{ij}A_{jk}^4
\end{align*}
and
\begin{align}
    \sum_{i,j,k}F_{2+2+1}^{(5)} & = 2\onebm^\top A^{\circ3}A^{\circ2}\onebm+\tr(A(A^{\circ2})^2)+\onebm^\top A^{\circ4}A\onebm.
\end{align}

For $F_{3+1+1}^{(5)}$, there is one 3-clique in $H$, $E^1_H=\{q-k-i-q\}$ so that 
\begin{align}
    F_{3+1+1}^{(5)} & = A_{ij}^3A_{ik}^2, \quad \sum_{i,j,k}F_{3+1+1}^{(5)} = \onebm^\top A^{\circ3}A^{\circ2}\onebm.
\end{align}
The last nonzero contraction $F_{3+2}^{(5)}$ is obtained by setting $q=k=i$ and $j=l$ that $F_{3+2}^{(5)} = A_{ij}^5$ and
\begin{align}\label{eq:5-d-F32}
    \sum_{i,j}F_{3+2}^{(5)} = \onebm^\top A^{\circ5}\onebm.
\end{align}
Using the formula in \Cref{sec:summary} and combining \eqref{eq:5-d-F}-\eqref{eq:5-d-F32}, we have that 
\begin{align*}
    &\sum_{i\neq j\neq k\neq l\neq q}A_{ij}A_{ij}A_{jk}A_{jq}A_{kl}A_{lq}\\
    =& \sum_{i,j,k,l,q}F-\sum_{i,j,k,l}F_{2+1+1+1}^{(5)}+\sum_{i,j,k}F_{2+2+1}^{(5)}+2\sum_{i,j,k}F_{3+1+1}^{(5)}-2\sum_{i,j}F_{3+2}^{(5)}\\
    = & \diag(A^4)^\top A \onebm+\onebm^\top \Big[ 4A^{\circ3}A^{\circ2}+A^{\circ4}A-(A^2)^{\circ2}\odot A-A\diag(A^{\circ2}\onebm)A^{\circ2}-A(A^{\circ2})^2-2A^{\circ5}\Big]\onebm \\
    & +\tr(A(A^{\circ2})^2)-2\tr(A^{\circ2}A^3).
\end{align*}
\item 
Let $F=A_{ij}A_{jk}A_{jl}A_{kq}A_{lq}A_{qi}$
We first compute the general sum of $F$ as 
\begin{align}\label{eq:5-e-F}
    \sum_{i,j,k,l,q}A_{ij}A_{jk}A_{jl}A_{kq}A_{lq}A_{qi} = \sum_{j,q}\sum_{i}A_{ij}A_{qi}\sum_{k}A_{jk}A_{kq}\sum_{l}A_{jl}A_{lq}=\onebm^\top (A^2)^{\circ3}\onebm.
\end{align}
Consider the graph $G = (V, E_G)$ with $E_G = \{i-j-k-q-i,j-l-q\}$ and its complement graph $H = (V, E_H)$ with $E_H = \{i-k-l-i,q-j\}$ (see Row 5 in \Cref{tab:result_5_nodes-1}). Using the same arguments in Part (a), the contraction $F_{2+1+1+1}^{(5)}$ can be obtained from the reduced graphs of $G$, shown in \Cref{fig:reduced_graph_5}, as follows:
\begin{align*}
    F_{2+1+1+1}^{(5)} & = 3A_{ij}^2A_{jk}A_{kl}A_{li}^2+A_{ij}^2A_{jk}^2A_{jl}^2,
\end{align*}
which has the sum 
\begin{align}
    \sum_{i,j,k,l}F_{2+1+1+1}^{(5)} = 3\tr((A^{\circ2})^2A^2)+\onebm^\top(A^{\circ2}\onebm)^{\circ3}.
\end{align}

\begin{figure}[!ht]
    \centering
    \subfloat[$i=k$]{%
        \begin{tikzpicture}[every node/.style={circle, draw, inner sep=0pt, minimum size=12pt,text centered},node distance=0.7cm]
            \node (i) at (0,2) {i};
            \node (j) [below=of i] {j};
            \node (l) [left=of j] {l};
            \node (q) [above=of l] {q};
            \draw (j)--(l)--(q);
            \draw (i) to[bend left=15] (j);
            \draw (i) to[bend right=15] (j);
            \draw (i) to[bend left=15] (q);
            \draw (i) to[bend right=15] (q);
        \end{tikzpicture}
    }
    \hspace{0.8cm} % space between subfigures
    \subfloat[$i=l$]{%
        \begin{tikzpicture}[every node/.style={circle, draw, inner sep=0pt, minimum size=12pt,text centered},node distance=0.7cm]
            \node (i) at (0,2) {i};
            \node (j) [below=of i] {j};
            \node (l) [left=of j] {k};
            \node (q) [above=of l] {q};
            \draw (j)--(l)--(q);
            \draw (i) to[bend left=15] (j);
            \draw (i) to[bend right=15] (j);
            \draw (i) to[bend left=15] (q);
            \draw (i) to[bend right=15] (q);
        \end{tikzpicture}
    }
    \hspace{0.8cm} % space between subfigures
    \subfloat[$i=k$]{%
        \begin{tikzpicture}[every node/.style={circle, draw, inner sep=0pt, minimum size=12pt,text centered},node distance=0.7cm]
            \node (i) at (0,2) {i};
            \node (j) [below=of i] {j};
            \node (k) [left=of j] {k};
            \node (q) [above=of k] {q};
            \draw (j)--(i)--(q);
            \draw (k) to[bend left=15] (j);
            \draw (k) to[bend right=15] (j);
            \draw (k) to[bend left=15] (q);
            \draw (k) to[bend right=15] (q);
        \end{tikzpicture}
    }
    \hspace{0.8cm} % space between subfigures
    \subfloat[$q=j$]{%
        \begin{tikzpicture}[every node/.style={circle, draw, inner sep=0pt, minimum size=12pt,text centered},node distance=0.7cm]
            \node (i) at (0,2) {i};
            \node (j) [below=of i] {j};
            \node (k) [left=of j] {k};
            \node (l) [above=of k] {l};
            \draw (i) to[bend left=15] (j);
            \draw (i) to[bend right=15] (j);
            \draw (l) to[bend left=15] (j);
            \draw (l) to[bend right=15] (j);
            \draw (k) to[bend left=15] (j);
            \draw (k) to[bend right=15] (j);
        \end{tikzpicture}
    }
    \caption{Reduced graphs of ($G$, Row 5 \Cref{tab:result_5_nodes-1}) by letting nodes of the edge in $E_H$ equal.}
    \label{fig:reduced_graph_5}
\end{figure}
For $F_{2+2+1}^{(5)}$, there are three disjoint subgraphs contained in $H$, given by $E_H^1=\{i-k,q-j\},E_H^2=\{i-l,q-j\},E_H^3=\{k-l,q-j\}$, which lead to 
\begin{align}
    F_{2+2+1}^{(5)} & = 3A_{ij}^4A_{jk}^2, \quad\sum_{i,j,k}F_{2+2+1}^{(5)} = 3\onebm^\top A^{\circ4}A^{\circ2}\onebm.
\end{align}
For $F_{3+1+1}^{(5)}$, there is one 3-clique in $H$, $E^1_H=\{l-i-k-l\}$ so that 
\begin{align}
    F_{3+1+1}^{(5)} & = A_{ij}^3A_{ik}^3, \quad \sum_{i,j,k}F_{3+1+1}^{(5)} = \onebm^\top (A^{\circ3})^2\onebm.
\end{align}
The last nonzero contraction $F_{3+2}^{(5)}$ is obtained by setting $q=k=i$ and $j=l$ that $F_{3+2}^{(5)} = A_{ij}^6$ and
\begin{align}\label{eq:5-e-F32}
    \sum_{i,j}F_{3+2}^{(5)} = \onebm^\top A^{\circ6}\onebm.
\end{align}
Using the formula in \Cref{sec:summary} and combining \eqref{eq:5-e-F}-\eqref{eq:5-e-F32}, we have that 
\begin{align*}
    &\sum_{i\neq j\neq k\neq l\neq q}A_{ij}A_{jk}A_{jl}A_{kq}A_{lq}A_{qi} \\
    =& \sum_{i,j,k,l,q}F-\sum_{i,j,k,l}F_{2+1+1+1}^{(5)}+\sum_{i,j,k}F_{2+2+1}^{(5)}+2\sum_{i,j,k}F_{3+1+1}^{(5)}-2\sum_{i,j}F_{3+2}^{(5)}\\
    = & \onebm^\top \Big[ (A^2)^{\circ3}+3A^{\circ4}A^{\circ2})+2(A^{\circ3})^2-2A^{\circ6}\Big]\onebm-\onebm^\top(A^{\circ2}\onebm)^{\circ3}-3\tr((A^{\circ2})^2A^2).
\end{align*}
\item Let $F=A_{ij}A_{jk}A_{jq}A_{kl}A_{qi}$.
We first compute the general sum of $F$ as 
\begin{align}\label{eq:5-f-F}
    \sum_{i,j,k,l,q}A_{ij}A_{jk}A_{jq}A_{kl}A_{qi} = \sum_{j,l}\sum_{i,q}A_{ij}A_{jq}A_{qi}\sum_{k}A_{jk}A_{kl}=\diag(A^3)^\top A^2\onebm.
\end{align}
Consider the graph $G = (V, E_G)$ with $E_G = \{j-i-q-j-k-l\}$ and its complement graph $H = (V, E_H)$ with $E_H = \{j-l-q-k-i-l\}$ (see Row 6 in \Cref{tab:result_5_nodes-1}). Using the same arguments in Part (a), the contraction $F_{2+1+1+1}^{(5)}$ can be obtained from the reduced graphs of $G$, shown in \Cref{fig:reduced_graph_6}, as follows:
\begin{align*}
    F_{2+1+1+1}^{(5)} & = 2A_{ij}A_{jk}^2A_{ki}A_{kl}+2A_{ij}A_{jk}A_{ki}A_{kl}A_{li}+A_{ij}A_{jl}A_{li}A_{jk}^2,
\end{align*}
which has the sum 
\begin{align}
    \sum_{i,j,k,l}F_{2+1+1+1}^{(5)} = 2\diag(A^{\circ2}A^2)^\top A\onebm+2\onebm^\top ((A^2)^{\circ2}\odot A)\onebm+\diag(A^3)^\top A^{\circ2}\onebm.
\end{align}
\begin{figure}[!ht]
    \centering
    \subfloat[$i=k$]{%
        \begin{tikzpicture}[every node/.style={circle, draw, inner sep=0pt, minimum size=12pt,text centered},node distance=0.7cm]
            \node (j) at (0,2) {j};
            \node (k) [below=of j] {k};
            \node (l) [left=of k] {l};
            \node (q) [above=of l] {q};
            \draw (l)--(k)--(q)--(j);
            \draw (k) to[bend left=15] (j);
            \draw (k) to[bend right=15] (j);
        \end{tikzpicture}
    }
    \hspace{0.8cm} % space between subfigures
    \subfloat[$i=l$]{%
        \begin{tikzpicture}[every node/.style={circle, draw, inner sep=0pt, minimum size=12pt,text centered},node distance=0.7cm]
            \node (j) at (0,2) {j};
            \node (k) [below=of j] {k};
            \node (l) [left=of k] {l};
            \node (q) [above=of l] {q};
            \draw (j)--(k)--(l)--(q)--(j)--(l);
        \end{tikzpicture}
    }
    \hspace{0.8cm} % space between subfigures
    \subfloat[$j=l$]{%
        \begin{tikzpicture}[every node/.style={circle, draw, inner sep=0pt, minimum size=12pt,text centered},node distance=0.7cm]
            \node (i) at (0,2) {i};
            \node (j) [below=of i] {j};
            \node (k) [left=of j] {k};
            \node (q) [above=of k] {q};
            \draw (j)--(i)--(q)--(j);
            \draw (k) to[bend left=15] (j);
            \draw (k) to[bend right=15] (j);
            
        \end{tikzpicture}
    }
    \hspace{0.8cm} % space between subfigures
    \subfloat[$k=q$]{%
        \begin{tikzpicture}[every node/.style={circle, draw, inner sep=0pt, minimum size=12pt,text centered},node distance=0.7cm]
            \node (i) at (0,2) {i};
            \node (j) [below=of i] {j};
            \node (k) [left=of j] {k};
            \node (l) [above=of k] {l};
            \draw (l)--(k)--(i)--(j);
            \draw (k) to[bend left=15] (j);
            \draw (k) to[bend right=15] (j);
        \end{tikzpicture}
    }
    \hspace{0.8cm} % space between subfigures
    \subfloat[$q=l$]{%
        \begin{tikzpicture}[every node/.style={circle, draw, inner sep=0pt, minimum size=12pt,text centered},node distance=0.7cm]
            \node (i) at (0,2) {i};
            \node (j) [below=of i] {j};
            \node (k) [left=of j] {k};
            \node (l) [above=of k] {l};
            \draw (l)--(k)--(j)--(l)--(i)--(j);
        \end{tikzpicture}
    }
    \caption{Reduced graphs of ($G$, Row 6 \Cref{tab:result_5_nodes-1}) by letting nodes of the edge in $E_H$ equal.}
    \label{fig:reduced_graph_6}
\end{figure}

For $F_{2+2+1}^{(5)}$, there are four disjoint subgraphs contained in $H$, given by $E_H^1=\{i-k,j-l\},E_H^2=\{i-k,l-q\},E_H^3=\{i-l,q-k\},E_H^4=\{j-l,q-k\}$, which lead to 
\begin{align}\label{eq:5-f-F221}
    F_{2+2+1}^{(5)} & = 2A_{ij}^3A_{jk}A_{ki}+2A_{ij}^2A_{jk}^2A_{ki}, \quad\sum_{i,j,k}F_{2+2+1}^{(5)} = 2\tr(A^{\circ3}A^2)+2\tr((A^{\circ2})^2A).
\end{align}
Using the formula in \Cref{sec:summary} and combining \eqref{eq:5-f-F}-\eqref{eq:5-f-F221}, we have that 
\begin{align*}
    &\sum_{i\neq j\neq k\neq l\neq q}A_{ij}A_{jk}A_{jq}A_{kl}A_{qi} = \sum_{i,j,k,l,q}F-\sum_{i,j,k,l}F_{2+1+1+1}^{(5)}+\sum_{i,j,k}F_{2+2+1}^{(5)}\\
    = & \diag(A^3)^\top A^2\onebm-2\diag(A^{\circ2}A^2)^\top A\onebm-\diag(A^3)^\top A^{\circ2}\onebm-2\onebm^\top ((A^2)^{\circ2}\odot A)\onebm \\
    & +2\tr(A^{\circ3}A^2)+2\tr((A^{\circ2})^2A).
\end{align*}

\item Let $F=A_{ij}A_{jk}A_{jq}A_{lq}A_{qi}$. We first compute the general sum of $F$ as 
\begin{align}\label{eq:5-g-F}
    \sum_{i,j,k,l,q}A_{ij}A_{jk}A_{jq}A_{lq}A_{qi} = \sum_{j,q}A_{jq}\sum_{i}A_{ji}A_{iq}\sum_kA_{jk}\sum_{l}A_{ql}=\onebm^\top A(A^2\odot A)A\onebm.
\end{align}
Consider the graph $G = (V, E_G)$ with $E_G = \{k-j-i-q-j,q-l\}$ and its complement graph $H = (V, E_H)$ with $E_H = \{i-k-l-i,q-k,l-j\}$ (see Row 7 in \Cref{tab:result_5_nodes-1}). Using the same arguments in Part (a), the contraction $F_{2+1+1+1}^{(5)}$ can be obtained from the reduced graphs of $G$, shown in \Cref{fig:reduced_graph_7}, as follows:
\begin{align*}
    F_{2+1+1+1}^{(5)} & = 2A_{ij}^2A_{jk}A_{ki}A_{kl}+2A_{ij}A_{jk}A_{jl}^2A_{li}+A_{ij}A_{jk}A_{kl}A_{li}A_{jl},
\end{align*}
which has the sum 
\begin{align}
    \sum_{i,j,k,l}F_{2+1+1+1}^{(5)} = 2\diag(A^{\circ2}A^2)^\top A\onebm+2\diag(AA^{\circ2}A)^\top A\onebm+\onebm^\top ((A^2)^{\circ2}\odot A)\onebm.
\end{align}
\begin{figure}[!ht]
    \centering
    \subfloat[$i=k$]{%
        \begin{tikzpicture}[every node/.style={circle, draw, inner sep=0pt, minimum size=12pt,text centered},node distance=0.7cm]
            \node (j) at (0,2) {j};
            \node (k) [below=of j] {k};
            \node (l) [left=of k] {l};
            \node (q) [above=of l] {q};
            \draw (k)--(q)--(j);
            \draw (l)--(q);
            \draw (k) to[bend left=15] (j);
            \draw (k) to[bend right=15] (j);
        \end{tikzpicture}
    }
    \hspace{0.8cm} % space between subfigures
    \subfloat[$i=l$]{%
        \begin{tikzpicture}[every node/.style={circle, draw, inner sep=0pt, minimum size=12pt,text centered},node distance=0.7cm]
            \node (j) at (0,2) {j};
            \node (k) [below=of j] {k};
            \node (l) [left=of k] {l};
            \node (q) [above=of l] {q};
            \draw (j)--(k);
            \draw (j)--(l);
            \draw (j)--(q);
            \draw (q) to[bend left=15] (l);
            \draw (q) to[bend right=15] (l);
        \end{tikzpicture}
    }
    \hspace{0.8cm} % space between subfigures
    \subfloat[$j=l$]{%
        \begin{tikzpicture}[every node/.style={circle, draw, inner sep=0pt, minimum size=12pt,text centered},node distance=0.7cm]
            \node (i) at (0,2) {i};
            \node (j) [below=of i] {j};
            \node (k) [left=of j] {k};
            \node (q) [above=of k] {q};
            \draw (q)--(i)--(j)--(k);
            \draw (q) to[bend left=15] (j);
            \draw (q) to[bend right=15] (j);
            
        \end{tikzpicture}
    }
    \hspace{0.8cm} % space between subfigures
    \subfloat[$k=q$]{%
        \begin{tikzpicture}[every node/.style={circle, draw, inner sep=0pt, minimum size=12pt,text centered},node distance=0.7cm]
            \node (i) at (0,2) {i};
            \node (j) [below=of i] {j};
            \node (l) [left=of j] {l};
            \node (q) [above=of l] {q};
            \draw (j)--(i)--(q)--(l);
            \draw (q) to[bend left=15] (j);
            \draw (q) to[bend right=15] (j);
        \end{tikzpicture}
    }
    \hspace{0.8cm} % space between subfigures
    \subfloat[$k=l$]{%
        \begin{tikzpicture}[every node/.style={circle, draw, inner sep=0pt, minimum size=12pt,text centered},node distance=0.7cm]
            \node (i) at (0,2) {i};
            \node (j) [below=of i] {j};
            \node (l) [left=of j] {l};
            \node (q) [above=of k] {q};
            \draw (i)--(j)--(l)--(q)--(i);
            \draw (q)--(j);
        \end{tikzpicture}
    }
    \caption{Reduced graphs of ($G$, Row 7 \Cref{tab:result_5_nodes-1}) by letting nodes of the edge in $E_H$ equal.}
    \label{fig:reduced_graph_7}
\end{figure}
For $F_{2+2+1}^{(5)}$, there are three disjoint subgraphs contained in $H$, given by $E_H^1=\{i-k,j-l\},E_H^2=\{i-l,k-q\},E_H^3=\{j-l,k-q\}$, which lead to 
\begin{align}
    F_{2+2+1}^{(5)} & = 2A_{ij}^2A_{jk}^2A_{ki}+A_{ij}A_{jk}^3A_{ki}, \quad\sum_{i,j,k}F_{2+2+1}^{(5)} = 2\tr((A^{\circ2})^2A)+\tr(A^{\circ3}A^2).
\end{align}

For $F_{3+1+1}^{(5)}$, there is one 3-clique in $H$, $E^1_H=\{l-i-k-l\}$ so that 
\begin{align}\label{eq:5-g-F311}
    F_{3+1+1}^{(5)} & = A_{ij}^2A_{jk}A_{ki}^2, \quad \sum_{i,j,k}F_{3+1+1}^{(5)} = \tr((A^{\circ2})^2A).
\end{align}

Using the formula in \Cref{sec:summary} and combining \eqref{eq:5-g-F}-\eqref{eq:5-g-F311}, we have that 
\begin{align*}
    &\sum_{i\neq j\neq k\neq l\neq q}A_{ij}A_{jk}A_{jq}A_{lq}A_{qi} = \sum_{i,j,k,l,q}F-\sum_{i,j,k,l}F_{2+1+1+1}^{(5)}+\sum_{i,j,k}F_{2+2+1}^{(5)}+2\sum_{i,j,k}F_{3+1+1}^{(5)}\\
    = & \onebm^\top \Big[ A(A^2\odot A)A-(A^2)^{\circ2}\odot A\Big]\onebm{-2\diag(A^{\circ2}A^2)^\top A\onebm-2\diag(AA^{\circ2}A)^\top A\onebm}+4\tr((A^{\circ2})^2A)+\tr(A^{\circ3}A^2).
\end{align*}
\item Let $F=A_{ij}A_{jk}A_{jq}A_{kl}A_{lq}A_{qi}$. We first compute the general sum of $F$ as 
\begin{align}\label{eq:5-h-F}
    \sum_{i,j,k,l,q}A_{ij}A_{jk}A_{jq}A_{kl}A_{lq}A_{qi} = \sum_{j,q}A_{jq}\sum_{i}A_{ji}A_{iq}\sum_{k,l}A_{jk}A_{kl}A_{lq}=\onebm^\top (A\odot A^2\odot A^3)\onebm.
\end{align}
Consider the graph $G = (V, E_G)$ with $E_G = \{i-j-q-i,j-k-l-q\}$ and its complement graph $H = (V, E_H)$ with $E_H = \{j-l-i-k-q\}$ (see Row 1 in \Cref{tab:result_5_nodes-2}). Using the same arguments in Part (a), the contraction $F_{2+1+1+1}^{(5)}$ can be obtained from the reduced graphs of $G$, shown in \Cref{fig:reduced_graph_8}, as follows:
\begin{align*}
    F_{2+1+1+1}^{(5)} & = 2A_{ij}A_{jk}A_{ki}A_{kl}^2A_{li}+2A_{ij}A_{jk}^2A_{ki}A_{jl}^2,
\end{align*}
which has the sum 
\begin{align}
    \sum_{i,j,k,l}F_{2+1+1+1}^{(5)} = 2\onebm^\top(A\odot A^2\odot (A^{\circ2}A))\onebm+2\diag(A^{\circ2}A^2)^\top A^{\circ2}\onebm.
\end{align}

\begin{figure}[!ht]
    \centering
    \subfloat[$i=k$]{%
        \begin{tikzpicture}[every node/.style={circle, draw, inner sep=0pt, minimum size=12pt,text centered},node distance=0.7cm]
            \node (j) at (0,2) {j};
            \node (k) [below=of j] {k};
            \node (l) [left=of k] {l};
            \node (q) [above=of l] {q};
            \draw (k)--(l)--(q)--(k);
            \draw (q)--(j);         
            \draw (k) to[bend left=15] (j);
            \draw (k) to[bend right=15] (j);
        \end{tikzpicture}
    }
    \hspace{0.8cm} % space between subfigures
    \subfloat[$i=l$]{%
        \begin{tikzpicture}[every node/.style={circle, draw, inner sep=0pt, minimum size=12pt,text centered},node distance=0.7cm]
            \node (j) at (0,2) {j};
            \node (k) [below=of j] {k};
            \node (l) [left=of k] {l};
            \node (q) [above=of l] {q};
            \draw (q)--(j)--(k);
            \draw (j)--(l)--(k);
            \draw (q) to[bend left=15] (l);
            \draw (q) to[bend right=15] (l);    
        \end{tikzpicture}
    }
    \hspace{0.8cm} % space between subfigures
    \subfloat[$j=l$]{%
        \begin{tikzpicture}[every node/.style={circle, draw, inner sep=0pt, minimum size=12pt,text centered},node distance=0.7cm]
            \node (i) at (0,2) {i};
            \node (j) [below=of i] {j};
            \node (k) [left=of j] {k};
            \node (q) [above=of k] {q};
            \draw (q)--(i)--(j);
            \draw (q) to[bend left=15] (j);
            \draw (q) to[bend right=15] (j);
            \draw (k) to[bend left=15] (j);
            \draw (k) to[bend right=15] (j);
        \end{tikzpicture}
    }
    \hspace{0.8cm} % space between subfigures
    \subfloat[$k=q$]{%
        \begin{tikzpicture}[every node/.style={circle, draw, inner sep=0pt, minimum size=12pt,text centered},node distance=0.7cm]
            \node (i) at (0,2) {i};
            \node (j) [below=of i] {j};
            \node (l) [left=of j] {l};
            \node (q) [above=of l] {q};
            \draw (q)--(i)--(j);
            \draw (q) to[bend left=15] (j);
            \draw (q) to[bend right=15] (j);
            \draw (q) to[bend left=15] (l);
            \draw (q) to[bend right=15] (l);
        \end{tikzpicture}
    }
    \caption{Reduced graphs of ($G$, Row 1 \Cref{tab:result_5_nodes-2}) by letting nodes of the edge in $E_H$ equal.}
    \label{fig:reduced_graph_8}
\end{figure}
For $F_{2+2+1}^{(5)}$, there are three disjoint subgraphs contained in $H$, given by $E_H^1=\{i-k,j-l\},E_H^2=\{i-l,k-q\},E_H^3=\{j-l,k-q\}$, which lead to 
\begin{align}\label{eq:5-h-F221}
    F_{2+2+1}^{(5)} & = 2A_{ij}^3A_{jk}^2A_{ki}+A_{ij}A_{jk}^4A_{ki}, \quad\sum_{i,j,k}F_{2+2+1}^{(5)} = 2\tr(AA^{\circ2}A^{\circ3})+\tr(A^{\circ4}A^2).
\end{align}
Using the formula in \Cref{sec:summary} and combining \eqref{eq:5-h-F}-\eqref{eq:5-h-F221}, we have that 
\begin{align*}
    &\sum_{i\neq j\neq k\neq l\neq q}A_{ij}A_{jk}A_{jq}A_{kl}A_{lq}A_{qi} = \sum_{i,j,k,l,q}F-\sum_{i,j,k,l}F_{2+1+1+1}^{(5)}+\sum_{i,j,k}F_{2+2+1}^{(5)}\\
    = & \onebm^\top \Big[ A\odot A^2\odot A^3-2A\odot A^2\odot (A^{\circ2}A)\Big]\onebm+2\tr(AA^{\circ2}A^{\circ3})+\tr(A^{\circ4}A^2) \\
    & -2\diag(A^{\circ2}A^2)^\top A^{\circ2}\onebm.
\end{align*}

\item Let $F=A_{jk}A_{jl}A_{jq}A_{kl}A_{lq}A_{qi}$. We first compute the general sum of $F$ as 
\begin{align}\label{eq:5-i-F}
    \sum_{i,j,k,l,q}A_{jk}A_{jl}A_{jq}A_{kl}A_{lq}A_{qi}& = \sum_{i,q}A_{iq}\sum_{j,l}A_{qj}A_{jl}A_{lq}\sum_{k}A_{jk}A_{kl}=\sum_{i,q}A_{iq}\sum_{j,l}A_{qj}[A\odot A^2]_{jl}A_{lq} \nonumber\\
    & = \diag(A(A\odot A^2)A)^\top A\onebm
\end{align}
Consider the graph $G = (V, E_G)$ with $E_G = \{i-q-l-k-j-q,l-j\}$ and its complement graph $H = (V, E_H)$ with $E_H = \{j-i-k,i-l,q-k\}$ (see Row 2 in \Cref{tab:result_5_nodes-2}). Using the same arguments in Part (a), the contraction $F_{2+1+1+1}^{(5)}$ can be obtained from the reduced graphs of $G$, shown in \Cref{fig:reduced_graph_9}, as follows:
\begin{align*}
    F_{2+1+1+1}^{(5)} & = 2A_{ij}A_{jk}A_{ki}A_{kl}A_{li}^2+A_{ij}A_{jk}A_{kl}A_{li}A_{ik}A_{jl}+A_{il}A_{jk}A_{kl}^2A_{lj}^2,
\end{align*}
which has the sum 
\begin{align}
    \sum_{i,j,k,l}F_{2+1+1+1}^{(5)} & = 2\tr(AA^{\circ2}(A\odot A^2))+\diag(A^{\circ2}AA^{\circ2})^\top A\onebm \nonumber\\
    & +\langle (\Acal_{(2)}\odot \Acal_{(3)})\times_{1,2}(\Acal_{(2)}\odot \Acal_{(3)})\times_{1,2}(\Acal_{(2)}\odot \Acal_{(3)}),\Ical_{n\times n\times n}\rangle,
\end{align}
where the last equality follows from Part (f),~\Cref{prop:4-node}.

\begin{figure}[!ht]
    \centering
    \subfloat[$i=j$]{%
        \begin{tikzpicture}[every node/.style={circle, draw, inner sep=0pt, minimum size=12pt,text centered},node distance=0.7cm]
            \node (j) at (0,2) {j};
            \node (k) [below=of j] {k};
            \node (l) [left=of k] {l};
            \node (q) [above=of l] {q};
            \draw (q)--(l)--(k)--(j)--(l);
            \draw (q) to[bend left=15] (j);
            \draw (q) to[bend right=15] (j);
        \end{tikzpicture}
    }
    \hspace{0.8cm} % space between subfigures
    \subfloat[$i=k$]{%
        \begin{tikzpicture}[every node/.style={circle, draw, inner sep=0pt, minimum size=12pt,text centered},node distance=0.7cm]
            \node (j) at (0,2) {j};
            \node (k) [below=of j] {k};
            \node (l) [left=of k] {l};
            \node (q) [above=of l] {q};
            \draw (q)--(j)--(k)--(l)--(q)--(k);
            \draw (j)--(l);  
        \end{tikzpicture}
    }
    \hspace{0.8cm} % space between subfigures
    \subfloat[$i=l$]{%
        \begin{tikzpicture}[every node/.style={circle, draw, inner sep=0pt, minimum size=12pt,text centered},node distance=0.7cm]
            \node (j) at (0,2) {j};
            \node (k) [below=of j] {k};
            \node (l) [left=of k] {l};
            \node (q) [above=of l] {q};
            \draw (q)--(j)--(l)--(k)--(j);
            \draw (q) to[bend left=15] (l);
            \draw (q) to[bend right=15] (l);
        \end{tikzpicture}
    }
    \hspace{0.8cm} % space between subfigures
    \subfloat[$k=q$]{%
        \begin{tikzpicture}[every node/.style={circle, draw, inner sep=0pt, minimum size=12pt,text centered},node distance=0.7cm]
            \node (i) at (0,2) {i};
            \node (j) [below=of i] {j};
            \node (l) [left=of j] {l};
            \node (q) [above=of l] {q};
            \draw (q)--(i);
            \draw (l)--(j);
            \draw (q) to[bend left=15] (j);
            \draw (q) to[bend right=15] (j);
            \draw (q) to[bend left=15] (l);
            \draw (q) to[bend right=15] (l);
        \end{tikzpicture}
    }
    \caption{Reduced graphs of ($G$, Row 2 \Cref{tab:result_5_nodes-2}) by letting nodes of the edge in $E_H$ equal.}
    \label{fig:reduced_graph_9}
\end{figure}
For $F_{2+2+1}^{(5)}$, there are two disjoint subgraphs contained in $H$, given by $E_H^1=\{i-j,k-q\},E_H^2=\{i-l,k-q\}$, which lead to 
\begin{align}\label{eq:5-i-F221}
    F_{2+2+1}^{(5)} & = 2A_{ij}A_{jk}^2A_{ki}^3, \quad\sum_{i,j,k}F_{2+2+1}^{(5)} = 2\tr(AA^{\circ2}A^{\circ3}).
\end{align}
Using the formula in \Cref{sec:summary} and combining \eqref{eq:5-i-F}-\eqref{eq:5-i-F221}, we have that 
\begin{align*}
    &\sum_{i\neq j\neq k\neq l\neq q}A_{jk}A_{jl}A_{jq}A_{kl}A_{lq}A_{qi} = \sum_{i,j,k,l,q}F-\sum_{i,j,k,l}F_{2+1+1+1}^{(5)}+\sum_{i,j,k}F_{2+2+1}^{(5)}\\
    = & \diag(A(A\odot A^2)A)^\top A\onebm-\diag(A^{\circ2}AA^{\circ2})^\top A\onebm -2\tr(AA^{\circ2}(A\odot A^2))+2\tr(AA^{\circ2}A^{\circ3})\\
    & -\langle (\Acal_{(2)}\odot \Acal_{(3)})\times_{1,2}(\Acal_{(2)}\odot \Acal_{(3)})\times_{1,2}(\Acal_{(2)}\odot \Acal_{(3)}),\Ical_{n\times n\times n}\rangle.
\end{align*}
\item Let $F=A_{ij}A_{jk}A_{jl}A_{kl}A_{kq}A_{lq}A_{qi}$. 
We first compute the general sum of $F$ as 
\begin{align*}
    \sum_{i,j,k,l,q}A_{ij}A_{jk}A_{jl}A_{kl}A_{kq}A_{lq}A_{qi}& = \sum_{j,k,l,q}A_{jk}A_{jl}A_{kl}A_{kq}A_{lq}\sum_iA_{qi}A_{ij} \\
    & = \sum_{i,j,k,l}[A^2]_{ij}A_{ik}A_{il}A_{jk}A_{jl}A_{kl},
\end{align*}
which can be evaluated using the same arguments in Part (f),~\Cref{prop:4-node}. For any fixed $l$, define $[T^l_1]_{ij}=A_{ij}A_{li}$ and $[T^l_2]_{ij}=[A^2]_{ij}A_{li}$ and we have that 
\begin{align*}
    \sum_l\sum_{i,j,k}[A^2]_{ij}A_{ik}A_{il}A_{jk}A_{jl}A_{kl}=\sum_l\sum_{i,j,k}([A^2]_{ij}A_{li})(A_{jk}A_{lj})(A_{ki}A_{lk})=\sum_l\tr((T_1^l)^2T_2^l).
\end{align*}
Let $\Acal_{(2)}, \Acal_{(3)},\Bcal_{(2)}\in\Rbb^{n\times n\times n}$ be the tensors with $[\Acal_{(2)}]_{:k',:}=[\Acal_{(3)}]_{::k'}=A$ and $[\Bcal_{(2)}]_{:k',:}=A^{2}$ for all $1\leq k'\leq n$. Then $\Tcal_1=\Acal_{(2)}\odot \Acal_{(3)}$ and $\Tcal_2=\Bcal_{(2)}\odot \Acal_{(3)}$ are the tensors that $[\Tcal_1]_{::l}=T_1^l$ and $[\Tcal_2]_{::l}=T_2^l$. Hence, it follows that 
\begin{align}\label{eq:5-j-F}
    \sum_{i,j,k,l,q}A_{ij}A_{jk}A_{jl}A_{kl}A_{kq}A_{lq}A_{qi}& = \langle (\Acal_{(2)}\odot\Acal_{(3)})\times_{1,2}(\Acal_{(2)}\odot\Acal_{(3)})\times_{1,2}(\Bcal_{(2)}\odot\Acal_{(3)}),\Ical_{n\times n\times n}\rangle.
\end{align}

Consider the graph $G = (V, E_G)$ with $E_G = \{i-j-k-l-q-i,q-k,l-j\}$ and its complement graph $H = (V, E_H)$ with $E_H = \{l-i-k,q-j\}$ (see Row 3 in \Cref{tab:result_5_nodes-2}). Using the same arguments in Part (a), the contraction $F_{2+1+1+1}^{(5)}$ can be obtained from the reduced graphs of $G$, shown in \Cref{fig:reduced_graph_10}, as follows:
\begin{align*}
    F_{2+1+1+1}^{(5)} & = 2A_{ij}^2A_{jk}A_{ki}A_{kl}A_{li}^2+A_{ij}^2A_{jk}^2A_{kl}A_{lj}^2,
\end{align*}
which has the sum 
\begin{align}
    \sum_{i,j,k,l}F_{2+1+1+1}^{(5)} & = 2\tr([(A^{\circ2}A)\odot A]AA^{\circ2})+\diag(A^{\circ2}AA^{\circ2})^\top A^{\circ2}\onebm.
\end{align}

\begin{figure}[!ht]
    \centering
    \subfloat[$i=k$]{%
        \begin{tikzpicture}[every node/.style={circle, draw, inner sep=0pt, minimum size=12pt,text centered},node distance=0.7cm]
            \node (i) at (0,2) {i};
            \node (j) [below=of i] {j};
            \node (k) [left=of j] {l};
            \node (l) [above=of k] {q};
            \draw (l)--(k)--(j);
            \draw (k)--(i);
            \draw (i) to[bend left=15] (j);
            \draw (i) to[bend right=15] (j);
            \draw (i) to[bend left=15] (l);
            \draw (i) to[bend right=15] (l);
        \end{tikzpicture}
    }
    \hspace{0.8cm} % space between subfigures
    \subfloat[$i=l$]{%
        \begin{tikzpicture}[every node/.style={circle, draw, inner sep=0pt, minimum size=12pt,text centered},node distance=0.7cm]
            \node (i) at (0,2) {i};
            \node (j) [below=of i] {j};
            \node (k) [left=of j] {k};
            \node (l) [above=of k] {q};
            \draw (l)--(k)--(j);
            \draw (k)--(i);
            \draw (i) to[bend left=15] (j);
            \draw (i) to[bend right=15] (j);
            \draw (i) to[bend left=15] (l);
            \draw (i) to[bend right=15] (l);
        \end{tikzpicture}
    }
    \hspace{0.8cm} % space between subfigures
    \subfloat[$q=j$]{%
        \begin{tikzpicture}[every node/.style={circle, draw, inner sep=0pt, minimum size=12pt,text centered},node distance=0.7cm]
            \node (i) at (0,2) {i};
            \node (j) [below=of i] {j};
            \node (k) [left=of j] {k};
            \node (l) [above=of k] {l};
            \draw (l)--(k);
            \draw (i) to[bend left=15] (j);
            \draw (i) to[bend right=15] (j);
            \draw (k) to[bend left=15] (j);
            \draw (k) to[bend right=15] (j);
            \draw (l) to[bend left=15] (j);
            \draw (l) to[bend right=15] (j);
        \end{tikzpicture}
    }
    \caption{Reduced graphs of ($G$, Row 3 \Cref{tab:result_5_nodes-2}) by letting nodes of the edge in $E_H$ equal.}
    \label{fig:reduced_graph_10}
\end{figure}
For $F_{2+2+1}^{(5)}$, there are two disjoint subgraphs contained in $H$, given by $E_H^1=\{i-k,j-q\},E_H^2=\{i-l,j-q\}$, which lead to 
\begin{align}\label{eq:5-j-F221}
    F_{2+2+1}^{(5)} & = 2A_{ij}^4A_{jk}^2A_{ki}, \quad\sum_{i,j,k}F_{2+2+1}^{(5)} = 2\tr(AA^{\circ2}A^{\circ4}).
\end{align}
Using the formula in \Cref{sec:summary} and combining \eqref{eq:5-j-F}-\eqref{eq:5-j-F221}, we have that 
\begin{align*}
    &\sum_{i\neq j\neq k\neq l\neq q}A_{ij}A_{jk}A_{jl}A_{kl}A_{kq}A_{lq}A_{qi} = \sum_{i,j,k,l,q}F-\sum_{i,j,k,l}F_{2+1+1+1}^{(5)}+\sum_{i,j,k}F_{2+2+1}^{(5)}\\
    = & \langle (\Acal_{(2)}\odot\Acal_{(3)})\times_{1,2}(\Acal_{(2)}\odot\Acal_{(3)})\times_{1,2}(\Bcal_{(2)}\odot\Acal_{(3)}),\Ical_{n\times n\times n}\rangle-2\tr([(A^{\circ2}A)\odot A]AA^{\circ2}) \\
    &-\diag(A^{\circ2}AA^{\circ2})^\top A^{\circ2}\onebm+2\tr(AA^{\circ2}A^{\circ4}).
\end{align*}
\item Let $F = A_{ij}A_{ik}A_{il}A_{qi}$ and $H=\{V,E\}$ with $E=\{(j,k),(k,l),(l,q),(q,j),(l,j),(q,k)\}$. Except the node $i$, the remaining nodes $j,k,l,q$ forms a 4-clique in $H$. Using the similar arguments in Part (a), we can verify that the components in $F_{2+1+1+1}^{(5)},F_{2+2+1}^{(5)},F_{3+1+1}^{(5)}$ are equivalent. Hence we have the following nonzero contractions:
\begin{align*}
    F_{2+1+1+1}^{(5)} & \sim 6F(i,j,k,l,j)=6A_{ij}^2A_{ik}A_{il}, \\
    F_{2+2+1}^{(5)}& \sim 3F(i,j,k,k,j) = 3A_{ij}^2A_{ik}^2, \\
    F_{3+1+1}^{(5)} & \sim 4F(i,j,k,j,j) = 4A_{ij}^3A_{ik}, \\
    F_{4+1}^{(5)} & \sim F(i,j,j,j,j) = A_{ij}^4,
\end{align*}
Using the formula in \Cref{sec:summary}, we have that 
\begin{align*}
    & \sum_{i\neq j\neq k\neq l\neq q}A_{ij}A_{ik}A_{il}A_{iq} \\
    & = \sum_{i,j,k,l,q}F-\sum_{i,j,k,l}F_{2+1+1+1}^{(5)}+\sum_{i,j,k}F_{2+2+1}^{(5)}+2\sum_{i,j,k,l}F_{3+1+1}^{(5)}-6\sum_{i,j}F_{4+1}^{(5)}\\
    & =  \sum_{i,j,k,l,q}A_{ij}A_{ik}A_{il}A_{iq}-6\sum_{i,j,k,l}A_{ij}^2A_{ik}A_{il}+3\sum_{i,j,k}A_{ij}^2A_{ik}^2+8\sum_{i,j,k}A_{ij}^3A_{ik}-6\sum_{i,j}A_{ij}^4 \\
    & = \onebm^\top(A\onebm)^{\circ4}-6\onebm^\top A^{\circ2}(A\onebm)^{\circ2}+3\onebm^\top(A^{\circ2})^{2}\onebm+8\onebm^\top A^{\circ3}A\onebm-6\onebm^\top A^{\circ4}\onebm.
\end{align*}

\item Let $F=A_{ij}A_{ik}A_{jk}A_{il}A_{qi}$. We first compute the general sum of $F$ as 
\begin{align}\label{eq:5-l-F}
    \sum_{i,j,k,l,q}A_{ij}A_{ik}A_{jk}A_{il}A_{qi}& = \sum_{i}\sum_lA_{il}\sum_qA_{iq}\sum_k A_{ik}\sum_jA_{ij}A_{jk}=\onebm^\top (A\odot A^2) (A\onebm)^{\circ2}
\end{align}
Consider the graph $G = (V, E_G)$ with $E_G = \{i-j-k-i,q-i-l\}$ and its complement graph $H = (V, E_H)$ with $E_H = \{q-j-l-q-k-l\}$ (see Row 5 in \Cref{tab:result_5_nodes-2}). Using the same arguments in Part (a), the contraction $F_{2+1+1+1}^{(5)}$ can be obtained from the reduced graphs of $G$, shown in \Cref{fig:reduced_graph_12}, as follows:
\begin{align*}
    F_{2+1+1+1}^{(5)} & = 2A_{ij}^2A_{jk}A_{ki}A_{il}+2A_{ij}A_{jk}A_{ki}^2A_{il}+A_{ij}A_{jk}A_{ki}A_{il}^2,
\end{align*}
which has the sum 
\begin{align}
    \sum_{i,j,k,l}F_{2+1+1+1}^{(5)} & = 4\diag(A^2A^{\circ2})^\top A\onebm+\diag(A^3)^\top A^{\circ2}\onebm.
\end{align}

\begin{figure}[!ht]
    \centering
    \subfloat[$j=l$]{%
        \begin{tikzpicture}[every node/.style={circle, draw, inner sep=0pt, minimum size=12pt,text centered},node distance=0.7cm]
            \node (i) at (0,2) {i};
            \node (j) [below=of i] {j};
            \node (k) [left=of j] {k};
            \node (q) [above=of k] {q};
            \draw (q)--(i)--(k)--(j);
            \draw (i) to[bend left=15] (j);
            \draw (i) to[bend right=15] (j);
        \end{tikzpicture}
    }
    \hspace{0.8cm} % space between subfigures
    \subfloat[$j=q$]{%
        \begin{tikzpicture}[every node/.style={circle, draw, inner sep=0pt, minimum size=12pt,text centered},node distance=0.7cm]
            \node (i) at (0,2) {i};
            \node (j) [below=of i] {j};
            \node (k) [left=of j] {k};
            \node (l) [above=of k] {l};
            \draw (l)--(i)--(k)--(j);
            \draw (i) to[bend left=15] (j);
            \draw (i) to[bend right=15] (j);
        \end{tikzpicture}
    }
    \hspace{0.8cm} % space between subfigures
    \subfloat[$k=l$]{%
        \begin{tikzpicture}[every node/.style={circle, draw, inner sep=0pt, minimum size=12pt,text centered},node distance=0.7cm]
            \node (i) at (0,2) {i};
            \node (j) [below=of i] {j};
            \node (k) [left=of j] {k};
            \node (q) [above=of k] {q};
            \draw (q)--(i)--(j)--(k);
            \draw (i) to[bend left=15] (k);
            \draw (i) to[bend right=15] (k);
        \end{tikzpicture}
    }
    \hspace{0.8cm} % space between subfigures
    \subfloat[$k=q$]{%
        \begin{tikzpicture}[every node/.style={circle, draw, inner sep=0pt, minimum size=12pt,text centered},node distance=0.7cm]
            \node (i) at (0,2) {i};
            \node (j) [below=of i] {j};
            \node (k) [left=of j] {k};
            \node (q) [above=of k] {l};
            \draw (q)--(i)--(j)--(k);
            \draw (i) to[bend left=15] (k);
            \draw (i) to[bend right=15] (k);
        \end{tikzpicture}
    }
    \hspace{0.8cm} % space between subfigures
    \subfloat[$q=l$]{%
        \begin{tikzpicture}[every node/.style={circle, draw, inner sep=0pt, minimum size=12pt,text centered},node distance=0.7cm]
            \node (i) at (0,2) {i};
            \node (j) [below=of i] {j};
            \node (k) [left=of j] {k};
            \node (l) [above=of k] {l};
            \draw (k)--(i)--(j)--(k);
            \draw (i) to[bend left=15] (l);
            \draw (i) to[bend right=15] (l);
        \end{tikzpicture}
    }
    \caption{Reduced graphs of ($G$, Row 5 \Cref{tab:result_5_nodes-2}) by letting nodes of the edge in $E_H$ equal.}
    \label{fig:reduced_graph_12}
\end{figure}
For $F_{2+2+1}^{(5)}$, there are two disjoint subgraphs contained in $H$, given by $E_H^1=\{j-l,k-q\},E_H^2=\{j-q,k-l\}$, which lead to 
\begin{align}
    F_{2+2+1}^{(5)} & = 2A_{ij}^2A_{jk}A_{ki}^2, \quad\sum_{i,j,k}F_{2+2+1}^{(5)} = 2\tr((A^{\circ2})^2A).
\end{align}

For $F_{3+1+1}^{(5)}$, there are two disjoint subgraphs contained in $H$, given by $E_H^1=\{j-l-q-j\},E_H^2=\{q-k-l-q\}$, which lead to 
\begin{align}\label{eq:5-l-F311}
    F_{3+1+1}^{(5)} & = 2A_{ij}^3A_{jk}A_{ki}, \quad\sum_{i,j,k}F_{3+1+1}^{(5)} = 2\tr(A^{\circ3}A^2).
\end{align}

Using the formula in \Cref{sec:summary} and combining \eqref{eq:5-l-F}-\eqref{eq:5-l-F311}, we have that 
\begin{align*}
    &\sum_{i\neq j\neq k\neq l\neq q}A_{ij}A_{ij}A_{ik}A_{jk}A_{il}A_{qi} = \sum_{i,j,k,l,q}F-\sum_{i,j,k,l}F_{2+1+1+1}^{(5)}+\sum_{i,j,k}F_{2+2+1}^{(5)}+2\sum_{i,j,k}F_{3+1+1}^{(5)}\\
    = & \onebm^\top (A\odot A^2) (A\onebm)^{\circ2}-4\diag(A^2A^{\circ2})^\top A\onebm-\diag(A^3)^\top A^{\circ2}\onebm+2\tr((A^{\circ2})^2A)+4\tr(A^{\circ3}A^2).
\end{align*}
\item Let $F=A_{ij}A_{ik}A_{jk}A_{il}A_{lq}A_{qi}$. We first compute the general sum of $F$ as 
\begin{align}\label{eq:5-m-F}
    \sum_{i,j,k,l,q}A_{ij}A_{ik}A_{jk}A_{il}A_{lq}A_{qi}& = \sum_{i}\sum_{j,k}A_{ij}A_{jk}A_{ki}\sum_{l,q}A_{il}A_{lq}A_{qi}=\tr((A^3)^{\circ2}).
\end{align}
Consider the graph $G = (V, E_G)$ with $E_G = \{i-j-k-i-l-q-i\}$ and its complement graph $H = (V, E_H)$ with $E_H = \{j-l-k-q-j\}$ (see Row 6 in \Cref{tab:result_5_nodes-2}). Using the same arguments in Part (a), the contraction $F_{2+1+1+1}^{(5)}$ can be obtained from the reduced graphs of $G$, shown in \Cref{fig:reduced_graph_13}, as follows:
\begin{align*}
    F_{2+1+1+1}^{(5)} & = 2A_{ij}^2A_{jk}A_{ki}A_{jl}A_{li}+2A_{ij}A_{jk}A_{kl}A_{li}A_{ik}^2,
\end{align*}
which has the sum 
\begin{align}
    \sum_{i,j,k,l}F_{2+1+1+1}^{(5)} & = 4\onebm^\top ((A^2)^{\circ2}\odot A^{\circ2})\onebm.
\end{align}

\begin{figure}[!ht]
    \centering
    \subfloat[$j=l$]{%
        \begin{tikzpicture}[every node/.style={circle, draw, inner sep=0pt, minimum size=12pt,text centered},node distance=0.7cm]
            \node (i) at (0,2) {i};
            \node (j) [below=of i] {j};
            \node (k) [left=of j] {k};
            \node (q) [above=of k] {q};
            \draw (q)--(i)--(k)--(j)--(q);
            \draw (i) to[bend left=15] (j);
            \draw (i) to[bend right=15] (j);
        \end{tikzpicture}
    }
    \hspace{0.8cm} % space between subfigures
    \subfloat[$j=q$]{%
        \begin{tikzpicture}[every node/.style={circle, draw, inner sep=0pt, minimum size=12pt,text centered},node distance=0.7cm]
            \node (i) at (0,2) {i};
            \node (j) [below=of i] {j};
            \node (k) [left=of j] {k};
            \node (q) [above=of k] {l};
            \draw (q)--(i)--(k)--(j)--(q);
            \draw (i) to[bend left=15] (j);
            \draw (i) to[bend right=15] (j);
        \end{tikzpicture}
    }
    \hspace{0.8cm} % space between subfigures
    \subfloat[$k=l$]{%
        \begin{tikzpicture}[every node/.style={circle, draw, inner sep=0pt, minimum size=12pt,text centered},node distance=0.7cm]
            \node (i) at (0,2) {i};
            \node (j) [below=of i] {j};
            \node (k) [left=of j] {k};
            \node (q) [above=of k] {q};
            \draw (q)--(i)--(j)--(k)--(q);
            \draw (i) to[bend left=15] (k);
            \draw (i) to[bend right=15] (k);
        \end{tikzpicture}
    }
    \hspace{0.8cm} % space between subfigures
    \subfloat[$k=q$]{%
        \begin{tikzpicture}[every node/.style={circle, draw, inner sep=0pt, minimum size=12pt,text centered},node distance=0.7cm]
            \node (i) at (0,2) {i};
            \node (j) [below=of i] {j};
            \node (k) [left=of j] {k};
            \node (q) [above=of k] {l};
            \draw (q)--(i)--(j)--(k)--(q);
            \draw (i) to[bend left=15] (k);
            \draw (i) to[bend right=15] (k);
        \end{tikzpicture}
    }
    \caption{Reduced graphs of ($G$, Row 6 \Cref{tab:result_5_nodes-2}) by letting nodes of the edge in $E_H$ equal.}
    \label{fig:reduced_graph_13}
\end{figure}
For $F_{2+2+1}^{(5)}$, there are two disjoint subgraphs contained in $H$, given by $E_H^1=\{j-l,k-q\},E_H^2=\{j-q,k-l\}$, which lead to 
\begin{align}\label{eq:5-m-F221}
    F_{2+2+1}^{(5)} & = 2A_{ij}^2A_{jk}^2A_{ki}^2, \quad\sum_{i,j,k}F_{2+2+1}^{(5)} = 2\tr((A^{\circ2})^3).
\end{align}
Using the formula in \Cref{sec:summary} and combining \eqref{eq:5-m-F}-\eqref{eq:5-m-F221}, we have that 
\begin{align*}
    &\sum_{i\neq j\neq k\neq l\neq q}A_{ij}A_{ik}A_{jk}A_{il}A_{lq}A_{qi} = \sum_{i,j,k,l,q}F-\sum_{i,j,k,l}F_{2+1+1+1}^{(5)}+\sum_{i,j,k}F_{2+2+1}^{(5)}\\
    = & \tr((A^3)^{\circ2})-4\onebm^\top ((A^2)^{\circ2}\odot A^{\circ2})\onebm+2\tr((A^{\circ2})^3).
\end{align*}

\item Let $F=A_{ij}A_{jk}A_{jl}A_{jq}A_{kq}A_{lq}$. The corresponding graph of $F$ is given by Row 7 of \Cref{tab:result_5_nodes-2} and its complement is $H=\{V,E\}$ where $E=\{(i,k),(i,l),(l,k),(i,q)\}$. 

Thus there are three nonzero nontrivial contraction function $F_{2+1+1+1}^{(5)},F_{2+2+1}^{(5)}$ and $F_{3+1+1}^{(5)}$. For $F_{2+1+1+1}^{(5)}$, the unique 2-cliques exists in $H$ which are $(i,k),(i,l),(l,k),(i,q)$, thus, we obtain that
\begin{align}\label{eq:F_2_1_1_1}
    F_{2+1+1+1}^{(5)} & = F(i,j,i,l,q)+F(i,j,k,i,q)+F(i,j,k,l,i)+F(i,j,k,k,q) \nonumber\\
    &\sim F(i,j,i,l,k)+F(i,j,k,i,l)+F(i,j,k,l,i)+F(i,j,k,k,l) \nonumber\\
    & = 2A_{ij}^2A_{il}A_{jk}A_{jl}A_{kl}+A_{ij}^2A_{ik}A_{il}A_{jk}A_{jl}+A_{ij}A_{jk}^2A_{jl}A_{kl}^2, \nonumber\\
    \sum_{i,j,k,l}F_{2+1+1+1}^{(5)}(i,j,k,l) & = \sum_{i,j,k,l}(2A_{ij}^2A_{il}A_{jk}A_{jl}A_{kl}+A_{ij}^2A_{ik}A_{il}A_{jk}A_{jl}+A_{ij}A_{jk}^2A_{jl}A_{kl}^2) \nonumber\\
    & = 2\onebm^\top (A^{\circ2}A)\odot A^2\odot A\onebm+\onebm^\top (A^{\circ2}\odot A^2\odot A^2)\onebm \nonumber \\
    & +\onebm^\top A\diag ((A^{\circ2})^2 A).
\end{align}
To find $F_{2+2+1}^{(5)}$ and its summation, we let $q=i$ and $l=k$ and obtain that 
\begin{align}\label{eq:F_2_2_1}
    F_{2+2+1}^{(5)} & = F(i,j,k,k,i) = A_{ij}^2A_{ik}^2A_{jk}^2, \nonumber\\
    \sum_{i,j,k}F_{2+2+1}^{(5)}(i,j,k) & = \sum_{i,j,k}A_{ij}^2A_{jk}^2A_{ki^2} = \tr((A\odot A)^3).
\end{align}
By letting $l=k=i$, the last nonzero contraction function is given by
\begin{align}\label{eq:F_3_1_1}
    F_{3+1+1}^{(5)} & = F(i,j,i,i,q)\sim F(i,j,i,i,k) = A_{ij}^3A_{ik}^2A_{jk},\nonumber \\
    \sum_{i,j,k}F_{3+1+1}^{(5)}(i,j,k) & = \sum_{i,j,k}A_{ij}^3A_{ik}^2A_{jk} = \tr[(A\odot A\odot A)(A\odot A)A].
\end{align}
Using the formula in \Cref{sec:summary} and combing \eqref{eq:F_2_1_1_1}--\eqref{eq:F_3_1_1}, we have that 
\begin{align*}
    &\sum_{i\neq j\neq k\neq l\neq q}A_{ij}A_{jk}A_{jl}A_{jq}A_{kq}A_{lq}\\
    = & \sum_{i,j,k,l,q}F(i,j,k,l,q)-\sum_{i,j,k,l}F_{2+1+1+1}^{(5)}(i,j,k,l)+\sum_{i,j,k}F_{2+2+1}^{(5)}(i,j,k)+2\sum_{i,j,k}F_{3+1+1}^{(5)}(i,j,k) \\
    = &  \onebm^\top [A((A^2)^{\circ2}\odot A)-2(A^{\circ2}A)\odot A^2\odot A-A^{\circ2}\odot (A^2)^{\circ2}]\onebm-\onebm^\top A\diag ((A^{\circ2})^2 A)\\
    & +\tr[(A^{\circ2})^3+2A^{\circ3}A^{\circ2}A].
\end{align*}

\item Let $F=A_{ij}A_{jl}A_{jq}A_{kl}A_{kq}A_{lq}A_{qi}$. We first compute the general sum of $F$ as 
\begin{align}\label{eq:5-o-F}
    \sum_{i,j,k,l,q}A_{ij}A_{jl}A_{jq}A_{kl}A_{kq}A_{lq}A_{qi}& = \sum_{j,l,q}A_{jl}A_{lq}A_{qj}\sum_{i}A_{qi}A_{ij}\sum_{k}A_{lk}A_{kq}\nonumber\\
    &=\sum_{j,l,q}A_{jl}[A\odot A^2]_{lq}[A\odot A^2]_{qj}=\tr(A(A\odot A^2)^2).
\end{align}
Consider the graph $G = (V, E_G)$ with $E_G = \{q-i-j-q-l-k-q,l-j\}$ and its complement graph $H = (V, E_H)$ with $E_H = \{l-i-k-j\}$ (see Row 1 in \Cref{tab:result_5_nodes-3}). Using the same argument in Part (a), the contraction $F_{2+1+1+1}^{(5)}$ can be obtained from the reduced graphs of $G$, shown in \Cref{fig:reduced_graph_14}, as follows:
\begin{align*}
    F_{2+1+1+1}^{(5)} & = A_{ij}A_{jk}A_{ki}A_{jl}^2A_{kl}A_{li}+A_{ik}^2A_{jk}A_{jl}A_{kl}^2A_{li}+A_{ij}A_{jk}^2A_{kl}A_{jl}^2A_{li}.
\end{align*}
The sum of the first term above can be evaluated using the same argument in Part (f),~\Cref{prop:4-node}. For any fixed $l$, define $[T^l_1]_{ij}=A_{ij}A_{li}$ and $[T^l_2]_{ij}=A_{ij}A_{li}^2$ and we have that 
\begin{align*}
    \sum_l\sum_{i,j,k}A_{ij}A_{jk}A_{ki}A_{jl}^2A_{kl}A_{li}=\sum_l\sum_{i,j,k}(A_{ij}A_{li})(A_{jk}A_{lj}^2)(A_{ki}A_{lk})=\sum_l\tr((T_1^l)^2T_2^l).
\end{align*}

Let $\Ccal_{(3)}\in\Rbb^{n\times n\times n}$ be the tensor satisfying $[\Ccal_{(3)}]_{::k'}=A^{\circ2}$ for all $1\leq k'\leq n$. Then $\Tcal_1=\Acal_{(2)}\odot \Acal_{(3)}$ and $\Tcal_2=\Acal_{(2)}\odot \Ccal_{(3)}$ are the tensors that $[\Tcal_1]_{::l}=T_1^l$ and $[\Tcal_2]_{::l}=T_2^l$. Hence, it follows that 
\begin{align}\label{eq:5-o-clique-reduced}
    \sum_{i,j,k,l,q}A_{ij}A_{jk}A_{ki}A_{jl}^2A_{kl}A_{li}& = \langle (\Acal_{(2)}\odot\Acal_{(3)})\times_{1,2}(\Acal_{(2)}\odot\Acal_{(3)})\times_{1,2}(\Acal_{(2)}\odot\Ccal_{(3)}),\Ical_{n\times n\times n}\rangle.
\end{align}
From \eqref{eq:5-o-clique-reduced}, the sum of $F_{2+1+1+1}^{(5)}$ becomes
\begin{align}
    \sum_{i,j,k,l}F_{2+1+1+1}^{(5)} & = \langle (\Acal_{(2)}\odot\Acal_{(3)})\times_{1,2}(\Acal_{(2)}\odot\Acal_{(3)})\times_{1,2}(\Acal_{(2)}\odot\Ccal_{(3)}),\Ical_{n\times n\times n}\rangle \nonumber\\
    &+\tr(AA^{\circ2}(A^{\circ2}\odot A^2))+\tr(A^2((A^{\circ2}A)\odot A^{\circ2})).
\end{align}

\begin{figure}[!ht]
    \centering
    \subfloat[$i=k$]{%
        \begin{tikzpicture}[every node/.style={circle, draw, inner sep=0pt, minimum size=12pt,text centered},node distance=0.7cm]
            \node (j) at (0,2) {j};
            \node (k) [below=of j] {k};
            \node (l) [left=of k] {l};
            \node (q) [above=of l] {q};
            \draw (j)--(q)--(l)--(k)--(j)--(l);
            \draw (q) to[bend left=15] (k);
            \draw (q) to[bend right=15] (k);
        \end{tikzpicture}
    }
    \hspace{0.8cm} % space between subfigures
    \subfloat[$i=l$]{%
        \begin{tikzpicture}[every node/.style={circle, draw, inner sep=0pt, minimum size=12pt,text centered},node distance=0.7cm]
            \node (j) at (0,2) {j};
            \node (k) [below=of j] {k};
            \node (l) [left=of k] {l};
            \node (q) [above=of l] {q};
            \draw (l)--(k)--(q)--(j);
            \draw (q) to[bend left=15] (l);
            \draw (q) to[bend right=15] (l);
            \draw (j) to[bend left=15] (l);
            \draw (j) to[bend right=15] (l);            
        \end{tikzpicture}
    }
    \hspace{0.8cm} % space between subfigures
    \subfloat[$j=k$]{%
        \begin{tikzpicture}[every node/.style={circle, draw, inner sep=0pt, minimum size=12pt,text centered},node distance=0.7cm]
            \node (i) at (0,2) {i};
            \node (j) [below=of i] {j};
            \node (l) [left=of j] {l};
            \node (q) [above=of l] {q};
            \draw (l)--(q)--(i)--(j);
            \draw (q) to[bend left=15] (j);
            \draw (q) to[bend right=15] (j);
            \draw (l) to[bend left=15] (j);
            \draw (l) to[bend right=15] (j);
            
        \end{tikzpicture}
    }
    \caption{Reduced graphs of ($G$, Row 1 \Cref{tab:result_5_nodes-3}) by letting nodes of the edge in $E_H$ equal.}
    \label{fig:reduced_graph_14}
\end{figure}
For $F_{2+2+1}^{(5)}$, there is one disjoint subgraph contained in $H$, given by $E_H^1=\{i-l,j-k\}$, which lead to 
\begin{align}\label{eq:5-o-F221}
    F_{2+2+1}^{(5)} & = A_{ij}^3A_{jk}^2A_{ki}^2, \quad\sum_{i,j,k}F_{2+2+1}^{(5)} = \tr(A^{\circ3}(A^{\circ2})^2).
\end{align}
Using the formula in \Cref{sec:summary} and combining \eqref{eq:5-o-F}-\eqref{eq:5-o-F221}, we have that 
\begin{align*}
    &\sum_{i\neq j\neq k\neq l\neq q}A_{ij}A_{jl}A_{jq}A_{kl}A_{kq}A_{lq}A_{qi} = \sum_{i,j,k,l,q}F-\sum_{i,j,k,l}F_{2+1+1+1}^{(5)}+\sum_{i,j,k}F_{2+2+1}^{(5)}\\
    = & \tr\Big[A(A\odot A^2)^2-AA^{\circ2}(A^{\circ2}\odot A^2)-A^2((A^{\circ2}A)\odot A^{\circ2})+A^{\circ3}(A^{\circ2})^2\Big] \\
    & -\langle (\Acal_{(2)}\odot\Acal_{(3)})\times_{1,2}(\Acal_{(2)}\odot\Acal_{(3)})\times_{1,2}(\Acal_{(2)}\odot\Ccal_{(3)}),\Ical_{n\times n\times n}\rangle.
\end{align*}

\item Let $F=A_{ij}A_{jk}A_{jl}A_{ki}A_{kq}A_{lq}A_{li}A_{qi}$. We first compute the general sum of $F$. For any fixed $l,q$, $F$ can be grouped by $F=(A_{ij}A_{il}A_{iq}A_{lq})(A_{jk}A_{jl})(A_{ki}A_{kq})$. Define $[T_1^{l,q}]_{ij}=A_{ij}A_{il}A_{iq}A_{lq}, [T^{l,q}_2]_{ij}=A_{ij}A_{il},[T^{l,q}_3]_{ij}=A_{ij}A_{iq}$ then we rewrite it as $F=[T_1^{l,q}]_{ij}[T_2^{l,q}]_{jk}[T_3^{l,q}]_{ki}$. Hence, one has that 
\begin{align*}
    \sum_{i,j,k,l,q} F = \sum_{l,q}\sum_{i,j,k}[T_1^{l,q}]_{ij}[T_2^{l,q}]_{jk}[T_3^{l,q}]_{ki} = \sum_{l,q}\tr(T_1^{l,q}T_2^{l,q}T_3^{l,q}). 
\end{align*}

We further define tensors $\mathscr{A}_{(1,2)},\mathscr{A}_{(1,3)},..,\mathscr{A}_{(3,4)}\in\Rbb^{n\times n\times n\times n}$ where the subscripts specify the uncontracted modes. For instance, $\mathscr{A}_{(3,4)}$ satisfies $[\mathscr{A}_{(3,4)}]_{::k'l'}=A$ for any $k',l'\in[n]$. We set $\mathscr{T}_1=\mathscr{A}_{(3,4)}\odot \mathscr{A}_{(2,4)}\odot \mathscr{A}_{(2,3)}\odot \mathscr{A}_{(1,2)}$ so that $[\mathscr{T}_1]_{ijlq}=[\mathscr{A}_{(3,4)}]_{ijlq}[\mathscr{A}_{(2,4)}]_{ijlq}[\mathscr{A}_{(2,3)}]_{ijlq}[\mathscr{A}_{(1,2)}]_{ijlq}=A_{ij}A_{il}A_{iq}A_{lq}$. Similarly, we define $\mathscr{T}_2=\mathscr{A}_{(3,4)}\odot \mathscr{A}_{(2,4)}$ and $\mathscr{T}_3=\mathscr{A}_{(3,4)}\odot \mathscr{A}_{(2,3)}$ so that $[\mathscr{T}_2]_{ijlq}=[\mathscr{A}_{(3,4)}]_{ijlq}[\mathscr{A}_{(2,4)}]_{ijlq}=A_{ij}A_{il}$ and $[\mathscr{T}_3]_{ijlq}=[\mathscr{A}_{(3,4)}]_{ijlq}[\mathscr{A}_{(2,3)}]_{ijlq}=A_{ij}A_{iq}$. Hence,
\begin{align}\label{eq:5-p-F}
    &\sum_{i,j,k,l,q} A_{ij}A_{jk}A_{jl}A_{ki}A_{kq}A_{lq}A_{li}A_{qi} \nonumber\\
    &= \langle (\mathscr{A}_{(3,4)}\odot \mathscr{A}_{(2,4)}\odot \mathscr{A}_{(2,3)}\odot \mathscr{A}_{(1,2)})\times_{1,2}(\mathscr{A}_{(3,4)}\odot \mathscr{A}_{(2,4)})\times_{1,2}(\mathscr{A}_{(3,4)}\odot \mathscr{A}_{(2,3)}),\Ical_{n\times n\times n\times n}\rangle. 
\end{align}

Consider the graph $G = (V, E_G)$ and its complement graph $H = (V, E_H)$ with $E_H = \{j-q,k-l\}$ (see Row 2 in \Cref{tab:result_5_nodes-3}). Using the same arguments in Part (a), the contraction $F_{2+1+1+1}^{(5)}$ can be obtained from the reduced graphs of $G$, shown in \Cref{fig:reduced_graph_15}, as follows:
\begin{align*}
    F_{2+1+1+1}^{(5)} & = A_{ij}^2A_{jk}^2A_{jl}A_{ki}A_{kl}A_{li}+A_{ij}A_{jk}^2A_{kl}^2A_{ki}^2A_{li},
\end{align*}
which has the sum
\begin{align}
    \sum_{i,j,k,l}F_{2+1+1+1}^{(5)} & = 2\onebm^\top((AA^{\circ2})^{\circ2}\odot A^{\circ2})\onebm.
\end{align}
\begin{figure}[!ht]
    \centering
    \subfloat[$q=j$]{%
        \begin{tikzpicture}[every node/.style={circle, draw, inner sep=0pt, minimum size=12pt,text centered},node distance=0.7cm]
            \node (i) at (0,2) {i};
            \node (j) [below=of i] {j};
            \node (k) [left=of j] {k};
            \node (l) [above=of k] {l};
            \draw (k)--(i)--(l);
            % \draw (j)--(l);
            \draw (i) to[bend left=15] (j);
            \draw (i) to[bend right=15] (j);
            \draw (k) to[bend left=15] (j);
            \draw (k) to[bend right=15] (j);
            \draw (l) to[bend left=15] (j);
            \draw (l) to[bend right=15] (j);
            
        \end{tikzpicture}
    }
    \hspace{0.8cm} % space between subfigures
    \subfloat[$k=l$]{%
        \begin{tikzpicture}[every node/.style={circle, draw, inner sep=0pt, minimum size=12pt,text centered},node distance=0.7cm]
            \node (i) at (0,2) {i};
            \node (j) [below=of i] {j};
            \node (k) [left=of j] {k};
            \node (q) [above=of k] {q};
            \draw (q)--(i)--(j);
            \draw (i) to[bend left=15] (k);
            \draw (i) to[bend right=15] (k);
            \draw (j) to[bend left=15] (k);
            \draw (j) to[bend right=15] (k);         
            \draw (q) to[bend left=15] (k);
            \draw (q) to[bend right=15] (k);
        \end{tikzpicture}
    }
    \caption{Reduced graphs of ($G$, Row 2 \Cref{tab:result_5_nodes-3}) by letting nodes of the edge in $E_H$ equal.}
    \label{fig:reduced_graph_15}
\end{figure}

For $F_{2+2+1}^{(5)}$, there is one disjoint subgraph contained in $H$, given by $E_H^1=\{j-q,k-l\}$, which lead to 
\begin{align}\label{eq:5-p-F221}
    F_{2+2+1}^{(5)} & = A_{ij}^4A_{jk}^2A_{ki}^2, \quad\sum_{i,j,k}F_{2+2+1}^{(5)} = \tr(A^{\circ4}(A^{\circ2})^2).
\end{align}
Using the formula in \Cref{sec:summary} and combining \eqref{eq:5-p-F}-\eqref{eq:5-p-F221}, we have that 
\begin{align*}
    &\sum_{i\neq j\neq k\neq l\neq q}A_{ij}A_{jk}A_{jl}A_{ki}A_{kq}A_{lq}A_{li}A_{qi} = \sum_{i,j,k,l,q}F-\sum_{i,j,k,l}F_{2+1+1+1}^{(5)}+\sum_{i,j,k}F_{2+2+1}^{(5)}\\
    = & \langle (\mathscr{A}_{(3,4)}\odot \mathscr{A}_{(2,4)}\odot \mathscr{A}_{(2,3)}\odot \mathscr{A}_{(1,2)})\times_{1,2}(\mathscr{A}_{(3,4)}\odot \mathscr{A}_{(2,4)})\times_{1,2}(\mathscr{A}_{(3,4)}\odot \mathscr{A}_{(2,3)}),\Ical_{n\times n\times n\times n}\rangle\\
    & +\tr(A^{\circ4}(A^{\circ2})^2)-2\onebm^\top((AA^{\circ2})^{\circ2}\odot A^{\circ2})\onebm.
\end{align*}

\item Let $F=A_{ij}A_{jk}A_{jl}A_{jq}A_{kl}A_{kq}A_{lq}$.
We first compute the general sum of $F$ as 
\begin{align*}
    \sum_{i,j,k,l,q}A_{ij}A_{jk}A_{jl}A_{jq}A_{kl}A_{kq}A_{lq}& = \sum_{j,k,l,q}A_{jk}A_{jl}A_{jq}A_{kl}A_{kq}A_{lq}\sum_iA_{ij} \\
    & = \sum_{i,j,k,l}[A\onebm]_iA_{ij}A_{ik}A_{il}A_{jk}A_{jl}A_{kl},
\end{align*}
which can be evaluated using the same arguments in Part (f),~\Cref{prop:4-node}. For any fixed $l$, define $[T^l_1]_{ij}=[A\onebm]_iA_{ij}A_{li}$ and $[T^l_2]_{ij}=A_{ij}A_{li}$ and we have that 
\begin{align*}
    \sum_l\sum_{i,j,k}[A\onebm]_iA_{ij}A_{ik}A_{il}A_{jk}A_{jl}A_{kl}=\sum_l\sum_{i,j,k}([A\onebm]_iA_{ij}A_{li})(A_{jk}A_{lj})(A_{ki}A_{lk})=\sum_l\tr((T_2^l)^2T_1^l).
\end{align*}
Let $\Acal_{(2)}, \Acal_{(3)},\Dcal_{(2,3)}\in\Rbb^{n\times n\times n}$ be the tensors with $[\Acal_{(2)}]_{:k':}=[\Acal_{(3)}]_{::k'}=A$ and $[\Dcal_{(2,3)}]_{:k'l'}=A\onebm$ for all $1\leq k',l'\leq n$. For $\Tcal_1=\Acal_{(2)}\odot \Acal_{(3)}\odot\Dcal_{(2,3)}$ and $\Tcal_2=\Acal_{(2)}\odot \Acal_{(3)}$, one has that $[\Tcal_1]_{ijl}=[\Acal_{(2)}]_{ijl}[\Acal_{(3)}]_{ijl}[\Dcal_{(2,3)}]_{ijl}=[A\onebm]_iA_{il}A_{ij}$ and $[\Tcal_2]_{ijl}=A_{il}A_{ij}$. Hence, it follows that 
\begin{align}\label{eq:5-q-F}
    &\sum_{i,j,k,l,q}A_{ij}A_{jk}A_{jl}A_{jq}A_{kl}A_{kq}A_{lq}\nonumber\\
    & = \langle (\Acal_{(2)}\odot\Acal_{(3)})\times_{1,2}(\Acal_{(2)}\odot\Acal_{(3)})\times_{1,2}(\Acal_{(2)}\odot\Acal_{(3)}\odot\Dcal_{(2,3)}),\Ical_{n\times n\times n}\rangle.
\end{align}
Consider the graph $G = (V, E_G)$ and its complement graph $H = (V, E_H)$ with $E_H = \{q-i-l,k-i\}$ (see Row 3 in \Cref{tab:result_5_nodes-3}). Using the same arguments in Part (a), the contraction $F_{2+1+1+1}^{(5)}$ has three equivalent components and can be obtained from the reduced graphs of $G$ that
\begin{align*}
    F_{2+1+1+1}^{(5)} & = 3A_{ij}^2A_{jk}A_{jl}A_{jq}A_{kl}A_{kq}A_{lq}=3(A_{ij}^2A_{li})(A_{jk}A_{lj})(A_{ki}A_{lk}).
\end{align*}
Let $\Ccal_{(3)}\in\Rbb^{n\times n\times n}$ be the tensor satisfying $[\Ccal_{(3)}]_{::k'}=A^{\circ2}$ for all $1\leq k'\leq n$. Then it follows that $[\Acal_{(2)}\odot\Ccal_{(3)}]_{ijl}=[\Acal_{(2)}]_{ijl}[\Ccal_{(3)}]_{ijl}=A^2_{ij}A_{il}$ and therefore
\begin{align}\label{eq:5-q-F2111}
    \sum_{i,j,k,l}F_{2+1+1+1}^{(5)} & = 3\langle (\Acal_{(2)}\odot\Acal_{(3)})\times_{1,2}(\Acal_{(2)}\odot\Acal_{(3)})\times_{1,2}(\Acal_{(2)}\odot\Ccal_{(3)}),\Ical_{n\times n\times n}\rangle.
\end{align}
Using the formula in \Cref{sec:summary} and combining \eqref{eq:5-q-F}-\eqref{eq:5-q-F2111}, we have that 
\begin{align*}
    &\sum_{i\neq j\neq k\neq l\neq q}A_{ij}A_{jk}A_{jl}A_{jq}A_{kl}A_{kq}A_{lq} = \sum_{i,j,k,l,q}F-\sum_{i,j,k,l}F_{2+1+1+1}^{(5)}\\
    = & \langle (\Acal_{(2)}\odot\Acal_{(3)})\times_{1,2}(\Acal_{(2)}\odot\Acal_{(3)})\times_{1,2}(\Acal_{(2)}\odot(\Acal_{(3)}\odot\Dcal_{(2,3)}-3\Ccal_{(3)})),\Ical_{n\times n\times n}\rangle.
\end{align*}
\item Let $F = A_{ij}A_{jk}A_{kq}A_{qi}A_{qj}A_{jl}A_{lq}$ and $H=\{V,E\}$ with $E=\{(i,k),(i,l),(k,l)\}$. There are only two nonzero and nontrivial contractions, $F_{2+1+1+1}^{(5)}$ and $F_{3+1+1}^{(5)}$. For the first contraction function, there are three components and, by letting $i=k$, $k=l$ and $i=l$ in $F$, can be found by 
\begin{align}\label{eq:d_F_2_1_1_1}
    F_{2+1+1+1}^{(5)} & = F(i,j,i,l,q)+F(i,j,k,k,q)+F(i,j,k,i,q)\sim 2F(i,j,i,k,l)+F(i,j,k,k,l) \nonumber\\
    & = 2A_{ij}^2A_{il}^2A_{jk}A_{jl}A_{kl}+A_{ij}A_{il}A_{jk}^2A_{jl}A_{kl}^2, \nonumber\\
    \sum_{i,j,k,l}F_{2+1+1+1}^{(5)}&=2\sum_{j,l}\left(\sum_iA_{ji}^2A_{il}^2\right)A_{jl}\sum_kA_{jk}A_{kl}+\sum_{j,l}\left(\sum_iA_{ji}A_{il}\right)A_{jl}\sum_kA_{jk}^2A_{kl}^2 \nonumber\\
    & = 3\onebm^\top [A^2\odot A\odot (A^{\circ2})^2]\onebm.
\end{align}
The second contraction and its sum are provided by
\begin{align}\label{eq:d_F_3_1_1}
    F_{3+1+1}^{(5)} & = F(i,j,i,i,q) \sim F(i,j,i,i,k) = A_{ij}^3A_{jq}A_{qi}^3, \nonumber\\
    \sum_{i,j,k}F_{3+1+1}^{(5)} & = \sum_{i,j,k}A_{ij}^3A_{jq}A_{qi}^3=\tr(A^{\circ3}A A^{\circ3}) = \tr((A^{\circ3})^2A).
\end{align}
Using the formula in \Cref{sec:summary} and \eqref{eq:d_F_2_1_1_1}--\eqref{eq:d_F_3_1_1}, we have that 
\begin{align*}
    & \sum_{i\neq j\neq k\neq l\neq q}A_{ij}A_{jk}A_{kq}A_{qi}A_{qj}A_{jl}A_{lq} \\
    & = \sum_{i,j,k,l,q}F-\sum_{i,j,k,l}F_{2+1+1+1}^{(5)}+2\sum_{i,j,k,l}F_{3+1+1}^{(5)}\\
    & =  \sum_{j,q}A_{jq}\sum_iA_{ji}A_{iq}\sum_kA_{jk}A_{kq}\sum_lA_{jl}A_{lq}-3\onebm^\top [A^2\odot A\odot (A^{\circ2})^2]\onebm+2\tr((A^{\circ3})^2A) \\
    & = \onebm^{\top} [(A^2)^{\circ3}\odot A]\onebm-3\onebm^\top [A^2\odot A\odot (A^{\circ2})^2]\onebm+2\tr((A^{\circ3})^2A) \\
    & = \onebm^{\top} \Big[(A^2)^{\circ3}\odot A-3A^2\odot A\odot (A^{\circ2})^2\Big]\onebm+2\tr((A^{\circ3})^2A).
\end{align*}

\item Let $F=A_{ij}A_{jk}A_{jl}A_{jq}A_{kl}A_{kq}A_{lq}A_{qi}$. 
We first compute the general sum of $F$ as 
\begin{align*}
    \sum_{i,j,k,l,q}A_{ij}A_{jk}A_{jl}A_{jq}A_{kl}A_{kq}A_{lq}A_{qi}& = \sum_{j,k,l,q}A_{jk}A_{jl}A_{jq}A_{kl}A_{kq}A_{lq}\sum_iA_{qi}A_{ij} \\
    & = \sum_{i,j,k,l}[A\odot A^2]_{ij}A_{ik}A_{il}A_{jk}A_{jl}A_{kl},
\end{align*}
which can be evaluated using the same arguments in Part (f), ~\Cref{prop:4-node}. Define $[T^l_1]_{ij}=[A\odot A^2]_{ij}A_{li}$ and $[T^l_2]_{ij}=A_{ij}A_{li}$ and we have that 
\begin{align*}
    \sum_l\sum_{i,j,k}([A\odot A^2]_{ij}A_{li})(A_{jk}A_{lj})(A_{ki}A_{lk})=\sum_l\tr((T_2^l)^2T_1^l).
\end{align*}
Let $\Ecal_{(3)}\in\Rbb^{n\times n\times n}$ be the tensors with $[\Ecal_{(3)}]_{::k'}=A\odot A^2$. Then one has that $[\Acal_{(2)}]_{ijl}[\Ecal_{(3)}]_{ijl}=[A\odot A^2]_{ij}A_{il}$. Hence, it follows that 
\begin{align}\label{eq:5-s-F}
    &\sum_{i,j,k,l,q}A_{ij}A_{jk}A_{jl}A_{jq}A_{kl}A_{kq}A_{lq}A_{qi}\nonumber \\
    &= \langle (\Acal_{(2)}\odot\Acal_{(3)})\times_{1,2}(\Acal_{(2)}\odot\Acal_{(3)})\times_{1,2}(\Acal_{(2)}\odot\Ecal_{(3)},\Ical_{n\times n\times n}\rangle.
\end{align}

\begin{figure}[!ht]
    \centering
    \subfloat[$i=k$]{%
        \begin{tikzpicture}[every node/.style={circle, draw, inner sep=0pt, minimum size=12pt,text centered},node distance=0.7cm]
            \node (j) at (0,2) {j};
            \node (k) [below=of j] {k};
            \node (l) [left=of k] {l};
            \node (q) [above=of l] {q};
            \draw (k)--(l)--(q)--(j)--(l);
            % \draw (j)--(l);
            \draw (k) to[bend left=15] (j);
            \draw (k) to[bend right=15] (j);
            \draw (q) to[bend left=15] (k);
            \draw (q) to[bend right=15] (k);
        \end{tikzpicture}
    }
    \hspace{0.8cm} % space between subfigures
    \subfloat[$i=l$]{%
        \begin{tikzpicture}[every node/.style={circle, draw, inner sep=0pt, minimum size=12pt,text centered},node distance=0.7cm]
            \node (j) at (0,2) {j};
            \node (k) [below=of j] {k};
            \node (l) [left=of k] {l};
            \node (q) [above=of l] {q};
            \draw (l)--(k)--(q)--(j)--(k);
            % \draw (j)--(l);
            \draw (j) to[bend left=15] (l);
            \draw (j) to[bend right=15] (l);
            \draw (q) to[bend left=15] (l);
            \draw (q) to[bend right=15] (l);
        \end{tikzpicture}
    }
    \caption{Reduced graphs of ($G$, Row 5 \Cref{tab:result_5_nodes-3}) by letting nodes of the edge in $E_H$ equal.}
    \label{fig:reduced_graph_16}
\end{figure}

Consider the graph $G = (V, E_G)$ and its complement graph $H = (V, E_H)$ with $E_H = \{l-i-q\}$ (see Row 5 in \Cref{tab:result_5_nodes-3}). Using the same arguments in Part (a), the contraction $F_{2+1+1+1}^{(5)}$ has two equivalent components and can be obtained from the reduced graphs of $G$ (shown in \Cref{fig:reduced_graph_16}) that
\begin{align*}
    F_{2+1+1+1}^{(5)} & = 2A_{ij}^2A_{jk}^2A_{jl}A_{jq}A_{kl}A_{kq}A_{lq}=2(A_{ij}^2A_{li})(A^2_{jk}A_{lj})(A_{ki}A_{lk}).
\end{align*}
Let $\Ccal_{(3)}\in\Rbb^{n\times n\times n}$ be the tensor satisfying $[\Ccal_{(3)}]_{::k'}=A^{\circ2}$ for all $1\leq k'\leq n$. We have that 
\begin{align}\label{eq:5-s-F221}
    \sum_{i,j,k,l}F_{2+1+1+1}^{(5)} & = 2\langle (\Acal_{(2)}\odot\Acal_{(3)})\times_{1,2}(\Acal_{(2)}\odot\Ccal_{(3)})\times_{1,2}(\Acal_{(2)}\odot\Ccal_{(3)}),\Ical_{n\times n\times n}\rangle.
\end{align}

Using the formula in \Cref{sec:summary} and combining \eqref{eq:5-s-F}-\eqref{eq:5-s-F221}, we have that 
\begin{align*}
    &\sum_{i\neq j\neq k\neq l\neq q}A_{ij}A_{jk}A_{jl}A_{jq}A_{kl}A_{kq}A_{lq}A_{qi} = \sum_{i,j,k,l,q}F-\sum_{i,j,k,l}F_{2+1+1+1}^{(5)}\\
    = & \langle (\Acal_{(2)}\odot\Acal_{(3)})\times_{1,2}\Big[(\Acal_{(2)}\odot\Acal_{(3)})\times_{1,2}(\Acal_{(2)}\odot\Ecal_{(3)})-2(\Acal_{(2)}\odot\Ccal_{(3)})\times_{1,2}(\Acal_{(2)}\odot\Ccal_{(3)}) \Big],\Ical_{n\times n\times n}\rangle.
\end{align*}

\item Let $F=A_{ij}A_{ik}A_{il}A_{iq}A_{jk}A_{jl}A_{jq}A_{kq}A_{lq}$.
In order to compute the general sum of $F$, we rewrite it as $F=(A_{ij}A_{il}A_{iq}A_{lq})(A_{jk}A_{jl}A_{jq})(A_{ki}A_{kq})$ and each term can be expressed in a $4$-way tensor. Define $\mathscr{T}_1=\mathscr{A}_{(3,4)}\odot \mathscr{A}_{(2,4)}\odot \mathscr{A}_{(2,3)}\odot \mathscr{A}_{(1,2)}, \mathscr{T}_2=\mathscr{A}_{(3,4)}\odot \mathscr{A}_{(2,4)}\odot \mathscr{A}_{(2,3)}$ and $\mathscr{T}_3=\mathscr{A}_{(3,4)}\odot \mathscr{A}_{(2,3)}$ so that $[\mathscr{T}_1]_{ijlq}=A_{ij}A_{il}A_{iq}A_{lq},[\mathscr{T}_2]_{ijlq}=A_{ij}A_{il}A_{iq},[\mathscr{T}_3]_{ijlq}=A_{ij}A_{iq}$. Hence, it follows that
\begin{align}\label{eq:5-t-F}
    &\sum_{i,j,k,l,q}(A_{ij}A_{il}A_{iq}A_{lq})(A_{jk}A_{jl}A_{jq})(A_{ki}A_{kq})\nonumber\\
    &=\sum_{l,q}\sum_{i,j,k}[\mathscr{T}_1]_{ijlq}[\mathscr{T}_2]_{jklq}[\mathscr{T}_3]_{kilq} =\sum_{l,q}\tr([\mathscr{T}_1]_{::lq}[\mathscr{T}_2]_{::lq}[\mathscr{T}_3]_{::lq}) \nonumber\\
    & = \langle (\mathscr{A}_{(3,4)}\odot \mathscr{A}_{(2,4)}\odot \mathscr{A}_{(2,3)}\odot \mathscr{A}_{(1,2)})\times_{1,2}(\mathscr{A}_{(3,4)}\odot \mathscr{A}_{(2,4)}\odot \mathscr{A}_{(2,3)})\times_{1,2}(\mathscr{A}_{(3,4)}\odot \mathscr{A}_{(2,3)}),\Ical_{n\times n\times n\times n}\rangle. 
\end{align}
There is one edge $\{k-l\}$ in the complement graph. Let $k=l$ and we have that 
\begin{align*}
    F_{2+1+1+1}^{(5)} = A_{ij}A_{ik}A^2_{il}A_{jk}A_{jl}^2A_{kl}^2 = (A_{ij}A_{li}^2)(A_{jk}A_{lj}^2)(A_{ki}A_{lk}^2).
\end{align*}
Let $\Ccal_{(2)}\in\Rbb^{n\times n\times n}$ be the tensor satisfying $[\Ccal_{(2)}]_{:k':}=A^{\circ2}$ for all $1\leq k'\leq n$. Then $[\Ccal_{(2)}\odot\Acal_{(3)}]_{ijl}=A_{ij}A_{li}^2$ and one has that
\begin{align}\label{eq:5-t-F221}
    \sum_{i,j,k,l}F_{2+1+1+1}^{(5)} & = \sum_l \sum_{i,j,k}[\Ccal_{(2)}\odot\Acal_{(3)}]_{ijl}[\Ccal_{(2)}\odot\Acal_{(3)}]_{jkl}[\Ccal_{(2)}\odot\Acal_{(3)}]_{kil} = \sum_{l}\tr(([\Ccal_{(2)}\odot\Acal_{(3)}]_{::l})^3) \nonumber\\
    & = \langle (\Ccal_{(2)}\odot\Acal_{(3)})\times_{1,2}(\Ccal_{(2)}\odot\Acal_{(3)})\times_{1,2}(\Ccal_{(2)}\odot\Acal_{(3)}),\Ical_{n\times n\times n}\rangle.
\end{align}
Using the formula in \Cref{sec:summary} and combining \eqref{eq:5-t-F}-\eqref{eq:5-t-F221}, we have that 
\begin{align*}
    &\sum_{i\neq j\neq k\neq l\neq q}A_{ij}A_{ik}A_{il}A_{iq}A_{jk}A_{jl}A_{jq}A_{kq}A_{lq} = \sum_{i,j,k,l,q}F-\sum_{i,j,k,l}F_{2+1+1+1}^{(5)}\\
    = & \langle (\mathscr{A}_{(3,4)}\odot \mathscr{A}_{(2,4)}\odot \mathscr{A}_{(2,3)}\odot \mathscr{A}_{(1,2)})\times_{1,2}(\mathscr{A}_{(3,4)}\odot \mathscr{A}_{(2,4)}\odot \mathscr{A}_{(2,3)})\times_{1,2} \\
    &(\mathscr{A}_{(3,4)}\odot \mathscr{A}_{(2,3)}),\Ical_{n\times n\times n\times n}\rangle-\langle (\Ccal_{(2)}\odot\Acal_{(3)})\times_{1,2}(\Ccal_{(2)}\odot\Acal_{(3)})\times_{1,2}(\Ccal_{(2)}\odot\Acal_{(3)}),\Ical_{n\times n\times n}\rangle.
\end{align*}

\item Let $F=A_{ij}A_{ik}A_{il}A_{iq}A_{jk}A_{jl}A_{jq}A_{kl}A_{kq}A_{lq}$. Rewrite it as $F=(A_{ij}A_{il}A_{iq}A_{lq})(A_{jk}A_{jl}A_{jq})\times$ $(A_{ki}A_{kl}A_{kq})$ and each term can be expressed in a $4$-way tensor. Define $\mathscr{T}_1=\mathscr{A}_{(3,4)}\odot \mathscr{A}_{(2,4)}\odot \mathscr{A}_{(2,3)}\odot \mathscr{A}_{(1,2)}, \mathscr{T}_2=\mathscr{A}_{(3,4)}\odot \mathscr{A}_{(2,4)}\odot \mathscr{A}_{(2,3)}$ so that $[\mathscr{T}_1]_{ijlq}=A_{ij}A_{il}A_{iq}A_{lq},[\mathscr{T}_2]_{ijlq}=A_{ij}A_{il}A_{iq}$. Hence, it follows that
\begin{align*}
    &\sum_{i,j,k,l,q}(A_{ij}A_{il}A_{iq}A_{lq})(A_{jk}A_{jl}A_{jq})(A_{ki}A_{kl}A_{kq}) \\
    &=\sum_{l,q}\sum_{i,j,k}[\mathscr{T}_1]_{ijlq}[\mathscr{T}_2]_{jklq}[\mathscr{T}_2]_{kilq} =\sum_{l,q}\tr([\mathscr{T}_1]_{::lq}([\mathscr{T}_2]_{::lq})^2) \nonumber\\
    & = \langle (\mathscr{A}_{(3,4)}\odot \mathscr{A}_{(2,4)}\odot \mathscr{A}_{(2,3)}\odot \mathscr{A}_{(1,2)})\times_{1,2}(\mathscr{A}_{(3,4)}\odot \mathscr{A}_{(2,4)}\odot \mathscr{A}_{(2,3)}) \\
    &\times_{1,2}(\mathscr{A}_{(3,4)}\odot \mathscr{A}_{(2,4)}\odot \mathscr{A}_{(2,3)}),\Ical_{n\times n\times n\times n}\rangle,
\end{align*}
which is equal to $\sum_{i\neq j\neq k\neq l\neq q}A_{ij}A_{ik}A_{il}A_{iq}A_{jk}A_{jl}A_{jq}A_{kl}A_{kq}A_{lq}$ since there is no edge in the complement graph.

\end{enumerate}
\end{proof}

\subsection{Number of weighted subgraphs with 6 nodes}{\label{sec:node6}}

We give one more example for counting weighted subgraphs with 6 nodes. 

The following identities hold:
\begin{align*}
    \displaystyle\sum_{i\neq j\neq k\neq l\neq q\neq r}A_{ij}A_{jk}A_{kl}A_{lq}A_{qr}A_{ri}=
\end{align*}

\begin{proof}
    Let $F=A_{ij}A_{jk}A_{kl}A_{lq}A_{qr}A_{ri}$. 
We first compute the general sum of $F$ as 
\begin{align*}
    \sum_{i,j,k,l,q,r}A_{ij}A_{jk}A_{kl}A_{lq}A_{qr}A_{ri}& = \tr(A^6).
\end{align*}

\begin{figure}[!ht]
    \centering
    \subfloat[$r=j$]{%
        \begin{tikzpicture}[every node/.style={circle, draw, inner sep=2pt, minimum size=8pt}]
        \node (i) at (0,0) {i};
        \node (j) at ([xshift=0.974cm, yshift=-0.709cm] i) {j};
        \node (k) at ([xshift=-0.372cm, yshift=-0.927cm] j) {k};
        \node (l) at ([xshift=-1.204cm, yshift=0cm] k) {l};
        \node (q) at ([xshift=-0.372cm, yshift=0.927cm] l) {q};

        \draw (i) to[bend left=15] (j);
        \draw (i) to[bend right=15] (j);
        \draw (j)--(k)--(l)--(q)--(j);
        \end{tikzpicture}
    }
    \hspace{0.8cm} % space between subfigures
    \subfloat[$r=k$]{%
        \begin{tikzpicture}[every node/.style={circle, draw, inner sep=2pt, minimum size=8pt}]
        \node (i) at (0,0) {i};
        \node (j) at ([xshift=0.974cm, yshift=-0.709cm] i) {j};
        \node (k) at ([xshift=-0.372cm, yshift=-0.927cm] j) {k};
        \node (l) at ([xshift=-1.204cm, yshift=0cm] k) {l};
        \node (q) at ([xshift=-0.372cm, yshift=0.927cm] l) {q};

        \draw (k)--(j)--(i)--(k)--(q)--(l)--(k);
        \end{tikzpicture}
    }
    \caption{Reduced graphs of ($G$, Row 1 \Cref{tab:result_6_nodes}) by letting nodes of the edge in $E_H$ equal.}
    \label{fig:reduced_graph_6-graph-1}
\end{figure}

Consider the graph $G = (V, E_G)$ and its complement graph $H = (V, E_H)$ with $E_H = \{i-k-q,j-l-r\}$ (see Row 1 in \Cref{tab:result_6_nodes}). There are five nonzero contractions, $F_{3+3}^{(6)}$, $F_{3+1+1+1}^{(6)}$, $F_{3+2+1}^{(6)}$, $F_{2+2+2}^{(6)}$, $F_{2+2+1+1}^{(6)}$, $F_{2+1+1+1+1}^{(6)}$. Using the same arguments of Part (a) in the proof of \Cref{prop:5-node}, the reduced graphs of $G$ by forcing the nodes of the edges in $H$ equal are provided in \Cref{fig:reduced_graph_6-graph-1}. 
The contraction $F_{2+1+1+1+1}^{(6)}$ has two equivalent components and can be computed as
\begin{align*}
    F_{2+1+1+1+1}^{(6)} & = 6A_{ij}^2A_{jk}A_{kl}A_{lq}A_{qj}+3A_{ij}A_{jk}A_{ki}A_{kl}A_{lq}A_{qk}.
\end{align*}

For the contraction $F_{2+2+1+1}^{(6)}$, there are four equivalent components by letting $(q=j,r=k)$ or $(q=k,r=j)$ or $(q=j,r=l)$ or $(q=k,r=l)$, whose reduced graphs are provided in \Cref{fig:reduced_graph_6-graph-2} and we have that
\begin{align*}
    F_{2+2+1+1}^{(6)} & = 3A_{ij}A_{jk}^2A_{ki}A_{jl}A_{kl}+3A_{ij}^2A_{jk}^2A_{kl}^2+6A_{ij}A_{jk}A_{kl}A_{li}A_{lj}^2+6A_{ij}A_{jk}A_{kl}^3A_{li}
\end{align*}

\begin{figure}[!ht]
    \centering
    \subfloat[$q=j,r=k$]{%
        \begin{tikzpicture}[every node/.style={circle, draw, inner sep=0pt, minimum size=12pt,text centered},node distance=0.7cm]
            \node (i) at (0,2) {i};
            \node (j) [below=of i] {j};
            \node (k) [left=of j] {k};
            \node (l) [above=of k] {l};
            
            \draw (i)--(j)--(l)--(k)--(i);
            \draw (k) to[bend left=15] (j);
            \draw (k) to[bend right=15] (j);
        \end{tikzpicture}
    }
    \hspace{0.8cm} % space between subfigures
    \subfloat[$q=k,r=j$]{%
        \begin{tikzpicture}[every node/.style={circle, draw, inner sep=0pt, minimum size=12pt,text centered},node distance=0.7cm]
            \node (i) at (0,2) {i};
            \node (j) [below=of i] {j};
            \node (k) [left=of j] {k};
            \node (l) [above=of k] {l};
            
            \draw (i) to[bend left=15] (j);
            \draw (i) to[bend right=15] (j);
            \draw (k) to[bend left=15] (j);
            \draw (k) to[bend right=15] (j);
            \draw (l) to[bend left=15] (k);
            \draw (l) to[bend right=15] (k);
        \end{tikzpicture}
    }
    \hspace{0.8cm} % space between subfigures
    \subfloat[$r=l,q=j$]{%
        \begin{tikzpicture}[every node/.style={circle, draw, inner sep=0pt, minimum size=12pt,text centered},node distance=0.7cm]
            \node (i) at (0,2) {i};
            \node (j) [below=of i] {j};
            \node (k) [left=of j] {k};
            \node (l) [above=of k] {l};
            \draw (i)--(j)--(k)--(l)--(i);
            \draw (l) to[bend left=15] (j);
            \draw (l) to[bend right=15] (j);
        \end{tikzpicture}
    }
    \hspace{0.8cm} % space between subfigures
    \subfloat[$r=l,q=k$]{%
        \begin{tikzpicture}[every node/.style={circle, draw, inner sep=0pt, minimum size=12pt,text centered},node distance=0.7cm]
            \node (i) at (0,2) {i};
            \node (j) [below=of i] {j};
            \node (k) [left=of j] {k};
            \node (l) [above=of k] {l};
            \draw (i)--(j)--(k)--(l)--(i);
            \draw (l) to[bend left=15] (k);
            \draw (l) to[bend right=15] (k);
        \end{tikzpicture}
    }
    \caption{Reduced graphs of ($G$, Row 1 \Cref{tab:result_6_nodes}) by letting nodes of two edges in $E_H$ equal.}
    \label{fig:reduced_graph_6-graph-2}
\end{figure}

For the other two contractions, it is straightforward to obtain
\begin{align*}
    F_{2+2+2}^{(6)}=4A_{ij}^2A_{jk}^2A_{ki}^2,\quad  F_{3+1+1+1}^{(6)} = 2A_{ij}^2A_{ik}^2A_{il}^2, \quad F_{3+2+1}^{(6)} = 6A_{ij}^2A_{ik}^4,\quad F_{3+3}^{(6)} = A_{ij}^6.
\end{align*}
Here $F_{2+2+2}^{(6)}$ is obtained by setting $(i=l,q=j,r=k)$, $(i=l,r=j,q=k)$, $(q=j,i=k,r=l)$ and $(r=k,q=i,j=l)$.
Thus all the summations of the nonzero contractions are given by
\begin{align}
    \sum_{i,j,k,l,q}F_{2+1+1+1+1}^{(6)} & = 6\sum_{i,j,k}A_{ij}^2A_{jk}\sum_{l,q}A_{kl}A_{lq}A_{qj}+3\sum_k\sum_{i,j}A_{ij}A_{jk}A_{ki}\sum_{l,q}A_{kl}A_{lq}A_{qk} \nonumber\\
    & = 6\sum_{i,j,k}A_{ij}^2A_{jk}[A^3]_{jk}+3\sum_k([A^3]_{kk})^2=6\onebm^\top A^{\circ2}(A\circ A^3)\onebm+3\|\diag(A^3)\|^2_2 \label{eq:6-1}
\end{align}
\begin{align}
    \sum_{i,j,k,l}F_{2+2+1+1}^{(6)} =& 3\sum_{j,k}A_{jk}^2\sum_iA_{ij}A_{ik}\sum_{l}A_{jl}A_{lk}+3\sum_{i,j,k,l}A_{ij}^2A_{jk}^2A_{kl}^2 \nonumber \\
    & +6\sum_{j,l}A_{jl}^2\sum_iA_{ij}A_{li}\sum_{k}A_{jk}A_{kl}+6\sum_{k,l}A_{kl}^3\sum_{i,j}A_{ij}A_{jk}A_{li} \nonumber\\
    =& 9\onebm^\top(A^{\circ2}\circ A^2\circ A^2)\onebm+3\onebm^\top (A^{\circ2})^3\onebm+6\onebm^\top(A^{\circ3}\circ A^3)\onebm \label{eq:6-2}\\
    \sum_{i,j,k}F_{2+2+2}^{(6)}=&4\tr((A^{\circ2})^3), \quad \sum_{i,j}F_{3+3}^{(6)} = \onebm^\top A^{\circ6}\onebm \\
    \sum_{i,j,k,l}F_{3+1+1+1}^{(6)} = &  2\onebm^\top(A^{\circ2}\onebm)^{\circ3}, \quad \sum_{i,j,k}F_{3+2+1}^{(6)} = 6\onebm^\top A^{\circ2} A^{\circ4}\onebm. \label{eq:6-3}
\end{align}

Using the formula in \Cref{sec:summary} and combining \eqref{eq:6-1}-\eqref{eq:6-3}, we have that 
\begin{align*}
    &\sum_{i\neq j\neq k\neq l\neq q\neq r}A_{ij}A_{jk}A_{kl}A_{lq}A_{qr}A_{ri} \\
    = & \sum_{i,j,k,l,q,r}F-\sum_{i,j,k,l,q}F_{2+1+1+1+1}^{(6)}+\sum_{i,j,k,l}F_{2+2+1+1}^{(6)}+2\sum_{i,j,k,l}F_{3+1+1+1}^{(6)}\\
    & -\sum_{i,j,k}F_{2+2+2}^{(6)}- 2\sum_{i,j,k}F_{3+2+1}^{(6)}+4\sum_{i,j}F_{3+3}^{(6)}\\
    =&\tr(A^6) -(6\onebm^\top[A^{\circ2}(A\circ A^3)]\onebm+3\|\diag(A^3)\|^2_2)+4\onebm^\top(A^{\circ2}\onebm)^{\circ3}-4\tr((A^{\circ2})^3) \\
    & + \onebm^\top\Big(9A^{\circ2}\circ A^2\circ A^2+3(A^{\circ2})^3+6A^{\circ3}\circ A^3-12A^{\circ2}A^{{\circ4}}+4A^{\circ6} \Big)\onebm
\end{align*}
\end{proof}

\begin{table}[!ht]
\centering
\resizebox{\textwidth}{!}{
\begin{tabular}{>{\centering\arraybackslash}m{3cm} >{\centering\arraybackslash}m{3cm} c p{7cm}}
\toprule
\textbf{Graph} & \textbf{Complement} & \textbf{Sum}  & \textbf{Matrix Formula}\\
\midrule

% --- Row 1 ---
\adjustbox{valign=m}{
\begin{tikzpicture}[every node/.style={circle, draw, inner sep=2pt, minimum size=8pt}]
\def\r{1.0cm}
\node (i) at (90:\r)   {i};
\node (j) at (30:\r)   {j};
\node (k) at (330:\r)  {k};
\node (l) at (270:\r)  {l};
\node (q) at (210:\r)  {q};
\node (r) at (150:\r)  {r};
\draw (i)--(j)--(k)--(l)--(q)--(r)--(i);
\end{tikzpicture}

}
&
\adjustbox{valign=m}{
\begin{tikzpicture}[every node/.style={circle, draw, inner sep=2pt, minimum size=8pt}]
\def\r{1.0cm}
\node (i) at (90:\r)   {i};
\node (j) at (30:\r)   {j};
\node (k) at (330:\r)  {k};
\node (l) at (270:\r)  {l};
\node (q) at (210:\r)  {q};
\node (r) at (150:\r)  {r};
\draw (i)--(k)--(q)--(i);
\draw (j)--(l)--(r)--(j);
\draw (i)--(l);
\draw (r)--(k);
\draw (j)--(q);
\end{tikzpicture}
}
& $\displaystyle\sum_{i\neq j\neq k\neq l\neq q\neq r}A_{ij}A_{jk}A_{kl}A_{lq}A_{qr}A_{ri}$ & {$\tr(A^6) -\Big(6\onebm^\top[A^{\circ2}(A\circ A^3)]\onebm+3\|\diag(A^3)\|^2_2\Big)+4\onebm^\top(A^{\circ2}\onebm)^{\circ3}-4\tr((A^{\circ2})^3)+\onebm^\top\Big(9A^{\circ2}\circ A^2\circ A^2+3(A^{\circ2})^3+6A^{\circ3}\circ A^3-12A^{\circ2}A^{\circ4}+4A^{\circ6} \Big)\onebm$}\\

\bottomrule
\end{tabular}
}
\caption{One connected 6-node subgraphs and its counting formulas given in an efficient form.}
\label{tab:result_6_nodes}
\end{table}

\section{Numerical Verification}

We evaluate the proposed formulas for the 5-node subgraphs in \Cref{prop:5-node}, labeled $(A)$ through $(U)$, across matrices of varying sizes. For each experiment, we generate a symmetric matrix $A$ of size $n$ with zero diagonals, i.e., $\diag(A)=\mathbf{0}$. The upper-triangular entries $A_{ij}$ (for $i > j$) are independently sampled from the uniform distribution $U[0,1]$. Each configuration is simulated 100 times, and we report the average runtime.

\begin{figure}[!ht]
    \centering
    \includegraphics[width=0.45\linewidth]{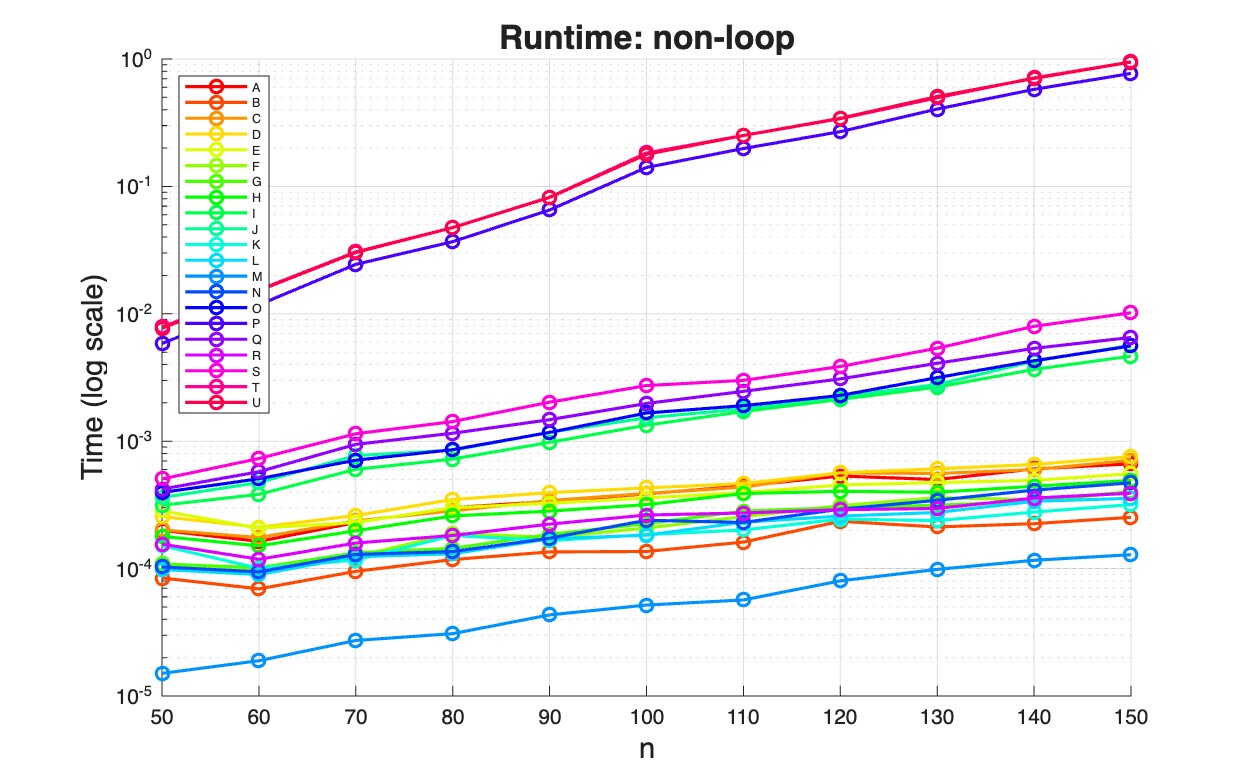}
    \includegraphics[width=0.45\linewidth]{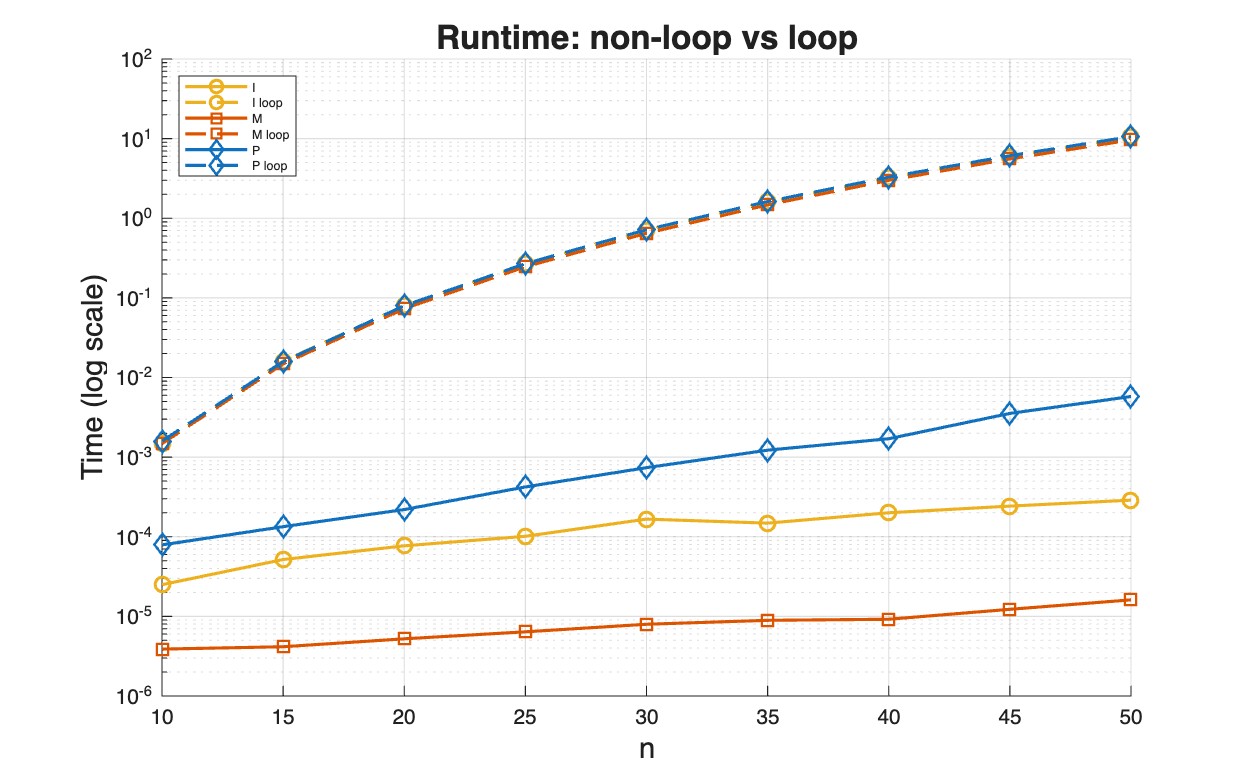}
    \caption{Runtime for the 5-node subgraph formulas (log-scale on the $y$-axis). Left: execution time of the proposed non-loop expressions. Right: comparison between the proposed formulas and their loop-based implementations for $(I),(M),(P)$.}

    \label{fig:runtime}
\end{figure}

The runtime of the corresponding \textsc{Matlab} implementations from \Cref{prop:5-node} is shown in the left subfigure of \Cref{fig:runtime}. Since the computations of $(P)$, $(T)$, and $(U)$ involve 4-way tensors, while those of $(I)$, $(J)$, $(O)$, $(Q)$, and $(S)$ involve 3-way tensors, the evaluation of $(P)$, $(T)$, and $(U)$ is slower than that of $(I)$, $(J)$, $(O)$, $(Q)$, and $(S)$. In turn, these 3-way tensor expressions require more time than subgraph formulas that rely solely on matrix operations. These trends are consistently confirmed in our simulations.

We further compare the runtime of the proposed non-loop formulas against their loop-based counterparts, also shown in the left subfigure of \Cref{fig:runtime}. We focus on subgraphs $(I)$, $(M)$, and $(P)$, which respectively depend on matrix, 3-way tensor, and 4-way tensor computations. For a matrix $A$ of size $n$, the proposed formula for $(P)$ is approximately $1{,}800\times$ faster than the loop version, while $(I)$ is about $36{,}000\times$ faster, and $(M)$ achieves an improvement of roughly $600{,}000\times$. The performance gap grows even further as $n$ increases.

\section{Discussion}\label{sec:discussion}
We introduce a simple and intuitive approach for counting subgraphs, based on the observation that the general block sum (GBS) $M_\sigma$ is often more computationally efficient to evaluate. However, the computation of $M_\sigma$ becomes increasingly demanding in both computation and memory as the graph includes larger cliques, particularly for $k\geq 6$. For instance, counting 4-clique requires handling 3-way tensors with a memory cost of $\Ocal(n^3)$, while counting 5-cliques involves 4-way tensors with a memory cost of $\Ocal(n^4)$, and so on. When memory usage becomes a significant bottleneck, implementing the algorithm in a lower-level language such as C++ is advantageous, as it provides finer control over memory management and computational efficiency.

% \section*{Acknowledgments}
% We would like to acknowledge the assistance of volunteers in putting
% together this example manuscript and supplement.

\bibliographystyle{plain}
\bibliography{ref}

@inproceedings{vassilevska2009finding,
  title={Finding, minimizing, and counting weighted subgraphs},
  author={Vassilevska, Virginia and Williams, Ryan},
  booktitle={Proceedings of the forty-first annual ACM symposium on Theory of computing},
  pages={455--464},
  year={2009}
}

@article{maugis2020testing,
  title={Testing for equivalence of network distribution using subgraph counts},
  author={Maugis, P.-A. G. and Olhede, Sofia C and Priebe, Carey E and Wolfe, Patrick J},
  journal={Journal of Computational and Graphical Statistics},
  volume={29},
  number={3},
  pages={455--465},
  year={2020},
  publisher={Taylor \& Francis}
}

@article{yin2018higher,
  title={Higher-order clustering in networks},
  author={Yin, Hao and Benson, Austin R and Leskovec, Jure},
  journal={Physical Review E},
  volume={97},
  number={5},
  pages={052306},
  year={2018},
  publisher={APS}
}

@article{gonen2011counting,
  title={Counting stars and other small subgraphs in sublinear-time},
  author={Gonen, Mira and Ron, Dana and Shavitt, Yuval},
  journal={SIAM Journal on Discrete Mathematics},
  volume={25},
  number={3},
  pages={1365--1411},
  year={2011},
  publisher={SIAM}
}

@article{ribeiro2021survey,
  title={A survey on subgraph counting: concepts, algorithms, and applications to network motifs and graphlets},
  author={Ribeiro, Pedro and Paredes, Pedro and Silva, Miguel EP and Aparicio, David and Silva, Fernando},
  journal={ACM computing surveys (csur)},
  volume={54},
  number={2},
  pages={1--36},
  year={2021},
  publisher={ACM New York, NY, USA}
}

@article{rota1964foundations,
  title={On the foundations of combinatorial theory: I. Theory of M{\"o}bius functions},
  author={Rota, Gian-Carlo},
  journal={Zeitschrift f{\"u}r Wahrscheinlichkeitstheorie und Verwandte Gebiete},
  volume={2},
  pages={340--368},
  year={1964}
}

@book{stanley2011enumerative,
  title={Enumerative Combinatorics, Volume 1},
  author={Stanley, Richard P},
  edition={2nd},
  series={Cambridge Studies in Advanced Mathematics},
  volume={49},
  publisher={Cambridge University Press},
  address={Cambridge},
  year={2011}
}

@inproceedings{pinar2017escape,
  title={Escape: Efficiently counting all 5-vertex subgraphs},
  author={Pinar, Ali and Seshadhri, Comandur and Vishal, Vaidyanathan},
  booktitle={Proceedings of the 26th international conference on world wide web},
  pages={1431--1440},
  year={2017}
}

@article{FFF15,
  title={Clique counting in mapreduce: Algorithms and experiments},
  author={Finocchi, Irene and Finocchi, Marco and Fusco, Emanuele G},
  journal={Journal of Experimental Algorithmics (JEA)},
  volume={20},
  pages={1--20},
  year={2015},
  publisher={ACM New York, NY, USA}
}

@article{AV14,
  title={Community detection in dense random networks},
  author={Arias-Castro, E. and Verzelen, N.},
  journal={Annals of Statistics},
  volume={42},
  number={3},
  pages={940--969},
  year={2014}
}

@InProceedings{GGC17,
  title = 	 {Two-Sample Tests for Large Random Graphs Using Network Statistics},
  author = 	 {Ghoshdastidar, Debarghya and Gutzeit, Maurilio and Carpentier, Alexandra and von Luxburg, Ulrike},
  booktitle = 	 {Proceedings of the 2017 Conference on Learning Theory},
  pages = 	 {954--977},
  year = 	 {2017},
  editor = 	 {Kale, Satyen and Shamir, Ohad},
  volume = 	 {65},
  series = 	 {Proceedings of Machine Learning Research},
  month = 	 {07--10 Jul},
  publisher =    {PMLR},
  url = 	 {https://proceedings.mlr.press/v65/ghoshdastidar17a.html}
}

@article{GL17,
  title={Testing for global network structure using small subgraph statistics},
  author={Gao, Chao and Lafferty, John},
  journal={arXiv preprint arXiv:1710.00862},
  year={2017}
}

@article{ABG18,
  title={Sublinear-time algorithms for counting star subgraphs via edge sampling},
  author={Aliakbarpour, Maryam and Biswas, Amartya Shankha and Gouleakis, Themis and Peebles, John and Rubinfeld, Ronitt and Yodpinyanee, Anak},
  journal={Algorithmica},
  volume={80},
  number={2},
  pages={668--697},
  year={2018},
  publisher={Springer}
}

@article{JKL21,
  title={Optimal adaptivity of signed-polygon statistics for network testing},
  author={Jin, Jiashun and Ke, Zheng Tracy and Luo, Shengming},
  journal={The Annals of Statistics},
  volume={49},
  number={6},
  pages={3408--3433},
  year={2021},
  publisher={Institute of Mathematical Statistics}
}

@article{YYS22,
  title={Hypothesis testing in sparse weighted stochastic block model},
  author={Yuan, Mingao and Yang, Fan and Shang, Zuofeng},
  journal={Statistical Papers},
  volume={63},
  number={4},
  pages={1051--1073},
  year={2022},
  publisher={Springer}
}

@article{YLFS22,
  title={Testing community structure for hypergraphs},
  author={Yuan, Mingao and Liu, Ruiqi and Feng, Yang and Shang, Zuofeng},
  journal={The Annals of Statistics},
  volume={50},
  number={1},
  pages={147--169},
  year={2022},
  publisher={Institute of Mathematical Statistics}
}

@article{YS22,
  title={Sharp detection boundaries on testing dense subhypergraph},
  author={Yuan, Mingao and Shang, Zuofeng},
  journal={Bernoulli},
  volume={28},
  number={4},
  pages={2459--2491},
  year={2022},
  publisher={Bernoulli Society for Mathematical Statistics and Probability}
}

@article{YR25,
  title={An uniformity index for random geometric graphs},
  author={Yuan, Mingao and Rahman, Irin},
  journal={Statistics},
  volume={59},
  number={4},
  pages={1063--1083},
  year={2025},
  publisher={Taylor \& Francis}
}

@article{YY25,
  title={Hypothesis testing for the dimension of random geometric graph},
  author={Yuan, Mingao and Yu, Feng},
  journal={arXiv preprint arXiv:2510.11844},
  year={2025}
}

@article{Y25,
  title={Hypothesis testing for the uniformity of random geometric graph},
  author={Yuan, Mingao},
  journal={arXiv preprint arXiv:2510.14210},
  year={2025}
}

@article{BR23,
  title={Number of cycles of small length in a graph},
  author={Barik, Sasmita and Reddy, Sane Umesh},
  journal={AKCE International Journal of Graphs and Combinatorics},
  volume={20},
  number={2},
  pages={134--147},
  year={2023}
}

@article{HM71,
  title={On the number of cycles in a graph},
  author={Harary, Frank and Manvel, Bennet},
  journal={Matematick{\'y} {\v c}asopis},
  volume={21},
  number={1},
  pages={55--63},
  year={1971}
}

@article{fan2022alma,
  title={{ALMA}: Alternating minimization algorithm for clustering mixture multilayer network},
  author={Fan, Xing and Pensky, Marianna and Yu, Feng and Zhang, Teng},
  journal={Journal of machine learning research},
  volume={23},
  number={330},
  pages={1--46},
  year={2022}
}

@article{N04,
  title={Analysis of weighted networks},
  author={Newman, M. E. J.},
  journal={Phys. Rev. E},
  volume={70},
  number={},
  pages={056131},
  year={2004}
}

\end{document}